# POINTWISE ERGODIC THEOREMS FOR ACTIONS OF GROUPS

AMOS NEVO

## CONTENTS

*Date*: April 5, 2005. Final version, to appear in : Handbook of Dynamical Systems, vol. 1B..

1991 *Mathematics Subject Classification*. Primary 22D40; Secondary 22E30, 28D10, 43A10, 43A90.

The author was supported in part by ISF Grant #126-01.

## 1. Introduction

In recent years a number of significant results and developments related to point-wise ergodic theorems for general measure-preserving actions of locally compact second countable (lcsc) groups have been established, including the solution of several long-standing open problems. The exposition that follows aims to survey some of these results and their proofs, and will include, in particular, an exposition of the following results.

(1) A complete solution to the ball averaging problem on Lie, and more generally lcsc, groups of polynomial volume growth. Namely, a proof that for any metric quasi-isometric to a word metric (and in particular, Riemmanian metrics on nilpotent Lie groups), the normalized ball averages satisfy the pointwise ergodic theorem in $L^1$. This brings to a very satisfactory close a long-standing problem in ergodic theory, dating at least to Calderon's 1952 paper on groups satisfying the doubling condition.

(2) In fact, two independent solutions will be described regarding the ball averaging problem in the case of connected Lie groups with polynomial volume growth, but both have the following in common. They resolve, in particular, a long-standing conjecture in the theory of amenable groups, dating at least to F. Greenleaf's 1969 book [Gr1]. The conjecture asserts that the sequence of powers of a neighbourhood on an amenable group constitute an asymptotically invariant sequence, namely has the Følner property. This conjecture was disproved for solvable groups with exponential growth, but has been now verified for groups with polynomial growth.

(3) The pointwise ergodic theorem in $L^1$ for a tempered sequence of asymptotically invariant sets was established recently, improving on the case of $L^2$ established earlier. This result resolves the long-standing problem of constructing *some* pointwise ergodic sequence in $L^1$ on an *arbitrary* amenable group. The ideas of the two available proofs will be briefly described.

(4) A new and streamlined account of the classical Dunford-Zygmund method will be described. This account allows the derivation of pointwise ergodic theorems for asymptotically invariant sequence on any lcsc amenable algebraic (or Lie) group over any local field, generalizing the Greenleaf-Emerson theorem. It also allows the construction of pointwise ergodic sequences on any lcsc algebraic (or Lie) group over any local field, generalizing Templeman's theorem for lcsc connected groups.

(5) A general spectral method will be described for the derivation of pointwise ergodic theorems for ball averages on Gelfand pairs. This method will be demonstrated for the ball averages on any lcsc simple algebraic group (over any local field). Pointwise theorems will be demonstrated also for the natural singular spherical averages on some of the Gelfand pairs.

(6) The proof of a pointwise ergodic theorem for actions of the free groups, generalizing Birkhoff's and Wiener's theorems for $\mathbb{Z}$ and $\mathbb{Z}^d$ will be described, using the general spectral method refered to above.

(7) The derivation of pointwise ergodic theorems for actions of simple algebraic groups with an explicit exponentially fast rate of convergence to the ergodic mean will be described. The same result will be described also for certain discrete lattice subgroups.

(8) Some ergodic theorems for semisimple Lie groups of real rank at least two will be described, which are in marked contrast to the results that Euclidean analogs might suggest, a contrast which has its roots in the exponential volume growth on semisimple groups.

Our goal is to elucidate some of the main ideas used in the proof of the pointwise ergodic theorems alluded to above. Our account of the pointwise ergodic theorems for groups with polynomial volume growth will be quite detailed, as these results are very recent and have not appeared before elsewhere. However, in the case of the spectral method, we have specifically attempted to give an account of the proofs which is as elementary as possible and demonstrated using the simplest available examples. The motivation for these choices is that the spectral methods which we employ require considerable background in the structure theory and representation theory of semisimple Lie groups, as well as classical singular integral theory. At this time these methods are not yet part of the standard tool kit in ergodic theory, and consequently it seems appropriate to give an exposition which focuses on the ergodic theorems and explains some of the main ideas in their proof, but requires as little as possible by way of background.

We have also tried to emphasized the pertinent open problems in the theory, many of which are presented along the way.

Ergodic theorems for actions of connected Lie groups, and particularly equidistribution theorems on homogeneous spaces and moduli spaces, have been developed and used in a rapidly expanding array of applications, many of which are presented in the two volumes of the present handbook. Thus it seems reasonable to limit the scope of our discussion in the present exposition and concentrate specifically on pointwise ergodic theorems, which have not been treated elsewhere.

We must note however that even whithin the more limited scope of pointwise ergodic theorems for general group actions our account has some important omissions. We mention some of these below, and offer as our rationale the fact that there already exist good expositions of these topics in the literature, some of which are referred to below. These omissions includes the analytic theory of homogeneous nilpotent Lie groups, and in particular the extensive theory of convolution operators, harmonic analysis, maximal functions and pointwise convergence theorems for diverse averages on Euclidean spaces, Heisenberg (-type) groups, homogeneous nilpotent groups and harmonic $AN$-groups (see [S2] and [CDKR] for an introduction to some of these topics). They also include the extensive results on equidistribution on homogeneous spaces (see [Da] and [Sta] for surveys, [GW] for some new results), as well as the general theory of actions of amenable locally compact second countable (lcsc) groups (see [OW], and [RW][L][We1] for more recent results). Another omission is the mean ergodic theorem for semisimple groups proved in [V], and other ergodic theorems on moduli spaces and their applications, which are described in detail in the present volume.

It is also natural to include in a discussion of maximal inequalities for group actions a discussion of convolution operators, particularly radial averages on general lcsc groups and their homogeneous spaces. This subject, for which the theory is very incomplete receives only very scant mention here, and we refer to [N7] for a short survey and many open problems.

## 2. Averaging along orbits in group actions

### 2.1. Averaging operators.
Let $G$ be a locally compact second countable (lcsc) group, $X$ an lcsc space on which $G$ acts (continuously) as a group of homemomorphisms, or more generally, a standard Borel space on which $G$ acts (measurably) by Borel automorphisms. Let $m$ be a $G$-invariant $\sigma$-finite Borel measure on $X$. The $G$-action on $X$ gives rise to a representation $\pi$ of $G$ as a group of isometries of the Banach spaces $L^p(X)$, given by $\pi(g)f(x) = f(g^{-1}x)$, and for $1 \le p < \infty$ the representation $\pi : G \to Iso(L^p(X))$ of $G$ into the isometry group is strongly continuous.

For any Borel probability measure on $G$, we can consider the averaging operator given, for every $f \in L^p(X)$, $1 \leq p \leq \infty$, by $\pi(\mu)f(x) = \int_G f(g^{-1}x)d\mu(g) = \int_G \pi(g)f(x)d\mu(g)$. The last equation is well-defined, and does indeed determine unambiguously an element of $L^p(X)$, and let us very briefly recall the well-known arguments proving this fact. Fix two Borel measurable functions $f_i$, $i = 1, 2$ on $X$, which have finite $L^p$-norm and are equal almost everywhere, and $h \in L^q(X)$, $q$ the dual exponent. Then by Fubini's theorem for $m$-almost all $x \in X$ the two functions $g \mapsto f_i(g^{-1}x)h(x)$ are $\mu$-integrable and equal, and so the values $\int_G f_i(g^{-1}x)d\mu(g) = \pi(\mu)f_i(x)$ are equal for $m$-almost all $x \in X$. Hence the latter integral, denoted by $\pi(\mu)f$, uniquely determines a function class (namely up to $m$-measure zero) for any $f \in L^p(X)$. Furthermore, $\pi(\mu)f$ has finite $L^p$-norm, and in fact by Hölder's inequality $\|\pi(\mu)f\|_p \leq \|f\|_p$. In addition, $\pi(\mu)f$ defined above coincides, as an element of $L^p(X)$, with the Lebesgue integral (w.r.t. $\mu$) of the Banach-space valued measurable function (strongly continuous if $1 \leq p < \infty$) given by $g \mapsto \pi(g)f$ from $G$ to $L^p(X)$. Detailed proofs of these well-known facts can be found in [DSI, Ch. III, §11, Thm 17, Ch. VIII, §7].

2.1.1. *The regular representation and the action by convolutions.* Consider the case where $X = G$, and the measure $m = m_G$ is left-invariant Haar measure. Then the operators $\lambda(g)f(x) = f(g^{-1}x)$ are isometric in every $L^p(G, m_G)$, and $g \mapsto \lambda(g)$ is the left regular representation. If $\mu$ is absolutely continuous with density $d\mu = b(g)dm_G(g)$, then the operator $\lambda(\mu)$ is the operator of left convolution by $\mu$ :

$$\lambda(\mu)f(x) = \int_G f(g^{-1}x)b(g)dm_G(g) = \mu * f(x)\,.$$

We can similarly consider the right regular representation of $G$, given by $\rho(g)f(x) = f(xg)$. Of course, in general $\rho(g)$ is *not* an isometric operator in $L^p(G, m_G)$, unless the left Haar measure $m_G$ is also right-invariant, namely unless $G$ is unimodular. Note also that the operator $\rho(\mu)$ is given by $\rho(\mu)f(x) = \int_G f(xg)b(g)dm_G(g)$, and is equal to the convolution operator $f * \mu^\vee$ if and only if $G$ is unimodular (here $\mu^\vee(A) = \mu(A^{-1})$). In particular if $G$ is unimodular and the measure $\mu$ is symmetric (namely satisfies $\mu^\vee = \mu$), then $\rho(\mu)f = f * \mu$.

2.2. **Ergodic theorems.** We will generally consider a family of probablity measures $\mu_t$ ($t \in \mathbb{R}_+$) on $G$, such that $t \mapsto \mu_t$ is $w^*$-continuous as a map into $M(G) = C_0(G)^*$. Usually each $\mu_t$ will have compact support, depending on $t$.

We will focus our attention on the case of probability-measure preserving actions, namely we assume that $(X, \mathcal{B}, m)$ is probability space, and $G$ is a group of measure-preserving transformations $g : X \to X$, and $m(gA) = m(A)$. Here the probability measures $\mu_t$ on $G$ can be regarded as family of averaging operators producing a sampling method along the group orbits $G \cdot x$ of $G$ in $X$, via $\pi(\mu_t)f(x) = \int_G f(g^{-1}x)d\mu_t(g)$.

Ergodicity of the $G$-action is defined as usual by the condition that every $G$-invariant integrable function is constant almost everywhere.

Our main goal below will be to establish a pointwise ergodic theorem in $L^p$ for interesting families of averages $\mu_t$ on $G$, in a general ergodic action $(X, m)$. By that we mean establishing the following convergence theorem :

$$\lim_{t \to \infty} \pi(\mu_t)f(x) = \int_X f \, dm$$

For $m$-almost every $x \in X$, and in the $L^p$-norm, for all $f \in L^p(X)$, where $1 \leq p < \infty$.

The ergodicity condition above is of course equivalent to the condition that the $\sigma$-algebra $\mathcal{I}$ of $G$-invariant sets is the trivial subalgebra of $\mathcal{B}$, consisting of sets of measure zero or one. In a general, not necessarily ergodic action we can consider the conditional expectation operator $\mathcal{E} : L^1(X, \mathcal{B}) \to L^1(X, \mathcal{I})$, and the sampling error along an orbit takes the form $|\pi(\mu_t)f(x) - \mathcal{E}f(x)|$. We recall that it is a well-known consequence of the standard ergodic decomposition theorem that in order to establish a pointwise ergodic theorem in an arbitrary action, it suffices to establish it for ergodic actions. In more detail, if for every $f \in L^p(X)$ $\lim_{t \to \infty} \pi(\mu_t)f(x) = \int_X f \, dm$ almost everywhere for every probability-preserving ergodic action, then $\lim_{t \to \infty} \pi(\mu_t)f(x) = \mathcal{E}f(x)$ almost everywhere for every probability-preserving action. We shall therefore often assume in what follows that the $G$-action is ergodic, when convenient.

In general, one would like to allow as many choices of sampling methods $\mu_t$ as possible, since different choices play absolutely crucial roles in different applications. A very basic distinction that arises is between the following cases :

(1) **$\mu_t$ is absolutely continuous** w.r.t. Haar measure on $G$. A fundamental question is to understand the case when $\mu_t = \beta_t$ are the normalized averages on a ball of radius $t$ and center $e$ w.r.t. an invariant metric on $G$. This includes for example an invariant Riemannian metric on a connected Lie group, or a word metric on an lcsc group.

(2) **$\mu_t$ is a singular mesure**, namely non-atomic and supported on a closed subset of $G$ of Haar measure zero. This includes, when $G$ is a Lie group, closed submanifolds of positive codimension in $G$, for example $\mu_t = \sigma_t$ the normalized averages on a sphere of $t$ and center $e$ defined using the restriction to the sphere of an invariant Riemannian metric.

(3) **$\mu_t$ is a discrete (atomic) measure**. This includes averages supported on discrete subgroups, for example when $\mu_t$ is supported on a lattice subgroup of $G$, or averages supported on integer points in a subvariety, when $G$ is a real algebraic group defined over $\mathbb{Q}$, for example.

We will begin by considering some absolutely continuous measures and some singular measures, and will comment on the important and interesting problem of discrete averages later, in §§10.5, 12.3 and Chapter 13. It is interesting to note that our discussion below will make it clear that the analysis of singular averages (for example, sphere averages) is a natural and indispensable ingredient in developing the theory of absolutely continuous averages, when the groups in question have exponential volume growth.

We note also that we will devote much of our attention in what follows to the study of averages defined geometrically via an invariant metric on the group, such as ball averages, shell averages and spherical averages. However, the spectral methods described in §9 -§12 below apply more generally and are not confined to radial averages. For semisimple groups, for example, we will consider in §12 also horospherical and many other non-radial averages.

2.3. **Maximal functions.** The maximal function associated with the family of averaging operators $\mu_t$ is defined by, for $f \in L^p(X)$:

$$M_\mu^* f(x) = f_\mu^*(x) = \sup_{t \in \mathbb{R}_+} |\pi(\mu_t) f(x)|$$

Let us hasten to note that in general it is not a-priori clear that the maximal function is well-defined and measurable for function classes in $L^p(X)$, and indeed, this is not always the case (see §2.3.1 below for further discussion).

We recall that a strong maximal inequality is an $L^p$-norm inequality for the maximal function, of the form $\left\| M_\mu^* f \right\|_p \le C_p \left\| f \right\|_p$, $\forall f \in L^p(X)$, where $1 < p \le \infty$. A weak-type maximal inequality is an estimate of the distribution function associated with the maximal function. Of particular interest is the weak-type $(1,1)$ maximal inequality given by :

$$m\{x \in X \,; \, f_\mu^*(x) > \delta\} \le \frac{C}{\delta} \left\| f \right\|_1 \quad, \;\; \forall f \in L^1(X)\,.$$

In the case of probability-measure-preserving actions it is natural to consider the maximal function given (for ergodic actions) by the following formula :

$$\tilde{f}_\mu^*(x) = \sup_{t \ge 0} \left| \pi(\mu_t) f(x) - \int_X f \, dm \right| \,.$$

The quantity $\int_X f \, dm$ is called the space average of the function $f$, and $\pi(\mu_t) f(x)$ are called the time averages along the orbit $G \cdot x$. Therefore, $\tilde{f}_\mu^*(x)$ is a bound for the largest error performed when sampling the values of $f$ along the $G$-orbit of $x \in X$, using $\mu_t$ as the sampling method. The strong $L^p$-Maximal inequality : $\left\| \tilde{f}_\mu^* \right\|_p \le C_p \left\| f \right\|_p$ bounds the total size of the error performed during the sampling process by the size of the function sampled, size being measured here by the $L^p$-norm. Similarly, the weak-type maximal inequality bounds the size of the set where the sampling error is larger than a fixed positive constant in terms of the size of the $L^1$-norm of $f$. In the case of general, not necessarily ergodic actions, the space average $\int_X f \, dm$ must be replaced by the conditional expectation $\mathcal{E} f$ of $f$ w.r.t. the $\sigma$-algebra of $G$-invariant functions. Of course, when $(X, m)$ is a probability space, the operators $\tilde{f}_\mu^*$ and $f_\mu^*$ satisfy exactly the same maximal inequalities and so we will continue with the analysis of just one of them, whose choice will be dictated by the problem at hand.

Finally, we recall that $L(\log^k L)(X)$ denotes the subspace of $L^1(X)$ consisting of functions for which $\int_X \left| f (\log^+ f)^k \right| dm(x)$ is finite.

2.3.1. *Measurability of the maximal function.* The maximal function is *not* necessarily a well-defined measurable function, even for very natural choices of the averaging operators $\mu_t$. Let us first note that if the averages $\mu_t$ are all absolutely continuous with respect to Haar measure on $G$, and $t \mapsto \mu_t$ is continuous w.r.t. the $L^1(G)$ norm, then for every $f \in L^p(X)$, and for $m$-almost every $x \in X$, $t \mapsto \pi(\mu_t) f(x)$ is a continuous function of $t$, see e.g. [DSI, Vol. I, Ch. VIII, §7, p.686]. It then follows that $f_\mu^*$ is indeed measurable, since the supremum may be taken of the countable set of rational numbers. However, absolute continuity is not a necessary condition for measurability of the maximal function, and we will discuss below many singular averages which give rise to a measurable maximal function.

Indeed absolute continuity is not necessary even for the continuity of $t \mapsto \pi(\mu_t)f(x)$, for almost every $x \in X$. Thus for example when $\mu_t = \sigma_t$ are the sphere averages on $\mathbb{R}^n$, $n \geq 3$ such a result for say bounded functions with bounded support, is established in [Co2, II.4], see also [S2, Ch. XI, §3.5 ].

Interestingly, for $\mathbb{R}^n$ (and other connected Lie groups) measurability of the maximal functions seems closely connected to considerations related to curvature, as is indicated by the following fact. Let $\partial Q_t$ denote the family of sets in the group $G = \mathbb{R}^2$ given by the boundaries of squares $Q_t$ centered at the origin, and let $q_t$ denote the usual (uniformly distributed, linear Lebesgue) probability measure on $\partial Q_t$. Let $\sigma_t$ be the rotation-invariant probability measure on the circle $S_t = \partial B_t$ of radius $t$ centered at the origin. Then the maximal function $f_q^*$ is *not* always measurable, even when $f$ is the characteristic function of a measurable set of finite measure, but $f_\sigma^*$ is measurable for such $f$ [Bn].

Clearly, when the action of $G$ on $X$ is a continuous action on a locally compact second countable (lcsc) space, the maximal function $f_\mu^*(x) = \sup_{t \in \mathbb{R}_+} |\pi(\mu_t)f(x)|$ is certainly measurable provided that $f \in C_c(X)$ and in addition $t \mapsto \int_G f(g^{-1}x)d\mu_t$ is continuous for every $f \in C_c(X)$. This condition will be satisfied by all the averages we will discuss below. A maximal inequality for $f_\mu^*$ when $f \in C_c(X)$ serves as an *a-priori* inequality, and is used to extend both the measurability of $f_\mu^*$ as well as the maximal inequality to the appropriate Lebesgue spaces. We refer to [N2, App. A] and [NS1, §2.2] for more information on these arguments.

### 2.4. A general recipe for proving pointwise ergodic theorems.
A proof of a pointwise ergodic theorem for a family of bounded operators $\pi(\mu_t) = T_t$ acting in $L^p(X, m)$ for some $1 \leq p < \infty$ can be obtained using the following four-step recipe.

(1) **Prove a mean ergodic theorem in $L^p(X)$** for the averages, namely show that $\lim_{t \to \infty} \left\| T_t f - \int_X f \, dm \right\|_p = 0$.

(2) **Find a dense subspace of functions $V \subset L^p(X)$ for which pointwise convergence holds**, namely $\lim_{t \to 0} T_t f = \int_X f \, dm$, for almost all $x \in X$, and for all $f \in V$.

(3) **Establish a strong $L^p$-maximal inequality** for the maximal function $f^*(x) = \sup_{t \geq 0} |T_t f(x)|$, where $f \in L^p(X)$, namely $\|f^*\|_p \leq C_p \|f\|_p$.

(4) **Use interpolation theory, either real or complex,** to establish a maximal inequality for the action of $T_t$ in $L^s(X)$, $s \neq p$.

We recall that the fact that (2) and (3) taken together imply pointwise convergence almost everywhere of $T_t f(x)$ for *every* $f \in L^p(X)$ is a formulation of the well-known Banach principle (see e.g. [Ga]). The identification of the limit is achieved in (1), by the mean ergodic theorem. Obviously, many variations are possible on this basic theme, for example the use of weak-type maximal inequalities, variational maximal inequalities, establishing the maximal inequality in (3) only on an a-priori dense subspace, as well as applying a wide array of interpolation methods.

Thus the basic problems we shall address below are establishing maximal inequalities for the family $\mu_t$, the corresponding pointwise convergence theorems for a dense subspace, the identification of the limit in the ergodic theorem, and applying interpolation techniques in the Lebesgue spaces $L^p(X)$.

## 3. Ergodic theorems for commutative groups

### 3.1. Flows of 1-parameter groups : Birkhoff's theorem.

In order to motivate the discussion below, let us start very concretely by considering one of the most basic example in ergodic theory. This will serve as our point of departure for several later developments.

**Example 3.1. Lines in $\mathbb{R}^2$ with an irrational slope**. Let $\ell = sw, s \in \mathbb{R}$, be a line in $\mathbb{R}^2$ with an irrational slope, and let $f$ be a function on the space $X = \mathbb{T}^2 = \mathbb{R}^2/\mathbb{Z}^2$ (i.e. $\mathbb{Z}^2$-periodic function on $\mathbb{R}^2$). Let $T_s = R_{sw} : X \to X$ denote the transformation given by translation on $X$, and let $m$ denote the translation-invariant probability measure on $X$. We note that $m$ is in fact the *unique* probability measure on $X$ invariant under the transformation $T_s$, $s \neq 0$ a fact which is immediate upon considering the Fourier coefficients of such an invariant measure. The time averages of $f$ along the orbit of $\ell$ passing through $x$ are defined by $\beta_t f(x) = \frac{1}{2t} \int_{-t}^{t} f(T_s x) ds$.

Recall the following classical results for the probability space $(X, m)$ and the transformations $T_s$ :

(1) **Weyl's Equidistribution Theorem.** $\forall f \in C(X)$, $\forall x \in X$ :

$$\lim_{t \to \infty} \beta_t f(x) = \int_X f(u) dm(u) = \text{ space average of } f .$$

(2) **Wiener's Differentiation Theorem**. $\forall f \in L^1(X)$ :

$$\lim_{t \to 0} \beta_t f(x) = f(x) , \text{ for almost all } x \in X .$$

(3) **Birkhoff's Pointwise Ergodic Theorem for flows.** $\forall f \in L^1(X)$ :

$$\lim_{t \to \infty} \beta_t f(x) = \int_X f dm , \text{ for almost all } x \in X .$$

(4) **Birkhoff's Pointwise ergodic theorem for invertible transformations.** The averages $\beta_n(\mathbb{Z}) f(x) = \frac{1}{2n+1} \sum_{k=-n}^{n} f(T_1^k(x))$, satisfy conclusion (1) and (3), namely $\frac{1}{2n+1} \sum_{k=-n}^{n} f(T_1^k x)$ have the same convergence properties as their continuous analogs.

Anticipating some later developments let us note that the classical ergodic averages $\beta_t$ are absolutely continuous measures on the line $\ell$, and in fact constitute the normalized averages on a ball of radius $t$ in $\mathbb{R}$, w.r.t. an invariant Riemannian metric (which is unique here up to scalar). $\beta_n(\mathbb{Z})$ are the normalized ball averages w.r.t. the induced invariant metric on the integer lattice $\mathbb{Z}$. This metric is here also a word metric on $\mathbb{Z}$ w.r.t. the set of generators $\{\pm 1\}$.

Let us further note the following :

(1) The equidistribution theorem stated in (1) holds in fact under the sole condition that $X$ be a compact metric space, and $T_s$ a homeomorphism (or a 1-parameter group of homeomorphisms) possesing a *unique* invariant probability measure. The case $X = \mathbb{T}^n$ was originally considered by H. Weyl.

(2) The differentiation theorem stated in (2) holds in fact for any standard measure space $(X, m)$, and any 1-parameter group of measure-preserving transformations. This result is due to N. Wiener [W, Thm. III'], and his proof combines Lebesgue's and Hardy-Littlewood's [HL] differentiation theorems on the real line with a principle of local transfer to a general measure-preserving flow (see the discussion in §5.4.2 for more details).

(3) The pointwise ergodic theorems stated in (3) and (4) hold in fact for any probability space $(X, m)$, and any invertible measure preserving transformation or 1-parameter group, under the sole condition of ergodicity. This is the content of Birkhoff's theorem [B].

Actions of the real line with an invariant probability measure exist in great abundance. We mention briefly the following examples, focusing mostly on those that will reappear again later.

**Example 3.2.**    (1) Any complete, divergence free vector field on a Riemannian manifold generates an associated volume-preserving flow. The total Riemannian volume if of course finite if the manifold is compact.

(2) In particular, the geodesic flow on a compact Riemannian manifold gives an $\mathbb{R}$-action preserving a finite volume. This includes the geodesic flow on a compact (or more generally, finite volume) locally symmetric space.

(3) The horocycle flow on a compact (or finite volume) surface of constant negative curvature also gives a volume-preserving flow. Similarly, analogous flows can be defined for all locally-symmetric spaces associated with semisimple Lie groups, by considering actions of one-parameter unipotent subgroups.

(4) Any 1-parameter subgroup of the connected Lie group $G$ acting by translations on $G/\Gamma$, where $\Gamma \subset$ is a discrete lattice subgroup of $G$.

Actions of $\mathbb{Z}$ with invariant measure are just as abundant. Indeed, both $\mathbb{Z}$ and $\mathbb{R}$ have the all-important property that each of their actions by homemorphisms of a compact metric space posseses at least one invariant ergodic measure, and this fact gives of course rise to a vast collection of examples. We mention here however only one, which will be particularly important in what follows. Let $G$ be any lcsc group, and $\Gamma \subset G$ a lattice. Then any element $g \in G$, as well as any subgroup $H \subset G$ act by measure-preserving transformations on the probability space $(G/\Gamma, m)$. This incudes the case where $\Gamma$ is a lattice in a simple non-compact algebraic group over a locally compact non-discrete field. In the latter case, the Howe-Moore mixing theorem asserts the remarkable fact that the action of every element $g$ is not only ergodic, but in fact mixing, provided only that the powers of $g$ are not confined to a compact subgroup, and that $G$ has no compact factors. We recall that mixing means that the correlations $\langle T^n f, f' \rangle$ converge to zero, if $f$ or $f'$ has zero integral.

We remark that equidistribution for *every* orbit does not hold in many of the examples above, e.g. for the case of geodesic flows on compact surfaces of constant negative curvature. Convergence for *every* orbit fails even if the function is assumed continuous, or even smooth. The restriction to *almost* every starting point is thus essential in the pointwise ergodic theorem.

3.2. **Flows of commutative multi-parameter groups : Wiener's theorem.**
Naturally, the next problem to consider is the generalization of the pointwise ergodic theorem from the case of a one-parameter flow to that of several commuting flows, and in particular from ball averages in actions of $\mathbb{R}$ (or $\mathbb{Z}$) to ball averages in actions of $\mathbb{R}^d$ (or $\mathbb{Z}^d$). This problem was solved by N. Wiener in [W], where he introduced several key ideas that came to play an important role in the further development of ergodic theory for groups with polynomial volume growth, and more generally amenable groups. We will focus below on Wiener's covering Lemma [W, Lemma C'], and the introduction of a transfer principle [W, Proof of Thm IV']. These arguments are geometric by nature and in particular, as already noted by Wiener, they apply equally well to ball averages on $\mathbb{R}^d$ or the lattice $\mathbb{Z}^d$, as well as to many other averages, e.g on Euclidean cells centered at the origin (see also Pitt's discussion [P]). Even further, they were generalized by A. Calderon [C1] to apply to certain non-commutative groups with polynomial volume growth. In §5 we will present these arguments and generalizations thereof, but before delving into the proofs, it may be instructive to consider some examples to which Wiener's (or Pitt's) theorem applies.

**Example 3.3.**    (1) Let $\mathcal{L}_n$ denote the space of unimodular lattices in $\mathbb{R}^n$. Clearly $SL_n(\mathbb{R})$ acts on $\mathcal{L}_n$, and the action is easily seen to be transitive. The stability group of the lattice $\mathbb{Z}^n$ is $SL_n(\mathbb{Z})$, so that we can identify $\mathcal{L}_n$ with the homogeneous space $SL_n(\mathbb{R})/SL_n(\mathbb{Z})$. The latter space carries a finite measure $m$ which is $SL_n(\mathbb{R})$-invariant, and we normalize it to be a probability measure. We fix a given non-trivial representation of $SL_3(\mathbb{R})$ in $SL_n(\mathbb{R})$, and let $A \cong \{a_t b_s | (t,s) \in \mathbb{R}^2\}$ denote the two dimensional vector group of diagonal matrices in $SL_3(\mathbb{R}) \subset SL_n(\mathbb{R})$, $n \geq 3$. Then for all $f \in L^1(\mathcal{L}_n)$ and for almost every $x \in \mathcal{L}_n$

$$\lim_{r \to \infty} \frac{1}{\pi r^2} \int_{t^2 + s^2 \leq r^2} f(a_t b_s x) dt ds = \int_{\mathcal{L}_n} f(x) dm$$

(2) Let $A_1, \ldots, A_d$ be $d$ commuting $n \times n$ matrices with integer entries and determinant $\{\pm 1\}$, and consider their natural action on $\mathbb{T}^n = \mathbb{R}^n/\mathbb{Z}^n$ by group automorphisms. Then, for every $f \in L^1(\mathbb{T}^n)$ and almost every $x \in \mathbb{T}^n$

$$\lim_{r \to \infty} \frac{1}{(2r+1)^d} \sum_{-r \leq i_j \leq r} f(A_1^{i_1} \cdots A_d^{i_d} x) = \int_{\mathbb{T}^n} f(x) dm$$

provided every integrable function invariant under $A_1, \ldots, A_d$ is constant. This is the case if at least one of the $A_i$ has no roots of unity as eigenvalues, for example.

(3) Consider the compact Abelian group $K = (\mathbb{Z}/2\mathbb{Z})^{\mathbb{Z}^2}$, with Haar measure given by the natural product measure. The $\mathbb{Z}^2$-action by coordinate shifts is an ergodic measure-preserving action by automorphisms of the compact Abelian group $K$. Denoting by $a_1$ and $b_1$ unit translations in the direction of the axes, we have : $\forall f \in L^1(K)$ and for almost every $x \in K$

$$\lim_{n \to \infty} \frac{1}{|B_n \cap \mathbb{Z}^2|} \sum_{(k,m) \in \mathbb{Z}^2 ; k^2 + m^2 \leq n^2} f(a_k b_m x) = \int_K f(x) dm$$

The example above can of course be greatly generalized, by replacing $\mathbb{Z}^2$ be other infinite Abelian groups, replacing $\mathbb{Z}_2$ by other compact groups, and also considering the extensive collection of closed shift-invariant subgroups arising in such actions. This gives rise to a wealth of ergodic measure preserving actions of multi-parameter groups by automorphisms of compact Abelian groups. For a discussion of this topic and its connections to commutative algebra and number theory we refer to [Sch].

In §5 we will considerably expand the scope of the discussion, formulate the ball averaging problem in ergodic theory and describe the complete solution of the ball averaging problem for all lcsc groups with polynomial volume growth. This result constitutes a common generalization of the pointwise and maximal ergodic theorems of Birkhoff, Wiener and Calderon, and brings to a satisfactory close a long line of development in ergodic theory. The result depends on some very recent developments and has not appeared before elsewhere. The proof is based on the arguments of Wiener and Calderon, together with one further ingredient, namely the asymptotic invariance (under translation) of the balls in these groups. This fact was recently established for connected Lie groups of polynomial volume growth by two interesting independent arguments which we will describe. In the general case of lcsc group asymptotic invariance of the balls ultimately depends also on considerations introduced by M. Gromov in his celebrated theorem [G] asserting that groups with polynomial volume growth are virtually nilpotent.

We will have on a number of occasions to use results on the (polynomial or exponential) volume growth of balls in most of the groups that will come under consideration below, and so we therefore now turn to this issue.

## 4. Invariant metrics, volume growth, and ball averages

4.1. **Growth type of groups.** Let $G$ be a compactly generated locally compact second countable (lcsc) group. If $V$ is a compact set generating $G$, namely $\cup_{n \in \mathbb{N}} V^n = G$, we can define a distance function $|g|_V$ on $G$ via $|g|_V = \min\{n \, ; \, g \in V^n\}$, where $V^0 = \{e\}$ and $V^n = V \cdot V \cdots V$ is the set of $n$-fold products of elements of $V$. The distance function is inversion-invariant if and only if the set $V$ is symmetric, namely $V = V^{-1}$, a condition equivalent to the symmetry of the function $d_V(g, h) = |g^{-1}h|_V$. We will assume $V = V^{-1}$ from now on, and then the associated function $d_V$ is a left-$G$-invariant metric on $G$, which we will call the word metric determined by $V$. Since $G$ is lcsc, for some $n \in \mathbb{N}$ $V^n$ contains an open set and then $m_G(V^{n+k}) > 0$ for all $k \in \mathbb{N}$, where $m_G$ is (any) Haar measure. Clearly $V^n V^k = V^{n+k}$, and thus the sequence $\log m_G(V^{n+k})$ is subadditive, and it follows that the limit $\lim_{n \to \infty} \frac{1}{n} \log m_G(V^n) = h_V$ exists. It is straightforward to see that if $h_V > 0$ for some $V$, then $h_{V'} > 0$ for any other compact generating set $V'$ (and any Haar measure). If $h_V > 0$ $G$ is called a group of exponential volume growth, and if $h_V = 0$ for some, or equivalently, all compact generating sets $V$, $G$ is called a group of subexponential growth (and it is then necessarily unimodular). In that case we can consider the quantity $\limsup_{n \to \infty} \frac{\log m_G(V^n)}{\log n} = q_V$, and we recall the well known fact that if $q_V < \infty$ for one compact generating set $V$, then $q_{V'} < \infty$ for any other compact generating set $V'$. In fact, in this case there exists a unique $0 < q(G) < \infty$ depending only on $G$, such that the ratio $m_G(V^n)/n^q$ is bounded for $n \in \mathbb{N}$, for any $V$ as above. This follows immediately from the fact that for large enough $k$, a ball of radius $k$ of one metric contains a ball of radius $ck$ of the second metric, and is contained is a ball of radius $Ck$ of the second metric, where $c \leq C$ are fixed positive constants. If $q_V < \infty$ for some compact generating set then $G$ is called a group of polynomial volume growth, and if $h_V = 0$ and $q_V = \infty$ for some compact generating set $V$, then $G$ is called a group of intermediate growth. As is well known, this possibility does in fact arise (see [dlH] for an accessible exposition of such a group, constructed by R. Grigorchuk). However, such groups do not arise as subgroups of connected Lie groups or algebraic groups.

We have used the metrics of the form $d_V$ described above to define the growth type of a group $G$. There are of course many other left-invariant metrics on a group, but since invariant metrics can be easily rescaled, not all of them will give the same growth type. For example, taking the metric $d(t,s) = \log(1 + |t - s|)$ on $G = \mathbb{R}$, we obtain an invariant metric whose balls centered at $0$ (denoted $B_t$) satisfy $m_{\mathbb{R}}(B_t) = \exp(t) - 1$, namely have exponential growth. It therefore natural to define a left-invariant metric $d$ on $G$ to be admissible, if the quantities $h_d = \lim_{t \to \infty} \frac{1}{t} \log m_G(B_t)$ and $q_d = \limsup_{t \to \infty} \frac{\log m_G(B_t)}{\log t}$ exhibit the same behavior as the quantities $h_V$ and $q_V$ defined by the metrics $d_V$ associated with compact generating sets $V$. In other words, $h_d > 0$ iff $h_V > 0$, and $q_d < \infty$ iff $q_V < \infty$. In the sequel, we will consider only admissible metrics on $G$.

4.2. **Invariant metrics.** In general, we can let $G$ be an lcsc group, and let $N : G \to \mathbb{R}_+$ be a continuous proper function satisfying $N(gh) \leq C(N(g) + N(h))$, namely a quasi-norm on $G$. Then $d_N(g, h) = N(g^{-1}h)$ defines a left-invariant quasi-metric on $G$. We can consider the sets $B_t^N = \{g \in G \,;\, N(g) \leq t\}$ and the normalized measures $\beta_t^N$ with density $\chi_{B_t^N}/m_G(B_t^N)$. When $G$ is say a connected Lie group and $N$ is sufficiently regular, it is also possible to define for every positive radius $t$ the natural probability measure $\sigma_t^N$ supported on $\partial B_t^N = S_t^N = \{g \in G \,;\, N(g) = t\}$. Our discussion below can be in principle extended to this more general context, but in the interest of simplicity we will restrict ourselves to a discussion of functions of the form $N(g) = d(e, g)$ where $d$ is an invariant admissible metric (so that the constant $C$ in the definition of a quasi-norm is equal to 1). We remark that one interesting natural set of examples consist of homogeneous invariant quasi-norms on a homogeneous connected nilpotent Lie groups. Such quasi-norms give rise to a quasi-distance which is not necessarily a metric, but in fact an equivalent homogeneous invariant norm can always be found. We refer e.g. to [S2](Ch. XIII, §5 and §7B) and [HS] for more on this topic.

Another important general example of a left-invariant metric is the function $d(g, h) = \log(1 + \|\tau(g)^{-1}\tau(h)\|)$, where $\|\cdot\|$ is a symmetric operator (or even just linear) norm on $M_n(\mathbb{R})$, and $\tau : G \to GL_n(\mathbb{R})$ is a faithful linear representation (and we assume that $d$ is admissible).

It is natural to introduce certain properties of metrics that facilitate the discussion and are relevant to questions of growth and ball averaging. Let us start by defining two notions of equivalence between metrics.

**Definition 4.1.** Consider two metric space $(X, d)$ and $(X', d')$ and a function $f : X \to X'$.

(1) $Y \subset X$ is called $C$-dense in $X$ if there exists $C \geq 0$ such that every $x \in X$ is at distance at most $C$ from $Y$.
(2) $f$ is called a **Quasi-isometry**, if $f$ satisfies, for some $B, b > 0$ and all $x, y \in X$

$$\frac{1}{B}d'(f(x), f(y)) - b \leq d(x, y) \leq Bd'(f(x), f(y)) + b$$

and in addition $f(X)$ is $C$-dense.
(3) $f$ is called a **Coarse-isometry**, if $f$ satisfies, for some $b \geq 0$ and all $x, y \in X$

$$d'(f(x), f(y)) - b \leq d(x, y) \leq d'(f(x), f(y)) + b$$

and in addition $f(X)$ is $C$-dense.

Next, consider the following possibilities regarding approximate analogs to geodesics in a metric space.

**Definition 4.2.** Let $(X, d)$ be a metric space. The metric $d$ is called

(1) **Discretely coarsely geodesic** [AM] if there exists $C \geq 0$ such that for any $x, y \in X$ it is possible to find a finite sequence of points $x = x_0, x_1 \ldots, x_n = y$ satisfying $d(x_{i-1}, x_i) \leq C$, $1 \leq i \leq n$, and $d(x, y) = \sum_{i=1}^{n} d(x_{i-1}, x_i)$.

(2) **Asymptotically geodesic** [Br] if for every $\varepsilon > 0$ there exists $C_\varepsilon > 0$, such that given $x, y \in X$, it is possible to find a finite sequence of points $x = x_0, x_1 \ldots, x_n = y$ satisfying $d(x_{i-1}, x_i) \leq C_\varepsilon$, $1 \leq i \leq n$, and $d(x, y) \geq (1 - \varepsilon) \sum_{i=1}^{n} d(x_{i-1}, x_i)$.

(3) **Monotone** [Te] if there exists $C \geq 1$ such that for any $x, y \in X$ it is possible to find a finite sequence of points $x = x_0, x_1 \ldots, x_n = y$ satisfying $d(x_{i-1}, x_i) \leq C$, and $d(x, x_{i-1}) + 1 \leq d(x, x_i)$, $0 < i \leq n$.

The precise behavior of the volume growth function $m_G(B_t)$ for balls defined by an admissible metric $d$ on an lcsc group $G$ is a very important characteristic from several perspectives, and is fundamental in consideration related to ergodic theorems, as we shall see below. The volume growth problem has seen some important recent progress, which we will describe in §§4.4 and 4.5. But before doing so, let us first formulate the following basic problem in ergodic theory.

### 4.3. The ball averaging problem in ergodic theory.
Our discussion so far, including of course the description of the basic pointwise ergodic theorems of Birkhoff and Wiener for Abelian groups, leads naturally to the formulation of the following problems.

Let $G$ be an lcsc group, $d$ an admissible metric on $G$, $B_t$ the corresponding balls of radius $t$ and center $e$, and let $\beta_t$ be normalized ball averages on $G$. By this we mean that $\beta_t$ are the absolutely continuous probability measures on $G$ whose density with respect to left Haar measure $m_G$ on $G$ is the characteristic function of $B_t$, normalized by $m_G(B_t)$.

When $G$ is a connected Lie group, a particularly basic case arises when $d$ is the distance function on $G$ associated to a left-invariant Riemannian metric on $G$. Note that in this case, we can also consider the spheres $S_t$, and the probability measures $\sigma_t$ supported on $S_t$, which is given canonically by the volume form arising from the restriction of the Riemannian metric to the spheres. In general, when $G$ is lcsc and the metric assumes integer values, we can consider the sequence $\sigma_t$, $t \in \mathbb{N}$ of probability measures defined by the restriction of a left Haar measure to the sphere. Anticipating some of the developments of the succeeding sections, we formulate the following.

(1) **The ball averaging problem** : Establish when $\beta_t$ satisfies the pointwise ergodic theorem in $L^1$ (or at least in every $L^p$, $1 < p < \infty$), namely for every ergodic probability-preserving action

$$\lim_{t \to \infty} \pi(\beta_t) f(x) = \int_X f \, dm \text{ for almost all } x \in X$$

and in the $L^p$-norm.

(2) **Sphere averaging problem** : Establish when there exists $p_0 < \infty$ such that the sphere averages $\sigma_t$ satisfy the pointwise ergodic theorem in $L^p$, $p > p_0$, namely

$$\lim_{t \to \infty} \pi(\sigma_t) f(x) = \int_X f \, dm \text{ for almost all } x \in X$$

and in the $L^p$-norm.

(3) **Spherical differentiation.** For $G$ a connected Lie group, $d$ a Riemmanian metric, establish when the singular spherical differentiation theorem holds for $f \in L^p$, $p > p_0$, namely

$$\lim_{t \to 0} \pi(\sigma_t) f(x) = f(x) \text{ for a.e. x}.$$

It is evident that the first two problems raised above are inextricably linked with the study of volume growth for balls and spheres, to which we now turn.

### 4.4. **Exact volume growth.**

4.4.1. *Exact polynomial volume growth.* When $G$ has polynomial volume growth the discussion of §4.1 allows us only to conclude that there exists a positive number $q$, such that for any Haar measure and any admissible metric on $G$, we have the following estimate for the volume of balls defined by the metric : $m_G(B_t) \leq C(G, d) t^q$. When $G$ has exponential volume growth $h$, we can only conclude that $c_\varepsilon (h - \varepsilon)^t \leq m_G(B_t) \leq C_\varepsilon (h + \varepsilon)^t$ for every $\varepsilon > 0$. As we shall see below, these estimates are entirely unsatisfactory for many purposes, and so it is natural to introduce the following definitions (see [Gu, Def. I.2]).

**Definition 4.3.** (1) The pair $(G, d)$ is said to have exact $t^q e^{ct}$ volume growth if the balls $B_t$ defined by the left-invariant admissible metric $d$ satisfy

$$\lim_{t \to \infty} \frac{m_G(B_t)}{C t^q \exp ct} = 1$$

for some non-negative $q$, $c$ and $C$ (which are necessarily uniquely determined).

(2) The pair $(G, d)$ is said to have strict $t^q e^{ct}$ volume growth if

$$b t^q \exp ct \leq m_G(B_t) \leq B t^q \exp ct$$

for some non-negative $q$ and $c$ (which are necessarily uniquely determined) and two positive constants $b \leq B$.

Important recent progress has been obtained in establishing strict or exact volume growth in several context. We will describe the developments regarding exact growth in the present section, and regarding strict growth in the following one.

The most obvious examples of groups $G$ and metrics $d$ with exact $t^q$ (namely polynomial) volume growth are :

(1) $G = \mathbb{R}^k$, and $d$ the metric defined by any norm on $\mathbb{R}^k$.
(2) $G = H_n = \mathbb{C}^n \times \mathbb{R}$ the Heisenberg groups, with the metric determined by the homogeneous norm $N(z, t) = (\|z\|^4 + t^2)^{1/4}$
(3) More generally, $G$ a connected nilpotent homogeneous group, with $d$ the metric derived from a homogeneous norm.

In these examples the homogeneity of the metric on $G$, namely the fact that $B_t = \alpha_t(B_1)$ (where $\alpha_t$ is the dilation automorphism group) immediately implies that $m_G(B_t) = m_G(B_1) t^q$, where $q$ is the homogeneous dimension. Thus in particular the polynomial volume growth for homogeneous metrics is exact. Riemannian metrics are not homogeneous in general, and here, under one further assumption, the following result was established by P. Pansu.

**Theorem 4.4.** [Pa] *Let $G$ be a simply connected nilpotent Lie group admitting a lattice subgroup $\Gamma$. Then*

(1) *$G$ has exact polynomial volume growth w.r.t. any Riemannian metric on $G$ which is $\Gamma$-invariant.*

(2) *Any lattice $\Gamma$ in $G$ has exact polynomial volume growth with respect to any word metric.*

Pansu's theorem has the following two consequences [Pa] (see also the discussion in [dlH, CH. VII.C]).

**Corollary 4.5.**    (1) Every discrete group containing a finitely generated nilpotent group of finite index has exact polynomial volume growth with respect to any word metric. Indeed, by a well-known result of Malcev, $\Gamma$ can be embedded as a lattice in a simply connected nilpotent Lie group.

(2) Every discrete group $\Gamma$ with polynomial volume growth has exact polynomial volume growth with respect to any word metric. Indeed, by Gromov's theorem [G], $\Gamma$ is virtually nilpotent, namely contains a nilpotent subgroup of finite index.

We note that a sharper form of Pansu's theorem has been established by M. Stoll [Sto], namely that $|\#B_t - Ct^q| \leq Bt^{q-1}$ for 2-step nilpotent groups.

For a simply connected nilpotent Lie group $G$, the condition that $G$ contains a lattice subgroup is equivalent to the condition that the Lie algebra of $G$ admits a basis with rational structure constants. Thus there exists only a countable set of such groups. So while Pansu's theorem (in combination with Gromov's) imply exact volume growth for every discrete group with polynomial growth, it leaves much to be desired in the case of general connected nilpotent Lie groups (and more generally lcsc groups with polynomial growth). Very recently, E. Breuillard has obtained a solution to this problem in the simply connected nilpotent case, in the following form.

**Theorem 4.6. Exact polynomial volume growth for simply connected nilpotent Lie groups**[Br]. *Let $G$ be a simply connected nilpotent Lie group $G$, let $H$ be a closed co-compact subgroup (e.g. $H = G$), and let $d$ be an $H$-invariant Riemannian or word metric on $G$. Then $G$ has exact polynomial volume growth with respect to $d$. In fact, the conclusion holds for any locally bounded, proper, asymptotically geodesic metric $d$.*

The method of proof used in [Br] is motivated by [Pa], and is based on showing that the ratio of the volume of the balls defined by $d$ is aymptotically the same as the volume of the balls defined by an appropriate Carnot-Caratheodory metric on $G$. For the latter the scaling property of the volume is clear, and thus exactness for the volume of the $d$-balls follows. Furthermore, exactness is also proved in [Br] for metrics on some connected non-nilpotent Lie groups of polynomial volume growth. Here use is made of the construction of the nil-shadow of such a group, and this allows the reduction of the problem to the nilpotent case. In view of the structure theorem for general lcsc groups of polynomial volume growth which will be discussed further in §5.5, the methods developed in [Br] may well lead to a complete solution of the problem of establishing exact growth on an arbitrary lcsc group of polynomial volume growth.

4.4.2. *Exact exponential volume growth.* Important recent progress on exact growth has recently been obtained also for a class of groups with exponential volume growth, namely semisimple Lie groups. Let us first note that a very natural choice for a metric on semisimple groups is the bi-$K$-invariant Riemannian metric on the group which is associated with the Killing form ($K$ a maximal compact subgroup). The exact $t^{r-1} \exp 2 \|\rho\| t$ volume growth (see §10.4 for notation) of balls in this case follows, for example, from the sharper results of [Kn2] on the asymptotic volume of the spheres, but can also be proved more directly.

Another natural family of metrics is given as follows. Fix a linear representation $\tau$ of $G$ into $GL_n(\mathbb{R})$, fix any (vector-space) norm $\|\cdot\|$ on $M_n(\mathbb{R})$, and let $d(g,h) = \log(1 + \|\tau(g^{-1})\tau(h)\|)$. $d$ will be symmetric if $\|g^{-1}\| = \|g\|$, and satisfy the triangle inequality if $\|gh\| \leq \|g\| \|h\|$, e.g. if $\|\cdot\|$ is a symmetric operator norm of $M_n(\mathbb{R})$. The following general result has been established by A. Gorodnik and B. Weiss.

**Theorem 4.7. Exact volume growth for semisimple Lie groups**[GW]. *Let $G$ be a connected semisimple Lie group with finite center, $\tau : G \to GL_n(\mathbb{R})$ a linear representation, and $\|\cdot\|$ any norm on $M_n(\mathbb{R})$. Then the sets given by $G_t = \{g \in G ; d(e,g) \leq t\}$ have exact $t^q \exp ct$ volume growth, for some $q$ and $c > 0$ depending on the representation $\tau$ and the norm $\|\cdot\|$.*

### 4.5. Strict volume growth.

4.5.1. *Strict polynomial growth.* A thorough investigation of growth as well as strict growth has been conducted by Y. Guivarc'h [Gu] (see also [Je] and [Ba]). For connected Lie groups with polynomial volume growth, and so in particular for nilpotent Lie groups, the fundamental result on strict $t^q$-volume growth for word metrics is as follows.

**Theorem 4.8. Strict volume growth for nilpotent groups**[Gu, Thm. II.3]. *Let $G$ be any connected Lie group of polynomial growth. Then any word metric, and thus also any metric quasi-isometric to a word metric, and in particular invariant Riemannian metric has strict polynomial volume growth.*

Furthermore, it is noted in [Gu] that since a finitely generated torsion free nilpotent group can always be embedded in a simply connected Lie group, strict volume growth holds for word metrics on countable nilpotent groups. Thus, using Gromov's theorem [G], strict growth holds for countable discrete groups of polynomial volume growth.

4.5.2. *The volume doubling condition.* One fundamental consequence of strict polynomial volume growth is the volume doubling condition, introduced by Calderon. Since our discussion centers around the ball averaging problem, we formulate it in the following form.

**Definition 4.9. Calderon's doubling volume condition**[C1]. *$G$ is said to satisfy the doubling volume condition w.r.t. the invariant admissible metric $d$ if the balls $B_t = \{g ; d(g,e) \leq t\}$ satisfy $m_G(B_{2t}) \leq C(G,d)m_G(B_t)$, for all $t > 0$. Here $m_G$ denote some (left or right) Haar measure on $G$.*

*Remark* 4.10.       (1) Clearly the volume doubling condition immediately implies that $m_G(B_{2^n}) \leq C(G,d)^n m_G(B_1)$, and so it follows that $(G,d)$ has polynomial volume growth and is therefore unimodular.

(2) Clearly when $(G, d)$ satisfies exact polynomial volume growth, the balls satisfy the doubling condition. But in fact, for the doubling condition already strict polynomial volume growth is sufficient. Thus for connected Lie groups of polynomial growth or for discrete nilpotent groups, volume doubling follows from Theorem 4.8.

4.5.3. *Strict exponential volume growth.* For semisimple Lie groups, using the polar coordinates on semisimple group, it can be easily established that any bi-$K$-invariant Riemmanian metric ($K$ a maximal compact subgroup) have strict $t^q \exp ct$ volume growth, with $q$ and $c > 0$ depending on the metric. More generally, a bi-$K$-invariant metric on $G$ which restricts to a Weyl-group-invariant norm on the Lie algebra of a (split) Cartan subgroup (see §10.4 for the definitions), will be called norm-like. The following recent result of H. Abels and G. Margulis [AM] establishes that norm-like metrics are the yardsticks for very general metrics on $G$, as follows.

**Theorem 4.11. Word metrics on semisimple groups are coarsely isometric to norm-like metrics**[AM]. *Let $G$ be a connected semisimple Lie group with finite center, and $d_V$ a word metric on $G$, associated with a bounded symmetric open set. Then there exists a bi-$K$-invariant norm-like metric $d$ such that $|d_V - d|$ is bounded on $G \times G$.*

We note that in fact the last result holds more generally for all reductive algebraic groups, and for every left-invariant coarsely geodesic quasi-metric on $G$ satisfying certain natural properness and boundedness condition, and refer to [AM] for the details.

**Corollary 4.12.** *The balls associated with a norm-like distance on a semisimple Lie group have strict $t^q \exp ct$ volume growth, and hence also the balls defined by $d_V$ have the same property.*

As to discrete groups with exponential volume growth, let us recall the following result due to M. Coornaert [Coo]. Let $\Gamma$ be a word-hyperbolic group, and $S$ a finite symmetric set of generators. Then the spheres, and hence also the balls, have strict $\exp ct$ volume growth, namely $b \exp ct \leq |S_t| \leq B \exp ct$.

Finally, let us note that establishing exact, and even strict volume growth is an open problem for most other groups and metrics. For some further results in this direction regarding discrete groups we refer to the [dlH].

4.6. **Balls and asymptotic invariance under translations.** Before formulating the ergodic theorems for balls on groups with polynomial volume growth, consider the following properties of the family of balls in $G$. These properties obviously hold under the assumption of exact polynomial volume growth, and they are easily seen to be mutually equivalent.

**Proposition 4.13.** *Let $(G, d)$ have exact polynomial volume growth w.r.t. the balls $B_t$ defined by an admissible invariant metric $d$ on $G$. Then the following properties hold*

(1) *The volume of the ball is asymptotically stable, namely for every $r > 0$,*

$$\lim_{t \to \infty} \frac{m_G(B_{t+r})}{m_G(B_t)} = 1 .$$

(2) *The volume of a shell is asymptotically negligible when compared with the volume of the ball, namely for every $r > 0$,*

$$\lim_{t \to \infty} \frac{m_G(B_{t+r} \setminus B_t)}{m_G(B_t)} = 0 \,.$$

(3) *The ball averages are asymptotically invariant under (right) convolution, namely for every $r > 0$,*

$$\lim_{t \to \infty} \|\beta_t * \beta_r - \beta_t\|_{L^1(G)} = 0 \,.$$

(4) *The balls are asymptotically uniformly invariant under (right) traslations, namely for every compact set $Q \subset G$,*

$$\lim_{t \to \infty} \frac{m_G(B_t \cdot Q \Delta B_t)}{m_G(B_t)} = 0 \,.$$

(5) *The balls are asymptotically invariant under (right) translations, namely for every $g \in G$,*

$$\lim_{t \to \infty} \frac{m_G((B_t g) \Delta B_t)}{m_G(B_t)} = 0 \,.$$

This very pleasant property of balls is discussed in F. Greenleaf [Gr1], where the question is raised whether it holds for balls defined by a word metric, for all lcsc groups admitting *any* sequence of asymptotically invariant compact sets of positive finite Haar measure. It is referred to in [Gr1] as the localization conjecture, as it locates a specific asymptotically invariant sequence in the group - namely the powers of a compact neighborhood of the identity. In this generality the conjecture is false, and in fact fails already for the $ax + b$ group, as shown in [Mi]. In fact, a general result is that for connected exponential solvable Lie groups, no subsequence of the sequence of balls w.r.t. a Riemannian metric (say) can be asymptotically invariant, as shown in [Pit].

Nevertheless, for some groups of polynomial growth the localization conjecture was very recently given two very interesting independent solutions, each yielding significantly more information that just asymptotic invariance. One solution is due to E. Breuillard , applies to all word metrics on simply connected nilpotent Lie groups (as well as some further Lie groups of polynomial volume growth) and in fact gives the sharper result that the volume growth of balls is exact, and hence they are of course asymptotically invariant. In fact the result applies to more general metrics - see Theorem 4.6.

Another solution is due to R. tessera [Te], applies also to all word metrics on connected Lie groups with polynomial volume growth, and also yields a result sharper than asymptotic invariance. Indeed, the fact that the volume of the shell of width $r$ namely $m_G(B_{t+r} \setminus B_t)$ is asymptotically negligible when compared to the volume of the ball $m_G(B_t)$, can be given a precise quantitative form, as follows.

**Theorem 4.14. Balls satisfying the doubling condition are asymptotically invariant**[Te]. *Let $G$ be an lcsc group, and let $d$ be a word metric satisfying the doubling condition. Then there exist positive constants $\delta$ and $C$, such that $m_G(B_{t+1} \setminus B_t) \leq C t^{-\delta} m_G(B_t)$ for all $t \in \mathbb{N}$. In particular, the sequence of balls $B_t$ is asymptotically invariant.*

In fact, in [Te] a more general result is proved, namely the same estimate is established for every monotone metric on a metric-measure space. We remark that the doubling condition follows from strict growth, a result established by Y. Guivarc'h for connected Lie groups of polynomial volume growth and stated in Theorem 4.8.

For completeness, we indicate here the elegant elementary argument given in [Te] for word metrics in the group case, which proves Theorem 4.14.

*Proof.* Let $V \subset G$ be a compact symmetric neighbourhood of the identity, and $d$ the corresponding left invariant word metric. Let $B_n$ be the sequence of ball centered at the identity, and let $C_{n,n+k} = B_{n+k} \setminus B_n$ be the shell of width $k$ and inner radius $n$. We are interested in comparing the size of the shell $C_{n-1,n}$ and the ball $B_n$. To do that, first note that the doubling condition is equivalent to the statement that for some fixed positive constant $c$ independent of $n$, $|B_n| = |C_{0,n}| \geq c\,|C_{n,2n}| = |B_{2n} \setminus B_n|$.

Thus a natural generalization of the doubling condition from balls to shells would be the estimate $|C_{m-k,m}| \geq c\,|C_{m,m+k}|$, for some fixed constant $c$ and all $k \leq m$.

Assuming this estimate for the moment, consider the following disjoint union of shell whose width doubles at each step :

$$D_j = C_{n-1,n} \cup C_{n-2,n-1} \cup C_{n-4,n-2} \cup \cdots \cup C_{n-2^j,n-2^{j-1}} = C_{n-2^j,n}$$

$D_{j+1}$ is obtained by adding one more shell to $D_j$, and the foregoing estimate (taking $k = 2^j$, $m = n-2^j$) implies that the last shell added to $D_j$ to produce $D_{j+1}$ satisfies

$$\left|C_{n-2^{j+1},n-2^j}\right| \geq c\,\left|C_{n-2^j,n}\right| = c\,|D_j|$$

Thus $|D_{j+1}| \geq (1+c)\,|D_j|$ and hence $|D_j| \geq (1+c)^j\,|C_{n-1,n}|$. Letting $i = [\log_2 n]$, note that clearly $D_i \subset B_n$ and hence

$$|B_n| \geq |D_i| \geq (1+c)^i\,|C_{n-1,n}| \geq \frac{1}{2}(1+c)^{\log n}\,|C_{n-1,n}| = \frac{1}{2}n^{\log_2(1+c)}\,|C_{n-1,n}|$$

Thus the size of the shell does indeed satisfy $|C_{n-1,n}| \leq B n^{-\delta}\,|B_n|$.

Now to obtain the estimate $|C_{m-k,m}| \geq c\,|C_{m,m+k}|$, one uses the doubling condition, as follows. First, note that by the triangle inequality, for $1 \leq k \leq n/4$

$$B_k \cdot C_{n-2k,n-2k+1} \subset C_{n-4k,n} \quad \text{and} \quad C_{n,n+4k} \subset B_{8k}C_{n-2k,n-2k+1}$$

Let $x_i, i \in I$ be a maximal $k$-net in $C_{n-2k,n-2k+1}$, namely such that $B_k(x_i) \cap B_k(x_j) = \emptyset$ if $i \neq j$. By maximality, we have $C_{n-2k,n-2k+1} \subset \cup_{i \in I} B_{2k}(x_i)$ and therefore $C_{n,n+4k} \subset \cup_{i \in I} B_{10k}(x_i)$. However, the doubling condition clearly implies that $|B_{10k}| \leq A\,|B_k|$, and thus we conclude that

$$|C_{n-4k,n}| \geq |C_{n-2k,n-2k+1}| \geq \sum_{i \in I}|B_k(x_i)| \geq \frac{1}{A}\sum_{i \in I}|B_{10k}(x_i)| \geq \frac{1}{A}\,|C_{n,n+4k}|$$

The inequality $|C_{m-k,m}| \geq c\,|C_{m,m+k}|$ for all $k \leq m$ follows similarly. $\qquad\square$

## 5. Pointwise ergodic theorems for groups of polynomial volume growth

Our main purpose in this section is to give a complete account of the proofs of the following results, which generalize Birkhoff's, Wiener's and Calderon's pointwise ergodic theorems. We start with the following basic result, which is a variation on the classical results developed by Wiener [W], Riesz [R] and Calderon [C1], and relies on their arguments.

**Theorem 5.1. Pointwise ergodic theorem for groups with volume doubling, asymptotically invariant balls.** *Let $G$ be locally compact second countable group of polynomial volume growth with respect to the admissible metric $d$. If the family of balls $B_t$ satisfies the doubling condition and is asymptotically invariant under translations, then the family of ball averages $\beta_t$ defined by $d$ satisfies the pointwise ergodic theorem in $L^p$, $1 \le p < \infty$.*

Recall now that according to Theorem 4.6, balls w.r.t. word metrics on (say) simply connected nilpotent Lie groups have exact polynomial growth. Thus they satisfy the doubling condition and are asymptotically invariant under translations. Together with Theorem 5.1 this proves the following result, due to E. Breuillard.

**Theorem 5.2. Pointwise ergodic theorem for simply connected nilpotent Lie groups.**[Br]

(1) *Let $G$ be a simply connected nilpotent Lie group Let $d$ be the distance function derived from a $G$-invariant Riemannian metric, a homogeneous $G$-invariant metric, or a word metric, and $\beta_t$ the corresponding ball averages. Then $\beta_t$ satisfy the pointwise ergodic theorem in every $L^p$, $1 \le p < \infty$.*

(2) *The same result holds for any finitely generated nilpotent group, with $d$ any word metric.*

We note that the result for discrete nilpotent groups follows already from Pansu's Theorem (Theorem 4.4).

An alternative proof of a generalization of Theorem 5.2 follows by combining Guivarch's results on growth, and Tessera's result on asymptotic invariance. Recall that Theorem 4.14 asserts that the doubling condition implies asymptotic invariance. Theorem 4.8 establishes that connected Lie groups of polynomial volume growth have strict polynomial growth and thus satisfy the doubling condition. Thus together with Theorem 5.1, these results imply pointwise convergence in $L^1$ of ball averages, for every Lie group of polynomial volume growth.

In fact, utilizing fully the results concerning growth in [Gu], together with Losert's structure theorem [Lo] for lcsc groups of polynomial volume growth based on Gromov's theorem [G], it is possible to give a complete solution to the ball averaging problem for all metrics (quasi-isometric to) word metrics on all lcsc groups of polynomial volume growth. We formulate this results as follows, and will outline its proof in §5.5.

**Theorem 5.3. Pointwise ergodic theorem for groups with polynomial volume growth.** *For every locally compact second countable group $G$ of polynomial volume growth, the family of ball averages $\beta_t$ defined by any word metric satisfies the pointwise ergodic theorem in $L^p$, $1 \le p < \infty$. The same holds true for the balls determined by any left-invariant metric quasi-isometric with a word metric.*

In the next four subsections of the present chapter we will give a complete proof of Theorem 5.1, demonstrating the four steps called for by the recipe of §2.3. In the fifth we will describe the ingredients needed to complete the proof of Theorem 5.3.

**Proof of Theorem 5.1.**

Given a group of polynomial volume growth with asymptotically invariant balls, the plan of the proof is to follow the four steps outlined in the recipe of §2.3, and we start with 5.1. **Step I : The mean ergodic theorem.** For the proof of the mean ergodic theorem we will use here a standard variant of F. Riesz's [R] classical proof of Von-Neumann's mean ergodic theorem.

First let us note that $L^2(X)$ is the direct (orthogonal) sum of the subspaces $\mathcal{I}$ consisting of $G$-invariant vectors, and $\mathcal{K} = \overline{\mathrm{span}}\{(\pi(g) - I)f\,;\, f \in L^2(X)\,,\, g \in G\}$. Now note that if $h \in \mathcal{K}$ is of the form $h = \pi(g)f - f$, then clearly,

$$\|\pi(\beta_t)h\|_{L^2(X)} = \|(\pi(\beta_t)\pi(g) - \pi(\beta_t))f\|_{L^2(X)} \leq \|\beta_t * \delta_g - \beta_t\|_{L^1(G)} \|f\|_{L^2(X)}$$

But since by Proposition 4.13 (4), as $t \to \infty$

$$\frac{1}{m_G(B_t)} m_G(B_t g \Delta B_t) \to 0$$

we have $\lim_{t\to\infty} \pi(\beta_t)h = 0$ in $L^2$-norm for a dense set of $h \in \mathcal{K}$, hence for all $h \in \mathcal{K}$, by an obvious approximation argument.

Since clearly $\pi(\beta_t)f = f$ for every $f \in \mathcal{I}$, we conclude that for every $f \in L^2(X) = \mathcal{K} \oplus \mathcal{I}$

$$\lim_{t\to\infty} \|\pi(\beta_t)f - \mathcal{E}f\|_{L^2(X)} = 0$$

where $\mathcal{E}$ is the projection on the space $\mathcal{I}$ of $G$-invariant vectors, so that the mean ergodic theorem holds in $L^2(X)$. The mean ergodic theorem also holds in every $L^p(X)$, $1 \leq p < \infty$ by standard approximation argument, using the fact that $L^\infty$ is norm-dense in every $L^p$.

5.2. **Step II : Pointwise convergence on a dense subspace.** Let $1 \leq p < \infty$, and consider the space

$$\mathcal{K}' = \mathrm{span}\{h = \pi(g)f - f\,;\, f \in L^\infty(X)\,,\, g \in G\}$$

and the space $\mathcal{I}$ of $G$-invariant functions in $L^\infty(X)$. The sum of these two spaces is dense in $L^p(X)$, as follows from the following two facts. First, every $u \in L^q(X)$ (where $\frac{1}{p} + \frac{1}{q} = 1$) which integrates to zero against every function in $\mathcal{K}'$ is a $G$-invariant function, since $L^\infty$ is norm-dense in every $L^p$. Second, since $u$ is measurable w.r.t. the $\sigma$-algebra of $G$-invariant functions, if it integrates to zero against the characteristic function of every $G$-invariant set, then necessarily $u = 0$.

Now again as in Riesz's argument of §5.2, if $h = \pi(g)f - f$ then for almost every $x \in X$

$$|\pi(\beta_t)h(x)| = |(\pi(\beta_t * \delta_g - \beta_t)f(x)| \leq \frac{2\|f\|_{L^\infty(X)}}{m_G(B_t)} m_G(B_t g \Delta B_t) \to 0\,.$$

Thus $\pi(\beta_t)f(x) \to \int_X f\,dm$ almost everywhere for every $f$ in the dense subspace $\mathcal{K}' \oplus \mathcal{I}$ of $L^p(X)$.

*Remark* 5.4. Anticipating some arguments needed in the sequel, we note that the proof of the mean ergodic theorem is based solely on property (5) in Proposition 4.13 above, namely asymptotic invariance under translation. The same remark applies to the proof of the existence of a dense subspace where pointwise convergence holds.

**5.3. Step III : The maximal inequality for ball averages.** The maximal inequality for ball averages will be established in two stages, in §5.4.1 and §5.4.2. First the maximal inequality for the special action of the ball averages by convolution on the group manifold will be established, using the volume doubling condition which holds in all groups of strict polynomial volume growth. Then a transfer principle will be formulated, which will allow us to deduce the maximal inequality for an arbitrary action from its validity for convolutions.

5.3.1. *Maximal inequality for convolutions : the volume doubling condition.* The method of proof of the weak-type $(1, 1)$ maximal inequality for convolutions which we will present originates in Wiener's proof for balls in $\mathbb{R}^d$. It was later observed by Calderon [C1] that the proof only depends on the fact that the volume of a Euclidean ball of a given radius is bounded by a fixed multiple of the volume of a ball of half the radius. In [C1] this volume doubling condition is introduced for more general families of sets $N_t \subset G$, but here we will continue to focus on the ball averaging problem.

We now turn to the covering Lemma, which we reproduce in the original finite form given in [W, Lemma C'].

**Lemma 5.5. Wiener-Calderon covering argument**[W][C1]. *Assume $d$ is an admissible metric on a lcsc group $G$, which satisfies the doubling volume condition. Then every finite family of balls $\{B_{t_i}, i \in I\}$ in $G$ contains a subfamily $\{B_{t_j}, j \in J \subset I\}$ of disjoint balls whose total volume $m_G(\cup_{j \in J} B_{t_j})$, is at least $\delta \cdot m_G(\cup_{i \in I} B_{t_i})$, where $\delta = \delta(G) > 0$. Thus the total volume of the disjoint subcover is at least a fixed fraction of the total volume of the original family.*

*Proof.* The volume doubling condition obviously implies that $m_G(B_t) \geq \delta m_G(B_{3t})$. $I$ being finite, choose one of the balls in the family which has maximal radius, and label it $B_{t_1}$. Consider now the subfamily $\{B_{t'_{i'}}; i' \in I' \subset I\}$ of all balls intersecting $B_{t_1}$. The union of all the balls $B_{t'_{i'}}, i' \in I'$ is contained in a ball of radius at most three times that of $B_{t_1}$, with the same center. Therefore keeping $B_{t_1}$ and deleting all other balls intersecting it, we keep at least a fraction $\delta$ of the total volume of the subfamily $\{B_{t'_{i'}}; i' \in I'\}$. We therefore put the index $t_1$ in $J$, and apply the same argument again to the family $\{B_{t_i}; i \in I \setminus I'\}$, which consists only of balls disjoint from $B_{t_1}$. Proceeding finitely many times, we obtain a disjoint sequence of balls whose total volume occupies at least a fraction $\delta$ of the total volume $m_G(\cup_{i \in I} B_{t_i})$. $\square$

A variant of the previous argument is the following covering Lemma, which can be proved in much the same way - see [M, Ch. IV, §1] for a more general formulation.

**Lemma 5.6. Vitali covering Lemma**. *Assume $d$ is an admissible metric on an lcsc group $G$, whose balls satisfy the doubling volume condition. Given any set $A$ of positive Haar measure, there exists a disjoint sequence of balls $B_{t_i}, i \in \mathbb{N}$, satisfying $m_G(A \setminus \cup_{i \in \mathbb{N}} B_{t_i}) = 0$.*

Using Lemma 5.5, the following maximal inequality was proved by Calderon [C1] (for more general families of sets), following [W]. It generalizes the Hardy-Littlewood maximal inequality for averages on the real line, as well as Wiener's maximal inequality for balls in Euclidean space.

**Theorem 5.7. Maximal inequality for convolutions with ball averages satisfying the doubling condition [W][C1].** *Assume $d$ is an admissible metric on an lcsc group $G$, whose balls satisfy the doubling volume condition. Then the family of ball averages $\beta_t$ satisfies the weak-type $(1,1)$-maximal inequality for convolutions, given by*

$$m_G\{g \in G\,;\, \sup_{0<t<\infty} |F * \beta_t(g)| > \varepsilon\} \le \frac{C(G)}{\varepsilon}\, \|F\|_{L^1(G)}$$

*Proof.* Since $G$ is unimodular and the balls are symmetric, each of the measures $\beta_t$ is symmetric. The convolution operators are thus given by (see §2.1)

$$F * \beta_t(g) = \frac{1}{m_G(B_t)} \int_{B_t} F(gh) dm_G(h) = \frac{1}{m_G(B_t)} \int_{y \in B_t(g)} F(y) dm_G(y)$$

where $B_t(g) = g B_t(e)$ is the ball of radius $t$ and center $g$. We denote as usual $F_\beta^* = \sup_{t>0} |F * \beta_t(g)|$ Since $|F * \beta_t(g)| \le |F| * \beta_t(g)$ we can and will assume that $F \ge 0$, without loss of generality. Let $U_\varepsilon = \{g \in G\,;\, F_\beta^*(g) > \varepsilon\}$, and let $W \subset U_\varepsilon$ be compact. By definition, for each $w \in W$ there is a ball $B_{r_w}(w) = w B_{r_w}(e)$ with center $w$ and radius $r_w$ satisfying

$$m_G(B_{r_w}(e)) = m_G(B_{r_w}(w)) < \frac{1}{\varepsilon} \int_{y \in B_{r_w}(w)} F(y) dm_G(y)$$

There exists a finite covering of $W$ using the collection of open balls $w B_{r_w}(e)$, and let us denote a finite covering family by $\{B_i\,;\, i \in I\}$. By Lemma 5.5 we can choose a finite disjoint subfamily $\{B_j\,;\, j \in J\}$ whose union retains at least a fraction $\delta(G)$ of the measure of the union of the balls $B_i$, $i \in I$. Combining these two estimates,

$$m_G(W) \le m_G(\cup_{i \in I} B_i) \le \frac{1}{\delta(G)} m_G(\cup_{j \in J} B_j) \le$$

$$\le \frac{1}{\delta(G)\varepsilon} \int_{\cup_{j \in J} B_j} F(g) dm_G(g) \le \frac{1}{\delta(G)\varepsilon} \|F\|_{L^1(G)}\,.$$

Taking the supremum over all compact sets contained in $U_\varepsilon$ we conclude that the same estimate holds for the measure of the set $U_\varepsilon$ and this concludes the proof of the weak-type $(1,1)$ maximal inequality. $\square$

We have thus established the maximal weak-type $(1,1)$ inequality for the family of operators $F \mapsto F * \beta_t$ of right convolutions by the ball averages $\beta_t$.

*Remark* 5.8. The covering arguments of Lemmas 5.6 and 5.7 described above fail for groups with exponential volume growth. However, the maximal inequality for ball averages acting by convolutions is often true, as we shall see below. Thus for semisimple Lie groups $G$ the action of the ball averages by convolution on the symmetric space $G/K$ satisfies the weak-type $(1,1)$ maximal inequality, as shown in [Str]. There does not seem to be an lcsc group for which the weak type $(1,1)$ maximal inequality for admissible ball averages is known to fail, for convolutions or otherwise. We refer to [N7] for more information on this subject.

5.3.2. *Maximal inequality for general actions : the transfer principle.* We now come to a basic observation whose origin is in Wiener's proof of the maximal inequality for ball averages in actions of $\mathbb{R}^d$, and was subsequently considerably generalized and expanded (as we shall see below). Wiener recognized that in order to prove the maximal inequality for ball averages in a measure-preserving action of the group on a *general space*, it is sufficient to prove the maximal inequality for the *action of the group on itself by translation*, provided that the balls are asymptotically invariant. This observation was later termed the transfer principle, and the idea underlying it is the following (see [W, Proof of Thm. IV']). Apply the operator $\beta_r$ to the function $f^t$, which is $f$ restricted to the subset $B_t \cdot x$ of the $G$-orbit of $x$ in $X$. Replacing $\beta_r f(x)$ by $\beta_r f^t(x)$ produces an error which is controlled by the (normalized) difference in volumes between $B_{r+t}$ and $B_t$. Since $(|B_{r+t} \setminus B_t|)/|B_t| \to 0$, it follows that we may consider in effect (almost) every orbit individually. Thus maximal inequalities for functions on $X$ are reduced to maximal inequalities for *convolution operators* on the group manifold. We will give here first a very simple formulation of the transfer principle for strong maximal inequalities, and defer a more general formulation to §6.

Explicitly, let $d(g, h)$ be a left-invariant admissible metric on $G$, and let $|g| = d(e, g)$. Then clearly $|gh| \le |g| + |h|$. Consider as usual the family $\beta_t$ of probability measures on $G$ defined by the balls.

**Theorem 5.9. The transfer principle for ball averages on groups with polynomial volume growth.** *Suppose $\rho(\beta_t)$, $0 \le t \le r$ satisfy the strong $L^p$-maximal inequality*

$$\left\| \sup_{0 < t \le r} \rho(\beta_t) F \right\|_p \le C_p \|F\|_p$$

*for $F \in L^p(G)$, and some $1 < p < \infty$, where $\rho$ is the right regular representation. Then*

(1) *$\pi(\beta_t)$ satisfy the maximal inequality*

$$\left\| \sup_{0 < t \le r} |\pi(\beta_t)f| \right\|_p \le C_p \left( \frac{m_G(B_R)}{m_G(B_{R-r})} \right)^{\frac{1}{p}} \|f\|_p$$

*in $L^p(X)$, for any measure preserving action $\pi$ of $G$ on a $\sigma$-finite measure space $X$. Here $R \ge 2r$ is any positive number.*

(2) *In particular, if $\lim_{R \to \infty} m_G(B_R)/m_G(B_{R-r}) = 1$, and $\beta_t$, $0 < t < \infty$ satisfies the maximal inequality for right convolutions on the group manifold, then $\beta_t$ satisfies the maximal inequality in any measure-preserving action.*

*Proof.* Given $f \in L^p(X)$, fix $x \in X$ and define : $F_x(g) = f(g^{-1}x) = \pi(g)f(x)$ if $|g| \le R$ and $F_x(g) = 0$ otherwise. Clearly, if $|g| \le R - r, |h| \le r$ , then $F_x(gh) = \pi(g)\pi(h)f(x)$. Now integrate over $h$ w.r.t. the measure $\beta_t$, $t \le r$. Then we clearly have $\pi(g)\pi(\beta_t)f(x) = \rho(\beta_t)F_x(g)$, as long as $|g| \le R - r$. Taking the supremum over $0 \le t \le r$ of the $p$-th power, we obtain :

$$\pi(g) \sup_{0 < t \le r} |\pi(\beta_t)f(x)|^p = \sup_{0 < t \le r} |\rho(\beta_t)F_x(g)|^p \ .$$

Now integrate w.r.t. $g \in B_{R-r}$ (recall that $G$ is unimodular), and use the maximal inequality for $\rho(\beta_t)$, to get :

$$\int_{B_{R-r}} \pi(g) \sup_{0 < t \le r} |\pi(\beta_t)f(x)|^p \, dg \le \int_G \sup_{0 < t \le r} |\rho(\beta_t)F_x(g)|^p \, dg \le C_p^p \int_G |F_x(g)|^p \, dg.$$

Since $F_x(g)$ has support in $B_R$, the last integral equals $C_p^p \int_{B_R} |F_x(g)|^p \, dg$. Finally, we integrate over $X$ and use the fact that $G$ is measure preserving, to obtain :

$$\int_{B_{R-r}} \int_X \pi(g) \sup_{0 < t \le r} |\pi(\beta_t)f(x)|^p \, dm \, dg \le C_p^p \int_{B_R} \int_X \left| f(g^{-1}x) \right|^p \, dm(x) \, dg .$$

Hence

$$\left\| \sup_{0 < t \le r} |\pi(\beta_t)f| \right\|_p \le C_p \left( \frac{|B_R|}{|B_{R-r}|} \right)^{\frac{1}{p}} \|f\|_p .$$

This concludes the proof of the maximal inequality stated in part (1). The proof of (2) is an immediate consequence. □

*Remark* 5.10. A similar argument establishes also the weak-type $(1,1)$ maximal inequality for the ball averages. We have elected to present first the proof of the strong $L^p$-maximal inequality, for simplicity. The transfer of the weak-type $(1,1)$ maximal inequality will be demonstrated in greater generality in §6 below.

5.4. **Step IV : Interpolation arguments.** Thus far, taking Remark 5.10 for granted, we have established the weak-type $(1,1)$ maximal inequality for general measure-preserving actions of $G$, together with the mean ergodic theorem and pointwise convergence on a dense subspace. According to the recipe of §2.3, Theorem 5.1 has therefore been established for $f \in L^1(X)$. To conclude the proof of Theorem 5.1 is an easy matter and it remains only to note the following.

The family $\beta_t$ consists of Markov operators, and clearly each $\pi(\beta_t)$ has norm bounded by 1 as an operator on $L^1(X)$ and $L^\infty(X)$. It is clear that

$$\sup_{t > 0} |\pi(\beta_t)f(x)| = f_\beta^*(x) \le \|f\|_{L^\infty(X)} .$$

Given the weak-type $(1,1)$-maximal inequality for ball averages in an arbitrary measure-preserving action, by Marcinkiewicz's interpolation theorem, $f_\beta^*$ satisfies the strong $L^p$-maximal inequality for $1 < p < \infty$.

This concludes the proof of Theorem 5.1. □

5.5. **Groups of polynomial volume growth : General case.** We now proceed with

*Proof of Theorem 5.3.*

According to Theorem 5.1, to show that the balls of a given admissible metric $d$ on a group $G$ with polynomial volume growth satisfy the pointwise ergodic theorem, it suffices to show that they are asymptotically invariant under translations, and volume doubling. By Theorem 4.14, in fact a sufficient condition for asymptotic invariance is that the balls satisfy the volume doubling property. As noted already above, if the balls have strict polynomial volume growth, then they clearly satisfy the doubling volume condition. Furthermore, if the balls for a given metric have strict polynomial growth, then it is clear that the balls with respect to any quasi-isometric metric also have strict polynomial growth, and thus also satisfy the doubling condition. Therefore Theorem 5.3 will be completely proved once we establish the following

**Proposition 5.11.** *Given any lcsc group $G$ with polynomial volume growth, $G$ has strict polynomial volume growth w.r.t. any metric quasi-isometric to a word metric.*

*Proof.* The argument utilizes the definitive results on strict polynomial growth of Lie groups, obtained by Y. Guivarc'h [Gu], together with a structure theorem for lcsc groups with polynomial volume growth obtained by V. Losert [V], generalizing Gromov's result [G] in the discrete case.

Thus according to [V, Thm 2] any lcsc group with polynomial volume growth $G$ admits a normal series $C \triangleleft R \triangleleft N \triangleleft G$ with $C$ and $G/N$ compact, $R/C$ a connected solvable Lie group of polynomial volume growth, and $N/R$ a finitely generated discrete nilpotent group.

By [Gu, Thm I.4], $G$ has strict growth if $G/C$ does, so we may assume $C = \{e\}$. $R$ is a closed subgroup of $G$, and thus has polynomial growth by [Gu, Thm. I.2]. Now by [Gu, Thm. I.4], since $N$ is normal and co-compact in $G$, it has a growth function equivalent to that of $G$, and hence $N$ has polynomial growth, and $G$ has strict growth if $N$ does. So it suffices to show that $N$ has strict growth, and since $R$ is solvable, and $N/R$ nilpotent, it follows that $N$ is a solvable Lie group with connected component $R$ and polynomial volume growth. Thus by [Gu, Thm III.5] it follows that $N$ has strict polynomial volume growth.     $\square$

We remark that Proposition 5.11, in combination with Theorem 4.14, implies of course the following fact, which we record for completeness.

**Proposition 5.12.** *Let $G$ be an lcsc group with polynomial volume growth, $d$ any metric quasi-isometric to a word metric, and $B_t$ the corresponding balls. Then*

(1) *The balls $B_t$ are asymptotically invariant under translation.*
(2) *The shells $C_t = B_{t+1} \setminus B_t$ satisfy $m_G(C_t) \leq C t^{-\delta} m_G(B_t)$, for some positive $\delta$ and $C$, for $t \geq 1$.*

5.5.1. **Subsequence theorem for groups with subexponential growth.** Let us note the following regarding pointwise ergodic theorems for subsequences of ball averages.

(1) Calderon's original formulation [C1] of his pointwise ergodic theorem did not prove or assume strict polynomial volume growth, but instead noted that the doubling condition implies the following property for the volume of the balls. There exists a set $D \subset \mathbb{R}$ of density 1, such that for any $s > 0$, we have $\lim_{t \to \infty} m_G(B_{t \pm s})/m_G(B_t) = 1$. Thus, using the arguments above, it is shown in [C1] that the pointwise ergodic theorem holds for $\pi(\beta_t)f(x)$, provided that as $t \to \infty$ it assumes only values from the set $D$. Thus it was shown in [C1] that in every group satisfying the doubling condition, there is a subsequence $\beta_{t_n}$ satisfying the pointwise ergodic theorem in $L^1$.
(2) Similarly, subexponential growth (which is equivalent to polynomial volume growth in the connected Lie group case but not in general) implies that a subsequence of the sequence of balls is asymptotically invariant. Hence in particular the mean ergodic theorem is true for the subsequence.
(3) For a discrete subgroup of exponential growth, it follows easily from the definition that no subsequence of balls can be a Følner sequence. We recall that it has been established in [Pit] that when $G$ is a connected solvable Lie group and has exponential volume growth, again no subsequence of the sequence of balls is asymptotically invariant.

## 6. Amenable groups : Følner Averages and their applications

6.1. **The transfer principle for amenable groups.** As we shall see presently, inspection of the proof of Theorem 5.9 reveals that the condition which is essential for the transfer principle to hold (for an arbitrary family $\mu_t$ of probability measures with compact supports, not just $\beta_t$) is the existence of a sequence of sets $F_n$ (given there by the balls $B_n$) with the following property.

**Definition 6.1. Følner conditions**. A sequence $F_n$ of compact sets of positive Haar measure in an lcsc group $G$ is called

(1) (right) **Følner sequence** if for every $g \in G$

$$\lim_{n \to \infty} \frac{\eta(F_n g \Delta F_n)}{\eta(F_n)} = 0 \,.$$

(2) (right) **uniform Følner sequence** if for any given compact set $Q \subset G$ we have,

$$\lim_{n \to \infty} \frac{\eta(F_n Q \Delta F_n)}{\eta(F_n)} = 0 \,.$$

It then follows also that for every compact set $Q$

$$\lim_{n \to \infty} \frac{\eta(F_n Q)}{\eta(F_n)} = 1$$

where $\eta$ is right Haar measure on $G$.

The fact that the existence of a Følner sequence is sufficient for the validity of a transfer principle (generalizing that of Wiener [W]) was proved in various different formulations by a number of authors, starting with Calderon [C2], followed by Coifman-Weiss [CW], Emerson [E1], Herz [He1] and Tempelman [T2].

The existence of a Følner sequence in an lcsc group $G$ is equivalent to $G$ being amenable, and as is well-known, polynomial volume growth implies amenability, but there are many amenable groups of exponential volume growth. Thus the Følner condition yields a significantly more general transfer principle than Theorem 5.9. We now turn to the formulation and proof of a version of the transfer principle for amenable groups that will be found useful below.

**Theorem 6.2. The transfer principle for amenable groups**. *Let $\mu_t$, $0 < t < \infty$ be probability measures with compact supports on an lcsc group $G$. Assume that for $t \le R$ we have $\mathrm{supp}(\mu_t) \subset Q$, where $Q$ is a compact subset (possibly depending on $R$). Suppose $\eta$ is right Haar measure on $G$, and that $\mu_t, 0 \le t \le R$, satisfy the following maximal inequality for $F \in L^p(G, \eta)$, where $1 < p < \infty$*

$$\left\| \sup_{0 < t \le R} \left| \int_G F(gh) d\mu_t(h) \right| \right\|_{L^p(G,\eta)} \le C_p \, \|F\|_{L^p(G,\eta)} \,.$$

(1) *If $A$ is any compact set of positive measure in $G$, then $\pi(\mu_t)$ satisfy the strong $L^p$-maximal inequality*

$$\left\| \sup_{0 < t \le R} |\pi(\mu_t)f| \right\|_p \le C_p \left( \frac{\eta(AQ)}{\eta(A)} \right)^{\frac{1}{p}} \|f\|_p$$

*for any measure preserving action $\pi$ of $G$ on a $\sigma$-finite measure space $(X, m)$.*

(2) *If the weak-type $(1,1)$-maximal inequality holds for the action by translation, namely for $F \in L^1(G, \eta)$*

$$\eta \left\{ g \in G \, ; \, \sup_{0 < t \leq R} |F(gh)d\mu_t(h)| > \delta \right\} < \frac{C}{\delta} \|F\|_{L^1(G, \eta)}$$

*Then for any measure-preserving action on $(X, m)$ the following weak-type $(1,1)$-maximal inequality holds (for any $A$ as in (1))*

$$m \left\{ x \in X \, ; \, \sup_{0 < t \leq R} |\pi(\mu_t)f(x)| > \delta \right\} < \frac{\eta(AQ)}{\eta(A)} \frac{C}{\delta} \|f\|_{L^1(X)}$$

(3) *When $A_n$ satisfies $\lim_{n \to \infty} \eta(A_n \cdot Q)/\eta(A_n) = 1$, the maximal inequalities for $\mu_t$ hold in $L^p(X)$, $1 < p < \infty$ (with the same bound as in $L^p(G)$) namely*

$$\left\| \sup_{0 < t \leq R} |\pi(\mu_t)f| \right\|_{L^p(X)} \leq C_p \|f\|_{L^p(X)} \, .$$

*and*

$$m \left\{ x \in X \, ; \, \sup_{0 < t \leq R} |\pi(\mu_t)f(x)| > \delta \right\} < \frac{C}{\delta} \|f\|_{L^1(X)}$$

(4) *Finally, if $A_n$ satisfy the foregoing condition for every compact $Q$ (namely if $A_n$ form a right uniform Følner sequence) and the assumption regarding the maximal inequalities for convolutions is satisfied for $R = \infty$, then so is the conclusion.*

*Proof. Proof of Part (1).*

Given $f \in L^p(X)$, $1 < p < \infty$, fix $x \in X$ and a compact set $A \subset G$ of positive measure, and define:

$$F_A(g) = f(g^{-1}x) = \pi(g)f(x) \text{ if } g \in A \, , \text{ and } F_A(g) = 0 \text{ otherwise.}$$

Clearly, if $k \in A$ and $h \in Q$, then $F_{AQ}(kh) = \pi(k)\pi(h)f(x)$.

By assumption $Q$ contains the support of $\mu_t$, $0 < t \leq R$, and we can therefore integrate the last equation over $h \in Q$ w.r.t. the measure $\mu_t$ and write, as long as $k \in A$

$$\pi(k)\pi(\mu_t)f(x) = \int_G F_{AQ}(kh)d\mu_t(h) \, .$$

Taking the supremum over $0 < t \leq R$ of the $p$-th power, we obtain :

$$\pi(k) \sup_{0 < t \leq R} |\pi(\mu_t)f(x)|^p = \sup_{0 < t \leq R} \left| \int_G F_{AQ}(kh)d\mu_t(h) \right|^p \, .$$

Now integrate over $k \in A$ using right-invariant Haar measure $\eta$, and extend the integration to all of $G$ on the right hand side. This yields the obvious inequality

$$\int_A \pi(k) \sup_{0 < t \leq R} |\pi(\mu_t)f(x)|^p \, d\eta(k) \leq \int_G \sup_{0 \leq t \leq R} \left| \int_G F_{AQ}(kh)d\mu_t(h) \right|^p d\eta(k)$$

using the strong $L^p$-maximal inequality which we assumed for the r.h.s., together with the fact that (by definition) $F_{AQ}$ is supported in $AQ \subset G$, we obtain that the last integral is bounded by :

$$C_p^p \int_G |F_{AQ}(g)|^p \, d\eta(g) = C_p^p \int_{AQ} |F_{AQ}(g)|^p \, d\eta(g) \, .$$

Finally, we integrate both sides of the inequality over $X$, and use Fubini's theorem to obtain :

$$\int_A \int_X \pi(k) \sup_{0 < t \leq R} |\pi(\mu_t) f(x)|^p \, dm(x) d\eta(k) \leq$$

$$\leq C^p \int_{AQ} \int_X \left| f(g^{-1} x) \right|^p dm(x) d\eta(g) \ .$$

Hence since the $G$-action is measure preserving :

$$\left\| \sup_{0 < t \leq R} |\pi(\mu_t) f| \right\|_{L^p(X)} \leq C_p \left( \frac{\eta(AQ)}{\eta(A)} \right)^{\frac{1}{p}} \|f\|_{L^p(X)} \ .$$

This concludes the proof of part (1) of the theorem.

*Proof of Part (2).*

As to part (2), fix a compact set $A$, and let us define the set

$$\mathcal{D}(\delta) = \left\{ (k, x) \in A \times X \, ; \, \sup_{0 < t \leq R} \pi(k) \, |\pi(\mu_t) f(x)| > \delta \right\} \ .$$

The first coordinate sections of $\mathcal{D}(\delta)$ are given, for each $k \in A$, by

$$\mathcal{D}^k(\delta) = \left\{ x \in X \, ; \, \sup_{0 < t \leq R} \pi(k) \, |\pi(\mu_t) f(x)| > \delta \right\} \ .$$

In particular, the set whose measure we are interested in estimating is

$$\mathcal{D}^e(\delta) = \left\{ x \in X \, ; \, \sup_{0 < t \leq R} |\pi(\mu_t) f(x)| > \delta \right\}$$

Note that clearly, $\mathcal{D}^k(\delta) = k(\mathcal{D}^e(\delta))$, and since each $k \in G$ is measure preserving, we have $m(\mathcal{D}^k(\delta)) = m(\mathcal{D}^e(\delta))$ for every $k \in A$.

The second coordinate sections of $\mathcal{D}(\delta)$ are given, for each $x \in X$, by

$$\mathcal{D}_x(\delta) = \left\{ k \in A \, ; \, \sup_{0 < t \leq R} \left| \int_G \pi(kh) f(x) d\mu_t(h) \right| > \delta \right\} \ .$$

By Fubini's theorem, we have

$$\eta \times m(\mathcal{D}_\delta) = \int_X \eta(\mathcal{D}_x(\delta)) dm(x) = \int_A m(\mathcal{D}^k(\delta)) d\eta(k) = \eta(A) m(\mathcal{D}^e(\delta))$$

Now by assumption, the action by translation satisfies the weak-type $(1,1)$-maximal inequality. Keeping the notation introduced in the proof of part (1), we have $|\pi(k) \pi(\mu_t) f(x)| \leq \int_G |F|_{AQ} (kh) d\mu_t(h)$ if $k \in A$, $h \in Q$. Hence

$$\eta(\mathcal{D}_x(\delta)) \leq \eta \left\{ k \in A \, ; \, \sup_{0 < t \leq R} \left| \int_G |F_{AQ}(kh)| \, d\mu_t(h) \right| > \delta \right\}$$

Combining the two foregoing arguments, we conclude

$$\eta(A) m(\mathcal{D}^e(\delta)) = \int_X \eta(\mathcal{D}_x(\delta)) dm(x) \leq \int_X \frac{C}{\delta} \|F_{AQ}\|_{L^1(G, \eta)} \, dm(x)$$

finally, using the fact that $G$ is measure-preserving and the definition of $F_{AQ}$, we obtain

$$m(\mathcal{D}^e(\delta)) = m \left\{ x \in X \, ; \, \sup_{0 < t \leq R} |\pi(\mu_t) f(x)| > \delta \right\} \leq \frac{\eta(AQ)}{\eta(A)} \frac{C}{\delta} \|f\|_{L^1(X)}$$

and the proof of part (2) is complete.

*Proof of Part (3).*

Part (3) follows immediately upon applying the following argument. If $A_n$ satisfies for every compact set $Q \subset G$

$$\lim_{n \to \infty} \frac{\eta(A_n Q)}{\eta(A_n)} = 1$$

we take the limit as $n \to \infty$ in part (1) and conclude that for each $Q$,

$$\left\| \sup_{0 < t \leq R} |\pi(\mu_t) f| \right\|_{L^p(X)} \leq C_p \|f\|_{L^p(X)} \ .$$

*Proof of Part (4).*

Since $C_p$ is fixed and independent of $Q$, we therefore choose a sequence of compact sets $Q_m \subset Q_{m+1}$ whose union is $G$. Thus $\mathrm{supp}(\mu_t) \subset Q_m$ for $t \leq R_m$, where $R_m \to \infty$. Applying the foregoing result to each set $Q_m$, a straightforward application of the monotone convergence theorem or Fatou's Lemma allows us to conclude that also

$$\|f_\mu^*\|_{L^p(X)} = \left\| \sup_{0 < t < \infty} |\pi(\nu_t) f| \right\|_{L^p(X)} \leq C_p \|f\|_{L^p(X)}$$

A similar argument proves the corresponding result for the weak-type $(1, 1)$ maximal inequality, and for the case $R \to \infty$.

This concludes the proof of the transfer principle of amenable groups.  $\square$

*Remark* 6.3. The formulation of the transfer principle in Theorem 6.2 is similar to the one given by Tempelman in [T, Ch. 5,§1.4]. It differs somewhat from those of Calderon [C1], Emerson [E1] and Coifman-Weiss [CW]. The latter formulations all consider the transfer of an arbitrary operator $T$ on $L^1_{loc}(G)$ satisfying the following properties.

(1) $T$ is sublinear,
(2) $T$ commutes with right translations,
(3) $T$ is semilocal, i.e. if $\mathrm{supp}(F) \subset C$ then $\mathrm{supp}(TF) \subset QC$, for some fixed compact set $Q$ depending on $T$,
(4) $T$ maps the space $L^1_{loc}$ into the space of continuous functions on $G$.

Under these condition, if the operator $T$ is bounded on $L^p(G)$, the transferred operator is defined and is bounded on $L^p(X)$ with the same bound, for an arbitrary measure-preserving action. The same holds for the maximal function associated with a sequence $T_n$ of such operators. In the context of convolution operators, condition (4) would usually require the absolute continuity of the measure $\mu_t$. However, we would like to emphasize that the transfer principle is valid also for any family $\mu_t$ of singular measures on $G$ (many of which will appear below) and does not require absolute continuity.

*Remark* 6.4. Let us note further that in the transfer principle formulated in Theorem 6.2 (as compared to Theorem 5.9)

(1) The family of measures $\mu_t$ whose maximal inequalities are being transfered is arbitrary, and in particular the measures are not required to be symmetric (or, as already noted, absolutely continuous).
(2) $G$ need not be unimodular, and furthermore the Følner sequence which guarantees the validity of the transfer principle can be arbitrary and no growth conditions (such as the doubling condition) on it are assumed.

(3) A principle of local transfer, namely when all the measures $\mu_t$ on $G$ have their support contained in a fixed compact set, holds without any restriction at all on the group $G$, which need not be amenable in this case.

(4) Similar results can be easily formulated for semi-group actions. Note that then we must consider the *anti-representation* $\pi'(g)f(x) = f(gx)$ (satisfying $\pi'(hg) = \pi'(g)\pi'(h)$), and the operators $\pi'(\mu_t)f(x) = \int_G f(gx)d\mu_t$.

## 6.2. Generalizations of the doubling condition : Regular Følner sequences.

In §6.1 it was established that the transfer principle for amenable groups is a direct consequence of the existence of a Følner sequence, and hence is valid for any amenable group. Thus to obtain a maximal inequality for $\mu_t$ in a general action of an amenable groups $G$ it suffices to establish a maximal inequality for the action of $G$ by right translation. This problem turns out to be a very difficult one for many of the most natural averages $\mu_t$. So far, we have seen in §5 how to prove such a result for balls w.r.t. an admissible metric, provided they satisfy the doubling volume condition, which then implies polynomial volume growth.

Calderon's original doubling volume condition [C1] is formulated for an increasing family $N_t$ of compact symmetric neighbourhoods of the identity, which generate $G$, and satisfies $N_t N_s \subset N_{t+s}$ together with $m_G(N_{2t}) \leq C m_G(N_t)$. A group possesing such a family is necessarily unimodular, so any Haar measure can be taken. Note that for balls $B_t$ w.r.t. an invariant metric, $B_{2t} = B_t \cdot B_t$, so Calderon's condition can be written as $m_G(B_t^{-1} \cdot B_t) \leq C m_G(B_t)$.

Thus, a natural generalization of the doubling condition of [C1] (as well as the conditions considered by Pitt [P] and Cotlar [Cot]), is given by the following condition, introduced by A. Tempelman [T1].

**Definition 6.5. Regular sequences**[T1]. A sequence $N_k$ of sets of positive finite measure in an lcsc group $G$ is called regular if

$$m_G(N_k^{-1} \cdot N_k) \leq C m_G(N_k)$$

for some $C$ independent of $k$, and a left Haar measure $m_G$ on $G$.

Regular sequences have been utilized to prove the following result, proved in [T2][Ch][Be] for the unimodular case, and in [E1] for general amenable groups.(For simplicity of notation, we switch to anti-representations of $G$ here).

**Theorem 6.6. Pointwise ergodic theorem for regular Følner sequences** [T1][Ch][Be][E1]. *Assume $G$ is an amenable lcsc group, $m_G$ left Haar measure, and $N_k \subset G$ is an increasing left Følner sequence, with $\cup_{k \in \mathbb{N}} N_k = G$, satisfying $m_G(N_k^{-1} N_k) \leq C m_G(N_k)$, i.e. a regular sequence. Then*

(1) *The maximal operator* $\sup_{k \in \mathbb{N}} \left| \frac{1}{m_G(N_k)} \int_{N_k} F(gh) dm_G(g) \right|$ *satisfies the weak-type $(1,1)$ and strong $L^p$ maximal inequalities for $F \in L^p(G, m_G)$.*

(2) *The operators* $\pi(\eta_k)f(x) = \frac{1}{m_G(N_k)} \int_{g \in N_k} f(gx) dm_G(g)$ *satisfy the weak-type $(1,1)$ and strong $L^p$ maximal inequalities, in every measure-preserving action of $G$ on $(X, m)$.*

(3) *The sequence $\eta_k$ satisfies the pointwise ergodic theorem in $L^1$, for every probability measure preserving action of $G$.*

As to the proof of Theorem 6.6, we note the following. Given the transfer principle of Theorem 6.2, to obtain the maximal inequalities stated in part (2) of Theorem 6.6, it suffices to prove the maximal inequality stated in part (1) for translations. The proof of the latter is similar to that of the corresponding result for ball averages discusssed in §5, and uses a natural generalization of the covering arguments employed by Wiener and Calderon to the present context. Namely, it is shown that under the condition of regularity, given a finite set of translates $\{N_k g_i \, ; \, k \in K, i \in I\}$ which covers a given compact set $F$, it is possible to select a subcover consisting of *disjoint* translates, which still covers a fixed fraction of $F$. The weak-type $(1, 1)$ inequality then follows as in §5.

To get the full pointwise ergodic theorem stated in Theorem 6.6(3), one needs also, according to the recipe of §2.3, a mean ergodic theorem and pointwise convergence on a dense subspace. We thus note the following

**Proposition 6.7. Mean ergodic theorem for Følner sequences**.

*For every Følner sequence $F_n$ on an amenable lcsc group $G$, the normalized averages on $F_n$ satisfy the mean ergodic theorem in $L^2(X)$ in every probability preserving action. Furthermore, there exists in every $L^p(X)$, $1 \le p < \infty$ a dense subspace on which pointwise convergence holds.*

*Proof.* Looking at the proofs given in §5.2 and §5.3 of the corresponding statements, it is clear that they are valid for every Følner sequence, as already noted there. □

*Remark* 6.8.     (1) A practical criterion for the existence of a regular Følner sequence has never been found. In particular it seems unknown which connected amenable Lie groups with exponential volume growth (if any) posses such a sequence. For the latter class a pointwise ergodic theorem will be proved for certain Følner sequences in §7, using a different approach.

(2) It was shown in [L] that there exists a discrete amenable group of exponential volume growth (the "lamplighter" group) for which no regular Følner sequence exists.

6.3. **Subsequence theorems : Tempered Følner sequences.** The discussion of §6.2 does not resolve the problem of the existence of a Følner sequence which satisfies the pointwise ergodic theorem in $L^1$ in an *arbitrary* lcsc amenable group. This problem was resolved in $L^2$ using a more general condition, which was introduced by A. Shulman (see [T, Ch. 5]) for this purpose, as follows.

**Definition 6.9. Tempered sequences.** A sequence of sets of positive finite measure $N_k$ in $G$ is called tempered if for some fixed $C$

$$m_G(\tilde{N}_n^{-1} \cdot N_{n+1}) \le C m_G(N_{n+1})$$

where $\tilde{N}_n = \cup_{k \le n} N_k$.

We note that this condition is very different from the regularity condition of §6.2. Indeed, regularity is primarily a growth condition, and it is often incompatibe with the Følner property, as noted in [GE]. Furthermore, regular Følner sequences have never been found and are unlikely to exist in groups with exponential volume growth. On the other hand, temperedness is primarily an invariance condition, implying that the $(n+1)$-th set is almost invariant under (inverse) left translations by all the $n$ previously chosen sets. Thus an easy induction argument shows that starting with any Følner sequence, one can choose a tempered subsequence [L, Prop. 1.4].

It was shown by A. Shulman that the averages associated with a tempered Følner sequence satisfy the pointwise ergodic theorem in $L^2$ (see [T, §5.6]). The complete result in $L^1$ has been established by E. Lindenstrauss [L], as follows.

**Theorem 6.10. Pointwise ergodic theorem for tempered Følner sequences**
[L]. *Let $N_k$ be a tempered left Følner sequence on an amenable lcsc group $G$. The normalized averages $\eta_k f(x) = \frac{1}{m_G(N_k)} \int_{N_k} f(gx) dg$ satisfy the pointwise ergodic theorem, and the weak-type $(1,1)$ maximal inequality in $L^1$, in every probability-preserving action of $G$. Thus every amenable lcsc group admits a Følner sequence satisfying the pointwise ergodic theorem in $L^1$.*

As to the proof of Theorem 6.10, we note that it proceeds as usual by establishing the weak-type $(1,1)$ maximal inequality for the operators $\eta_k$ acting on $X = G$, and then appeals to the transfer principle of Theorem 6.2. The maximal inequality uses a covering argument in $G$, where the covering sets are taken from the set of all translates $\mathcal{F} = \{N_n g \, ; \, n \in \mathbb{N}, g \in G\}$. The proof of the covering argument introduces an important new probabilistic technique to the discussion, as follows. Given a compact set $F \subset G$ to be covered, it is shown that a probability distribution can be introduced on the set of subcollections of $\mathcal{F}$, such that typically (w.r.t. the probability distribution) a random subcollection consists of almost disjoint sets, and these cover most of $F$, approximately evenly. This allows the construction of a subcover which retains at least a fixed fraction of the measure of the original set $F$, and which is almost disjoint, and thus the weak-type $(1,1)$-maximal follows in a manner similar to the proof of Theorem 5.9.

We note that another proof of the same result is due to B. Weiss [We1]. The proof establishes a weak-type $(1,1)$-maximal inequality for averaging with respect to a temepered Følner sequence, based of a covering argument for the family of translates of the sequence. The covering Lemma is based on an interesting direct combinatorial argument utlizing temperedness, and is thus a deterministic one.

*Remark* 6.11. The process of refining a given Følner sequence to a tempered sub-sequence may be rather drastic, namely the resulting subsequence may be very sparse. Thus in [L, Cor. 5.6] it is shown that any tempered Følner sequence on the lamplighter group must satisfy $\lim_{k \to \infty} |N_{k+1}| / |N_k| = \infty$, and in particular be of super-exponential growth. A question raised in [L] is whether there always exists a Følner sequence with exponential growth satisfying the pointwise ergodic theorem in $L^1$.

## 7. A NON-COMMUTATIVE GENERALIZATION OF WIENER'S THEOREM

The present section is devoted to an exposition of some non-commutative generalizations of Birkhoff's and Wiener's pointwise ergodic theorems. These results are based on a method introduced independently by Dunford [D] and Zygmund [Z], and first utilized in their proof of a fundamental result on (dominated) pointwise convergence of averages of product type on the product of several not-necessarily-commuting one-parameter flows. We will demonstrate below that the method can be applied to yield a large collection of pointwise convergence theorems and strong maximal inequalities for families of measures on groups, provided the groups can be represented as a product of more basic flows, for example one parameter flows. The Dunford-Zygmund method is thus an ideal tool for proving ergodic theorems for connected Lie groups, when they admit global Lie coordinates of the second kind, and as we shall see also for algebraic groups, which admit a variety of decomposition theorems. In particular, the method was used by Emerson and Greenleaf [GE] to give a proof of a pointwise ergodic theorem for certain sequences of Følner averages on any connected amenable group. Below we will present the proof of a generalization of the Dunford-Zygmund Theorem, as well as a sharper form of the Greenleaf-Emerson theorem. We will also give a proof of a pointwise ergodic theorems for groups with an Iwasawa decomposition, generalizing Tempelman's theorem for connected Lie groups. Finally we will indicate further results based on these methods, which apply for example to general algebraic groups (see [N9] for details).

We note however that the averages that the Dunford-Zygmund method apply to are of a very specific form, and in particular they usually bear no resemblance to ball averages w.r.t. an invariant metric on the group.

7.1. **The Dunford-Zygmund method.** Let us start by formulating and proving the Dunford-Zygmund theorem in the most basic special case. As will become clear below, however, this case already demonstrates the main ideas involved.

**Proposition 7.1. The pointwise ergodic theorem for two non-commuting flows.** [D][Z]. *Let $u_t$, $t \in \mathbb{R}$ and $v_s$, $s \in \mathbb{R}$ be two not-necessarily-commuting $\mathbb{R}$-flows, namely representations of $\mathbb{R}$ as measure preserving transformations of a probability space $(X, m)$. Then*

(1) *The strong $L^p$-maximal inequality, $1 < p \leq \infty$, holds for rectangle averages, namely*

$$\left\| \sup_{T, S > 0} \left| \frac{1}{4TS} \int_{-T}^{T} \int_{-S}^{S} u_t v_s f \, dt \, ds \right| \right\|_p \leq C_p \|f\|_p$$

(2) *Let $\mathcal{U}$ (resp. $\mathcal{V}$) be the conditional expectation w.r.t. the $\sigma$-algebra of sets invariant under every $u_t$, $t \in \mathbb{R}$ (resp. every $v_s$, $s \in \mathbb{R}$). Then for every $f \in L^p(X)$, $1 < p < \infty$ and for almost every $x \in X$*

$$\lim_{\min(T, S) \to \infty} \frac{1}{4TS} \int_{-T}^{T} \int_{-S}^{S} u_t v_s f(x) \, dt \, ds = \mathcal{U} \mathcal{V} f(x)$$

*and the convergence is also in the $L^p$-norm.*

(3) *Pointwise convergence also holds for $f \in L(\log L)(X)$.*

*Proof.* (1) The maximal inequality stated follows from the following simple observations. As usual assume without loss of generality that $f \geq 0$, and then clearly

$$\sup_{0 < T \leq T_0, \, 0 < S \leq S_0} \frac{1}{4TS} \int_{-T}^{T} \int_{-S}^{S} u_t v_s f(x) \, dt \, ds \leq$$

$$\sup_{0 < T \leq T_0} \frac{1}{2T} \int_{-T}^{T} u_t \left( \sup_{0 < S} \frac{1}{2S} \int_{-S}^{S} v_s f(x) \, ds \right) dt \leq M_U^* \left( M_V^* f \right)(x)$$

Now the maximal function $M_V^* f$, associated with the averages $1/2S \int_{-S}^{S} v_s ds$, has an $L^p$-norm bound, by the maximal inequality for one-parameter flows, which is a consequence of Theorem 5.7. The same of course holds for the maximal function $M_U^*$ which is associated with the averages $1/2T \int_{-T}^{T} u_t dt$. Hence the maximal inequality for averaging over $|t| \leq T$, $|s| \leq S$ follows. The desired strong maximal inequality then follows from the monotone convergence theorem.

(2) According to the recipe of §2.3, given the maximal inequality, pointwise almost everywhere convergence follows provided it can be established on a dense subspace of $L^p(X)$ $(p < \infty)$. We again apply the analog of Riesz's argument [R] used in the proof of Theorem 5.3 and define the space $\mathcal{K} = \text{span}\{v_s f - f \mid s \in \mathbb{R}, f \in L^\infty(X)\}$. The sum of $\mathcal{K}$ and the space of $\{v_s \,;\, s \in \mathbb{R}\}$-invariant functions is dense in $L^p(X)$, $1 \le p < \infty$. Now if $h = (v_{s_0} f - f) \in \mathcal{K}$, then for $S \ge |s_0|$, and every $T > 0$

$$\left| \frac{1}{2T} \int_{-T}^{T} u_t \left( \frac{1}{2S} \int_{-S}^{S} v_s (v_{s_0} f - f)(x) ds \right) dt \right| \le \frac{2 |s_0| \cdot \|f\|_\infty}{S} \longrightarrow 0$$

as $S \to \infty$. Therefore

$$\lim_{\min(S,T) \to \infty} \frac{1}{2T} \int_{-T}^{T} u_t \left( \frac{1}{2S} \int_{-S}^{S} v_s (v_{s_0} f - f)(x) ds \right) dt = 0 \,.$$

Of course, if $h$ is $v_s$-invariant, then the expression

$$\frac{1}{4TS} \int_{-T}^{T} \int_{-S}^{S} u_t v_s h(x) ds dt = \frac{1}{2T} \int_{-T}^{T} u_t h(x) dt$$

converges almost everywhere by Birkhoff's theorem. Thus pointwise convergence holds on a dense subspace, and using the maximal inequality, also for every $f \in L^p(X)$, $1 < p \le \infty$.

(3) The identification of the limit is obtained as follows. For any $f \in L^p(X)$ :

$$\left\| \frac{1}{4TS} \int_{-T}^{T} \int_{-S}^{S} u_t v_s f ds dt - \mathcal{U}\mathcal{V} f \right\|_p \le$$

$$\le \left\| \frac{1}{2T} \int_{-T}^{T} u_t \left( \frac{1}{2S} \int_{-S}^{S} v_s f ds - \mathcal{V} f \right) dt \right\|_p + \left\| \frac{1}{2T} \int_{-T}^{T} u_t \mathcal{V} f dt - \mathcal{U}\mathcal{V} f \right\|_p \,.$$

Using the norm convergence to $\mathcal{V} f$ of the averages $\frac{1}{2S} \int_{-S}^{S} v_s f ds$, and the fact that each operator $\frac{1}{2T} \int_{-T}^{T} u_t dt$ is a contraction we can estimate the norm of the first summand. The norm of the second summand is estimated by the norm convergence of $\frac{1}{2T} \int_{-T}^{T} u_t h dt$ to $\mathcal{U} h$, $h = \mathcal{V} f$. We conclude that the limit of the foregoing expression is 0, and the proof of convergence in $L^p$ is complete.

(4) For the proof of pointwise convergence for $f \in L(\log L)(X)$ we refer the reader to the original argument in [Z].

<div align="right">□</div>

*Remark* 7.2. **On the identification of the limit in Theorem 7.1.**

(1) It is obvious that Theorem 7.1 and its proof admit extensive generalizations. Thus for example we can replace $U$ and $V$ by any Abelian (not necessarily connected) Lie (semi-) group, and replace the interval averages by any other Følner averages satisfying the pointwise and maximal ergodic theorem. Furthermore, similar conclusions holds for any finite sequence $U_1, \ldots, U_k$ of Abelian (semi-) groups. We will comment on this fact and some further generalizations below.

(2) Anticipating some arguments that will occur later on, note that in the proof of part (2), the asymptotic invariance of the intervals $[-S, S] \subset V$ under translation plays an essential role in the proof. It provides an estimate in the pointwise convergence argument for the bounded functions in question which is *uniform* in $T$, and depends only on the size of $S$, allowing the argument to proceed.

(3) A key problem in utilizing Theorem 7.1 is to make the identification of the possible limits in Theorem 7.1 more precise. For example, if the intersection of the ranges of the two projections $\mathcal{U}$ and $\mathcal{V}$ reduces to the constant functions, when is it true that the limit in the theorem is the projection on the constants, namely $\int_X f dm$ ? Note that $\mathcal{U}f$ is not necessarily $\mathcal{V}$-invariant if $f$ is.

*Remark* 7.3. **Unrestricted convergence and $L(\log L)$-results**

(1) It was established already in [Z] that the existence of the pointwise limit stated in part (2) of Theorem 7.1 holds for functions $f \in L(\log L)^{k-1}(X)$, in the case of $k$ one-parameter flows. We refer to [Fa] for a more general result.

(2) On the other hand, We note that it is essential for the argument given in the proof of theorem 7.1 that we consider the strong maximal inequality - given by a *norm* inequality in $L^p$. Weak-type inequalities cannot be treated by the same argument, and indeed there is no weak-type $(1, 1)$ maximal inequality for general two non-commuting flows. For general multi-parameter flows, and in fact even for commuting ones, pointwise convergence in $L^1$ does not always hold. For a counterexample (for an action of $\mathbb{Z}^d$) see [Kr, §6.1]. This phenomenon is due to the fact that in Theorem 7.1 we allow *unrestricted convergence*, meaning that in the rectangle averages the side lengths $T$ and $S$ are chosen independently. Such general rectangles behave differently than squares - for which pointwise convergence in $L^1$ does hold for commuting flows, by Wiener's theorem. The situation here is analogous to the classical discussion of maximal inequalities for unrestricted rectangle averages on $\mathbb{R}^d$, and the maximal inequality of Theorem 7.1 can be viewed as a non-commutative generalization of the Jessen-Marcienkiewicz-Zygmund maximal inequality - see e.g [Fa] or [S2, Ch. X §2.2] for a discussion.

Let us make the following simple observation, which will be found useful below.

*Remark* 7.4. Under the condition of Theorem 7.1, if at least one of the groups $U$ or $V$ is ergodic on $(X, m)$ then the following pointwise ergodic theorem holds, for $f \in L^p(X)$, $1 < p < \infty$, and almost every $x \in X$ :

$$\lim_{\min(T,S)\to\infty} \frac{1}{4TS} \int_{-T}^{T} \int_{-S}^{S} u_t v_s f(x) dt ds = \int_X f dm \,.$$

Indeed, the assumption amounts to the fact that $\mathcal{V}$ or $\mathcal{U}$ coincide with the projection $f \mapsto \int_X f dm$. Since in any case $\int_X \mathcal{U}f dm = \int_X \mathcal{V}f dm = \int_X f dm$, the conclusion follows.

We now note the following natural generalization of Theorem 7.1, as well as [GE, Thm. 7.1], which will be used below.

**Theorem 7.5.** [N9] *Let $U$ and $V$ be lcsc group, and assume that $V$ is amenable. Let $E_T \subset U$ and $F_S \subset V$ be compact sets of finite Haar measure. Assume that $F_S$ is a Følner family, and that both families of normalized averages satisfy the pointwise ergodic theorem and strong maximal inequality in $L^p$, $1 < p < \infty$ . Let $U$ and $V$ act by measure-preserving transformations on a probability space $(X, m)$. Then*

(1) *The strong $L^p$-maximal inequality, $1 < p \leq \infty$, holds for the averages given by*

$$\left\| \sup_{T,S>0} \left| \frac{1}{m_U(E_T)m_V(F_S)} \int_{E_T} \int_{F_S} uvf \, dm_U(u) dm_V(v) \right| \right\|_p \leq C_p \|f\|_p$$

(2) *Let $\mathcal{U}$ (resp. $\mathcal{V}$) be the conditional expectation w.r.t. the $\sigma$-algebra of sets invariant under every $u \in U$, (resp. every $v \in V$). Then for every $f \in L^p(X)$, $1 < p < \infty$ and for almost every $x \in X$*

$$\lim_{\min(T,S)\to\infty} \frac{1}{m_U(E_T)m_V(F_S)} \int_{E_T} \int_{F_S} uvf(x) \, dm_U(u) dm_V(v) = \mathcal{U}\mathcal{V}f(x)$$

*and the convergence is also in the $L^p$-norm.*

(3) *If either $U$ or $V$ act ergodically, then the limit satisfies $\mathcal{U}\mathcal{V}f(x) = \int_X f \, dm$.*

(4) *If both $E_T$ and $F_S$ satisfy the weak-type maximal inequality in $L^1$, then pointwise convergence holds also for $f \in L(\log L)(X)$.*

*Remark* 7.6. The proof of Theorem 7.5 proceeds using an argument similar to the one given in the proof of Theorem 7.1 (and Remark 7.4). We note however that the Følner property of $V$ plays a crucial role in the pointwise result. If we are interested only in norm convergence, then the assumptions that $V$ is amenable and $F_S$ are Følner are not necessary, and $F_S$ can be any family of sets on $V$ for which the mean ergodic theorem holds. This fact follows easily using the argument in the proof of part (3) of Theorem 7.1.

We also note that if both families of averages $E_T$ and $F_S$ satisfy the weak-type maximal inequality in $L^1$, then the combined averages satisfy the maximal inequality in $L(\log L)(X)$. This is a general property of the composition of two maximal operators both satisfying the weak-type maximal inequality, which is due to [Fa]. We refer to [N9] for the details.

7.2. **The ergodic theory of semi-direct products.** Continuing now with our theme of establishing pointwise ergodic theorems for group actions, assume now that the two groups $U$ and $V$ that figure in Theorem 7.5 generate an lcsc group. Thus we let $G = UV$ be a semi-direct product, where $V \lhd G$ is a closed normal subgroup, and $U$ a closed subgroup of $G$, so that $G \cong U \ltimes V$. The multiplication in $U \ltimes V$ is given by $(u_1, v_1)(u, v) = (u_1 u, u^{-1} v_1 u \cdot v)$, and the explicit isomorphism between $U \ltimes V$ and $G$ is given by $\psi : (u, v) \mapsto uv$. The left Haar measure on $U \ltimes V$ is given by the following well-known recipe.

**Lemma 7.7.** *The map $\psi : U \ltimes V \to G$ is homeomorphism that maps the product of the two left Haar measures $m_U \times m_V$ on the product space $U \times V$ to a left Haar measure on $G \cong U \ltimes V$.*

*Proof.* The claim amount to the fact that for $f \in C_c(G)$, writing $f(uv) = F(u, v)$, the integral given by

$$\int_{G=UV} f(uv) d\psi_*(m_U \times m_V)(uv) = \int_{U \ltimes V} F(u, v) dm_U(u) dm_V(v)$$

is invariant under left translations. But since $V$ is normalized by $U$, and $m_U$, $m_V$ are left-invariant, translating by $g = u_1 v_1$ we get :

$$\int_{UV} f((u_1 v_1 uv)) d\psi_*(m_U \times m_V) = \int_U \left( \int_V F(u_1 u, u^{-1} v_1 u \cdot v) dm_V(v) \right) dm_U(u) =$$

$$= \int_U \left( \int_V F(u_1 u, v) dm_V(v) \right) dm_U(u) =$$

$$= \int_{U \times V} F(u, v) dm_U dm_V = \int_{UV} f(uv) d\psi_*(m_U \times m_V)$$

$\square$

Lemma 7.7 allows us to state a pointwise ergodic theorem for averages defined by certain sets of product type on $G = UV$, taken with the normalized left Haar measure as follows.

**Theorem 7.8. Pointwise ergodic theorem for semi-direct products with normal amenable subgroup**[N9]. *Let $G = UV$ be a semidirect product as above, where $V \triangleleft G$ is an amenable subgroup. Let $E_T \subset U$ and $F_S \subset V$ be as in Theorem 7.5. Assume that $G$ acts ergodically by measure-preserving transformations on a probability space $(X, m)$. Then for every $f \in L^p(X)$, $1 < p < \infty$ and almost every $x \in X$ we have*

$$\lim_{\min(T,S) \to \infty} \frac{1}{m_U(E_T) m_V(F_S)} \int_{u \in E_T} \int_{v \in F_S} uvf(x) dm_U(u) dm_V(v) = \int_X f dm =$$

$$= \lim_{\min(T,S) \to \infty} \frac{1}{m_G(Q_{T,S})} \int_{Q_{T,S}} gf(x) dm_G(g)$$

*where the convergence is also in the $L^p$-norm. Here $Q_{T,S} \subset G$ is the image of $E_T \times F_S \subset U \times V$ under the multiplication map $\psi$, and the average is w.r.t. the restriction of left-invariant Haar measure on $G$ to $Q_{T,S}$.*

*Furthermore, if the families $E_T$ and $F_S$ satisfy the weak-type maximal inequality in $L^1$, then the averages on $Q_{T,S}$ satisfy the maximal inequality in $L(\log L)$.*

*Proof.* The fact that the pointwise (and the norm) limit exist follows from Theorem 7.5, and the limit is given by $\mathcal{U}\mathcal{V}f$. The only remaining issue is to identify the limit. But since $V$ is a normal subgroup of $G$, for every function $h$ invariant under $V$, its translate $\pi(g)h$ is also invariant under $V$. Therefore the space $\mathcal{I}_V$ of $V$-invariant functions is $G$-invariant, and hence it is also invariant under the projection $\mathcal{U}$, which is the limit in the strong operator topology of $\frac{1}{m_U(E_T)} \int_{E_T} u dm_U(u)$. Hence $\mathcal{U}$ maps the space $\mathcal{I}_V$ into the space $\mathcal{I}_V \cap \mathcal{I}_U$ of functions invariant under both $U$ and $V$. Since $G$ is ergodic, the latter space consists of the constant functions only, and hence the function $\mathcal{U}\mathcal{V}f \in \mathcal{I}_U \cap \mathcal{I}_V$ is a constant. Its value must be $\int_X f dm$ since the operators $\pi(u)$ and $\pi(v)$ are measure-preserving, and $\mathcal{U}$ and $\mathcal{V}$ are conditional expectations. The assertion regarding $Q_{T,S}$ follows from Lemma 7.7. For the last assertion stated, see Remark 7.6. $\square$

We note that Theorem 7.8 generalizes [GE, Thm. 7.1], which considers the case where $G$ is amenable, and identifies the limit only when $Q_{T_n, S_n}$ is assumed to be a Følner sequence. In the case of amenable semi-direct products it is possible to use Følner families on the constituent groups to construct explicitly a Følner family on $G$, a fact due to F. Greenleaf [Gr, Thm 5.3], whose formulation follows.

**Theorem 7.9. Construction of Følner sequence in semi-direct products**[Gr, Thm. 5.3] *Assume $G$ is an amenable lcsc group, and $G = UV$ is a semi-dierct product, where $U$ and $V$ are closed subgroups with $V$ normal. Let $E_T \subset U$ and $F_S \subset V$ be Følner families. Then*

(1) *There exist subsequences $T_n$ and $S_n$ such that the sets $J_n = Q_{T_n,S_n} \subset G$ associated with $E_{T_n}$ and $F_{S_n}$ constitute a Følner sequence in $G$.*

(2) *$J_n$ can be chosen to be an increasing sequence whose union covers $G$, provided $E_T$ and $F_S$ have the same properties.*

7.3. **Structure theorems and ergodic theorems for amenable groups.** Theorem 7.8 and Theorem 7.9 can be exploited together with some structure theorems for various classes of lcsc groups to produce a variety of pointwise and maximal ergodic theorems. The first such result was established by F. Greenleaf and W. Emerson for connected amenable lcsc groups in [GE, Thm. 3.1], and was later extended to all connected lcsc groups by A. Tempelman [T, Ch. 6, Thm 8.4]). We will give a streamlined account of the connections between the structure theory of some classes of lcsc groups and the ergodic theorems they satisfy, and this will allow us to derive extensions of these results in several directions. Succinctly put, we can establish pointwise ergodic theorems for groups of the form $G = KP$, where $K$ is compact and $S$ solvable and normal, and also more generally for groups $G = KP$ where $K$ is compact and $P$ amenable, but not necessarily normal.

To begin with, we have the following corollary of Theorem 7.8.

**Theorem 7.10. Pointwise ergodic theorem for algebraically connected amenable algebraic groups and connected amenable Lie groups** [N9].

*Let $G$ be an amenable lcsc group, and $H = UV$ be any decomposition into a semi-direct product of two closed subgroup, where $V$ is normal. Assume that (keeping the notation of Theorem 7.8) $E_T \subset U$ and $F_S \subset V$ satisfy the pointwise ergodic theorem in $L^p$, $1 < p < \infty$ and that $F_S$ is Følner.*

(1) *The family $Q_{T,S} \subset H$ also satisfies the pointwise ergodic theorem in $L^p$, $1 < p < \infty$, as $\min\{T,S\} \to \infty$ independently, namely in every ergodic $H$-space*

$$\lim_{\min\{T,S\} \to \infty} \frac{1}{m_G(Q_{T,S})} \int_{Q_{T,S}} gf(x)dm = \int_X f\,dm\,.$$

(2) *Every algebraically connected amenable algebraic group $G$ over a locally compact non-discrete field admits a closed normal co-compact subgroup $H$ with a semi-direct product structure $H = UV$, with $E_T \subset U$ and $F_S \subset V$ satisfying the conclusion in (1). For a certain compact subgroup $K \subset G$, the normalized averages on the compact sets (of positive Haar measure in $G$) $Q'_{T,S} = KE_T F_S$, satisfy the conclusion in (1).*

(3) *The same conclusion holds for a connected amenable Lie group $G$, provided its maximal compact normal subgroup $C$ is trivial. In general, the conclusion holds for the averages $Q''_{T,S} \subset G$ which are the inverse images of $Q'_{T,S} \subset G/C$.*

*Proof.* The convergence statement in Part (1) is of course an immediate consequence of Theorem 7.8.

Parts (2) and (3) state, first, that the groups $G$ in the relevant category satisfy a structure theorem. For algebraic groups over $F$, this fact follows without difficulty from the fact that $G$ has homomorphic image with finite kernel, which has a finite index subgroup isomorphic to a subgroup $L$ of a proper parabolic subgroup of $GL(n,F)$ (see e.g. [Gu, Thm. IV.2]. For $L$ a decomposition of the form $L = KAN$ is valid, where $N$ is the unipotent radical (which is nilpotent), $A$ is Abelian (and thus $S = AN$ is solvable), and $K$ compact. The existence of $E_T$ and $F_S$ with the required properties on the component groups (or simply on $G$ itself) is thus clear as soon as the averages $F_S$ are shown to exist on unipotent groups. This can be done by induction on the dimension, for instance.

As to connected Lie groups, assume that the maximal compact normal subgroup is trivial. Then $G$ is well-known to be a compact extension of a solvable Lie group, i.e. $G = KS$, $K$ compact, $S$ closed normal solvable and co-compact (see e.g. [Gu, Thm IV.3], or [GE]). The existence of a semi-direct product decomposition for $S$ with $E_T$ and $F_S$ as required follows, e.g. as in [GE]

To complete the proof of (2) and (3), assume that the averages constructed from $E_T$ and $F_S$ on a co-compact normal subgroup $H = UV$ converge pointwise to an $H$-invariant function. If $G = KH$, $K$ compact, then we claim that the averages $Q'_{T,S}$ on $G$ converge pointwise to a $G$-invariant function. This follows since $H$ is normal, and hence its space of invariants is invariant under $G$, and hence under $K$. Thus if $f$ is an $H$-invariant limit of the averages on $H$, the average of its translates by $K$ is still an $H$-invariant function. THus $f$ is a $G$-invariant function, hence the constant $\int_G f \, dm$ when $G$ acts ergodically.

Finally, to complete the case of connected Lie groups, note that clearly any limit of the averages $Q''_{T,S} = m_C * Q''_{T,S} * m_C$ is a $C$-invariant function. We therefore consider the space of $C$-invariant functions, which is a $G$-invariant subspace on which $G$ acts via $G/C$. The preceding argument applies, and any limit of $Q''_{T,S}$ is indeed invariant under $G$.    □

*Remark* 7.11.      (1) We note that by Theorem 7.9 a subsequence of sets $Q'_{T_n, S_n}$ can be chosen which is an increasing sequence of Følner sets $J_n$ in $G$, whose union covers $G$, and which satisfies the pointwise ergodic theorem in $L^p$, $1 < p < \infty$. Thus Theorem 7.10 generalizes the pointwise ergodic theorem for connected amenable Lie groups due to F. Greenleaf and W. Emerson [GE].

(2) In the proof given in [GE, Thm. 3.1], the identification of the limit is achieved only by restricting the averages, and choosing sequences $T_n \to \infty$ and $S_n \to \infty$ where the averages $Q_{T_n, S_n}$ can be guaranteed to form a Følner sequence $J_n$. Then the strong limit of this sequence of averages must be projection onto $G$-invariant function (by Proposition 6.7), hence the constant $\int_X f \, dm$ in the ergodic case. The pointwise ergodic theorem is asserted in [GE, Thm. 3.1] only for averaging along such Følner sequences in connected amenable Lie groups.

(3) Nevertheless, in fact pointwise convergence to the ergodic mean holds more generally, for the unrestricted averages on $Q_{T,S}$, as Theorem 7.10 show. Furthermore the assumption of amenability of $G$ is also superfluous, as Theorem 7.8 shows. Indeed, as we saw the identification of the limit as a product of conditional expectations that the Dunford-Zygmund method provides, can replace the assumption of the existence of global Følner sets in $G$. Together with the existence of a normal amenable subgroup, this allows the conclusion that the limit must be the ergodic mean.

(4) Note also that the *existence problem* for such a Følner sequence in an *arbitrary* lcsc amenable group (not necessarily connected Lie) has been solved by the pointwise ergodic theorem for tempered Følner sequences, which holds in fact in $L^1$ (and not only $L^p$, $1 < p < \infty$), as stated in Theorem 6.10.

The Dunford-Zygmud method provides a great deal of useful information even in the case where the group $G$ in question does not have any semi-direct product structure, for example if $G$ is a simple group. We now turn to discuss this possibility.

7.4. **Structure theorems and ergodic theorems for non-amenable groups.** Let us first recall the following well-known result regarding Haar measure (see e.g. [Gaa, Ch. V, §3, Prop. 12]), which generalizes Lemma 7.7. We remark that for Lie groups (and by the same argument for algebraic groups over locally compact non-discrete fields) a very simple proof is given e.g. in [GV, Ch. 2, §2.4], see also [Kn, Ch. V §6].

**Lemma 7.12. Haar measure on general products.**
    Let $G$ be an lcsc group, and let $P$ and $K$ be two closed subgroups, such that $P \cap K$ is compact. Assume that $G = PK$. Then for every $f \in C_c(G)$ we have

(1)
$$\int_G f(g) dm_G(g) = \int_{P \times K} f(pk) \frac{\Delta_G(k)}{\Delta_K(k)} dm_P(p) dm_K(k)$$

where $m_G$, $m_P$ and $m_K$ are left-invariant Haar measures, and $\Delta_K$ and $\Delta_G$ are the modular functions of $K$ and $G$.

(2) In particular, if $G$ is unimodular, then

$$\int_G f(g) dm_G(g) = \int_{P \times K} f(pk) dm_P(p) \Delta_K(k^{-1}) dm_K(k) = \int_{P \times K} f(pk) dm_P(p) d\eta_K(k)$$

where $d\eta_K = \Delta_K(k^{-1}) dm_K$ is right-invariant Haar measure on $K$.

(3) If in addition $K$ is unimodular, then $dm_G = dm_P dm_K$.

Perhaps the simplest family of groups which are not semi-direct products and for which Lemma 7.12 applies is the family of groups $G$ with an Iwasawa decomposition. Namely $G$ contains two closed subgroup $P$ and $K$, with $P$ amenable and $K$ compact, such that $G = PK$, and neither subgroups is normal. This family however is extremely important, since as is well known, it contains all connected lcsc groups and all algebraically connected algebraic groups over locally compact non discrete fields.

Together with Theorem 7.8, Lemma 7.12 can be used to derive a variety of pointwise and maximal ergodic theorems for averaging on compact sets with respect to Haar-uniform measure on them. Let us start with the following result, which utilizes also Theorem 7.10. We note that the third and fourth part of the following theorem generalize results of A. Tempelman[T, Ch. 6, Thm 8.6, 8.7].

**Theorem 7.13. Pointwise ergodic theorem for groups with an Iwasawa decompostion**[N9]. Let $G$ be a non-amenable lcsc group with an Iwasawa decomposition, $G = PK$, with $P$ amenable and closed, and $K$ compact. Let $E_T \subset P$ be a family of sets, giving rise to averages (w.r.t. left Haar measure on $P$) satisfying the maximal and pointwise ergodic theorem in $L^p$, $1 < p < \infty$. Then

(1) The Haar uniform averages on the sets $R_T = E_T K \subset G$, (normalized by left Haar measure on $G$) satisfy the maximal inequality in $L^p$, $1 < p < \infty$.

(2) In every probability $G$-space in which $P$ acts ergodically, the averages in (1) converge pointwise to the ergodic mean, for every $f \in L^p$, $1 < p < \infty$.

(3) *Every algebraically connected algebraic group over a locally compact non-discrete field, and every connected lcsc group, admit an Iwasawa decomposition $G = PK$, and a family $E_T \subset P$ such that in every ergodic action of $G$ and for every $f \in L^p$, $1 < p < \infty$,*

$$\lim_{T \to \infty} \frac{1}{m_G(R_T)} \int_{R_T} gf(x) dm_G(g) = \int_X f dm \text{ for almost all } x \in X .$$

(4) *Under the same assumption as in (3), pointwise convergence and the weak-type $(1,1)$-maximal inequality hold in fact for $f \in L^1(X)$, provided that the averages on $E_T \subset P$ satisfy the same properties. Such averages $E_T$ can be chosen in every algebraically connected semisimple algebraic group.*

*Sketch of proof.* Parts (1) and (2) are immediate consequences of Lemma 7.12 and Theorem 7.5. (In fact Theorem 7.5 is superfluous in the present simple case, which can be deduced directly)

Part (3) involves three ingredients. First, every algebraically connected semisimple algebraic group has an Iwasawa decomposition, and in particular, so does every connected semisimple Lie group with finite center. A general algebraic group is a semidirect product of a semisimple component $M$ without compact factors, and an amenable (and algebraic) radical. Taking the product of the amenable radical and the minimal parabolic subgroup of the semisimple component $M$ we clearly obtain an Iwasawa decomposition. Similarly, every connected non-amenable lcsc group $G$ has a compact normal subgroup $K_0$ such that $G/K_0$ is a connected semisimple Lie group without compact factors and with trivial center. Clearly the inverse image of a minimal parabolic subgroup of $G/K_0$ again gives rise to an Iwasawa decomposition in $G$.

The second ingredient is the fact that the subgroup $P$ figuring in the Iwasawa decomposition described above does admit sets $E_T$ which satisfy the pointwise and maximal ergodic theorem in $L^p$, $1 < p < \infty$. This of course follows in both cases considered here from Theorem 7.10.

The third ingredient is the fact that $P$ in the Iwasawa decomposition above does indeed always act ergodically if $G$ does. This fact follows from the Howe-Moore mixing theorem [HM] (see also [HT]), as follows. The vanishing at infinity of matrix coefficients of unitary represenantions without invariant vectors of a simple algebraic group $H$ implies that every element whose powers are not confined to a compact subgroup of $H$ acts ergodically in ergodic $H$-spaces. This can be used to show that $P$ as defined above is always ergodic in ergodic $G$-spaces, in both cases. Indeed, every $P$-invariant function is invariant under the solvable radical $S$ of $P$. The latter is a normal subgroup of $G$, and so $G$ leaves invariant the space of $S$-invariant functions. On the latter space $G$ acts vis its factor group $G/S$, which is a semisimple group. Factoring further by the maximal compact normal factor of $G/S$ we are thus reduced to the case of a semisimple group without compact factors, where the Howe-Moore mixing theorem can be applied.

As to Part (4), we are considering the composition of the averages associated with $E_T \subset P$ (taken with normalized left-invariant Haar measure on $P$), composed with the contant bounded operator given by averaging on the compact set $K$. It is clear by definition that if the maximal function associated with averaging on $E_T$ satisfies the weak-type $(1,1)$ maximal inequality in its own right, then it will still have this property when composed with a fixed bounded operator. Now for algebraically connected algebraic groups, $G = PK = ANK$, where $A$ is a split torus. It follows that $A$ admits a family of Følner sets with a pointwise and weak-type maximal theorem in $L^1$, and thus we can compose the averages on $A$ with a fixed average on a compact set in $N$ and averaging on $K$. The integration functional $f \mapsto \int_X f \, dm$ is invariant under such operators, and so the limit is still the ergodic mean. We refer to [N9] for more details.                                           □

**Example 7.14. Groups of graph automorphisms with an Iwasawa decomposition.** We note that the class of Iwasawa groups is very extensive, and contains many lcsc groups which are not algebraic. For example, consider a closed non-compact boundary-transitive subgroup $G$ of the group of automorphisms of a bi-regular tree, or more generally of the automorphism group of a locally finite graphs with infinitely many ends. The stability group $P$ of a point in the boundary is a closed amenable subgroup, and it has a compact complemet $K$ with $G = KP$. It is not hard to show that $P$ itself has the structure of a semi-direct product $P = \mathbb{Z} \ltimes K_P$, where $K_P$ is compact (see [N0] for the details). Hence $P$ fulfills the hypotheses of Theorem 7.10 and so has Følner sets $E_T$ satisfying the pointwise and maximal theorems in $L^p$, $1 < p < \infty$. By Theorem 7.13, $R_T = E_T K$ satisfies the maximal and pointwise theorem in $G$-actions in which $P$ acts ergodically. However, by [LM], typically for such groups $G$ the subgroup $P$ is indeed ergodic, and even mixing in every ergodic action of $G$, so the conclusion of Part (3) of Theorem 7.13 holds in the context as well. A similar conclusion applies in other contexts as well, e.g. for certain groups of automorphisms of a product of bi-regular trees, which act transitively on the product of the boundaries of the trees.

*Remark* 7.15. We note that it is possible to use Theorem 7.8 to deduce a completely different pointwise and maximal ergodic theorems for the algebraic groups in question. Indeed, represent an algebraic group $G = ML$ as a semi-direct product of a semisimple algebraic group without compact factors and an amenable radical $L$. Take the sets $F_S \subset L$ provided by Theorem 7.10. If we can find $E_T \subset M$ which satisfy the maximal and pointwise ergodic theoem in $L^p$, $1 < p < \infty$, then by Theorem 7.8, the same will hold for the sets $Q_{T,S} \subset G$. Now, as we will see below, in every semisimple algebraic group without compact factors, the natural ball averages on $G$, bi-invariant under a maximal compact subgroup, do indeed satisfy the pointwise and maximal theorems. Thus another pointwise ergodic theorem is obtained from Theorem 7.8, different than Theorem 7.13. We refer to [N9] for more details.

*Remark* 7.16. The result stated in Theorem 7.13 can be greatly improved in many cases. For example, consider simple algebraic groups over a locally compact non-discrete field. Then for many ergodic actions of such a group (in fact for all actions if its split rank is at least two), it is possible to establish the following property, which is in striking contrast to the ergodic theorems we have been considering thus far. The convergence of the horospherical averages (see §12.2) associated with an Iwasawa decomposition to the ergodic mean takes place at an exponentially fast rate for almost every point. Furthermore this rate can be described explicitly, and (for groups of split rank at least two) is independent of the action altogether, namely it depends only on the group and the averages chosen. We will discuss this phenomenon in greater detail below in Theorem 12.6, and we refer to [N4] and [N9] for full details.

7.5. **Groups of bounded generation.** As noted in Remark 7.2, the Dunford-Zygmund theorem applies to any finite sequence of (say) one parameter flows. Thus the Dunford-Zygmund method produces pointwise and maximal theorems for averages on subsets of Lie groups which admit global coordinates of the second kind. A natural extension of this concept here is that $G$ can be parametrized by a product of (say) Abelian closed subgroups (none of which need to be normal or connected). Thus the method applies to an extensive class of lcsc groups, namely all lcsc groups of bounded generation, i.e. allowing global coordinates with coordinate subgroups isomorphic to $\mathbb{R}$ or $\mathbb{Z}$. By the same token, we could allow any finite sequence of amenable subgroups to be taken here, provided we take Følner averages satisfying the maximal and pointwise ergodic theorem.

It should be noted however, that in the present generality, the averages that are constructed by the Dunford-Zygmund method are usually no longer Haar-uniform averages on compact subsets with an explicit geometric or algebraic description. Giving up both the assumption of a *semi-direct* product decomposition, as well as unimodularity of the component groups, implies that the averages defined in terms of global coordinates will involve certain densities w.r.t. Haar measure, and will be supporterd on certain compact sets arising from the product. Both the sets and the densities may be difficult to compute, in general.

To illustrate the point, let us formulate the following consequence of the Dunford-Zygmund method for subsets of lcsc groups admitting global coordinates of the second kind (compare [T, Ch. 6, Thm 8.4]).

**Theorem 7.17. Pointwise ergodic theorem for boundedly generated lcsc groups**. *Let $G$ be an lcsc group, and suppose that $F \subset G$ contains $n$ closed amenable subgroups $U_i$, $1 \leq i \leq n$, (for example, isomorphic to $\mathbb{Z}$, $\mathbb{R}$, $\mathbb{Q}_p$ or $\mathbb{Q}_p^*$) such that the map $\psi : U_1 \times \cdots, \times U_n \to F$ given by $\psi(u_n, \ldots, u_1) = u_n \cdots u_1$ is surjective. Let $(X, \mathcal{B}, m)$ be a standard Borel probability $G$-space. Let $E_{R_i}^i \subset U_i$ be Følner sets satisfying the pointwise and maximal ergodic theorem in $L^p$, $1 < p < \infty$. Then as $\min \{R_i, 1 \leq i \leq n\} \to \infty$, for every $f \in L^p$, $1 < p < \infty$*

$$T(R_n, \ldots, R_1)f(x) =$$

$$= \frac{\int_{E_{R_n}^n} \cdots \int_{E_{R_1}^1} \pi(u_n \cdots u_1) f(x) dm_{U_n}(u_n) \cdots dm_{U_1}(u_1)}{m_{U_n}\left(E_{R_n}^n\right) \cdots m_{U_1}\left(E_{R_1}^1\right)} \longrightarrow \mathcal{U}_n \cdots \mathcal{U}_1 f(x)$$

*for $m$-almost all $x \in X$, and in the $L^p$-norm. The corresponding maximal function satisfies a strong maximal inequality in $L^p$, $1 < p < \infty$. Furthermore, pointwise convergence also holds for $f \in L (\log L)^{n-1} (X)$, provided the averages on $E_{R_i}^i$ satisfy the weak-type maximal inequality in $L^1$ for $1 \leq i \leq n$.*

We note that discrete groups of bounded generation, where each component group is cyclic, have attracted quite a bit of attention. It was shown by [CK] that for any ring of integers $\mathcal{O}$ of an algebraic number field, the group $SL_n(\mathcal{O})$ is of bounded generation, provided $n \geq 3$. In fact, one can take all the cyclic subgroups to be generated by elementary matrices, whose number is estimated by an explicit function of $n$ (for an estimate in the case of $SL_n(\mathbb{Z})$ see [AdM]). The problem of determining which arithmetic groups are boundedly generated has been shown to be closely connected to the congruence subgroup problem (see [Rap] and [Lu]). It was established in [Tav] that it is typically the case that $S$-arithmetic lattices in connected simply-connected absolutely simple algebraic groups defined over a number field, which have split rank at least two have the property of bounded generation.

Thus Theorem 7.17 gives the corollary that on all of these lattices the averages described in Theorem 7.17 satisfy a maximal inequality and converge pointwise. If at least one of the cyclic subgroups act ergodically, then the limit in Theorem 7.17 is $\int_X f \, dm$. This is the case for example if the action of the lattice $\Gamma \subset G$ is a restriction to $\Gamma$ of an ergodic action of $G$, as follows from the Howe-Moore mixing theorem.

*Remark* 7.18. Let us note that the following question remains unresolved by the discussion of the present chapter. Let $G$ be an lcsc group with the structure indicated in Theorem 7.17 Does the family of operators $T(R_1, \ldots, R_n)$ of Theorem 7.17 satisfy the unrestricted pointwise ergodic theorem in $L^p$, $1 < p < \infty$, namely

$$\lim_{\min(R_i) \to \infty} T(R_1, \ldots, R_n) f(x) \to \int_X f \, dm$$

pointwise almost everywhere and in $L^p$-norm, in every ergodic action of $G$ ? We remark that in the original Dunford-Zygmund formulation, the flows need not satisfy any relation, and thus no further information can be expected on the resulting limit operator. Here we assume that the component groups in question at least all lie in one and the same lcsc group. This does indeed imply that the limit operator is invariant under the group in favorable cases, but the general case remains unresolved. This problem is unresolved even for the case of connected amenable Lie groups which admit an iterated semi-direct product structure of the form $U_1 \ltimes (U_2 \ltimes (\cdots (\ltimes U_n)))$, see [GE]. It constitutes an interesting challenge even for simply-connected nilpotent Lie groups, particularly since the Emerson-Greenleaf theorem can be used in this context as an important ingredient in the proof by M. Ratner [Ra] of strict measure rigidity for actions of unipotent subgroups of solvable groups.

Let us also note that for a connected Lie group $G$, the averages denoted above by $T(R_n, \ldots, R_1)$ in Theorem 7.17 and in Theorem 7.5 are usually very different than balls w.r.t. an invariant Riemannian metric. Whether there exists some choices of the parameters $R_i$ which will give a family of averages which are comparable to balls w.r.t. an invariant Riemannian metric does not seem to be known, in general.

### 7.6. **From amenable to non-amenable groups : Some open problems.**
We have focused in the present Section on pointwise ergodic theorems for amenable and non-amenable lcsc groups which have a common origin, namely the Dunford-Zygmund method. In the succeeding chapters we will return to the theme of establishing pointwise and maximal theorems for radial (and other geometric) averages for group actions. In contrast to Chapter 5 which considered groups of polynomial volume growth, our emphasis below will be on non-amenable groups (and thus with exponential volume growth), and particularly on connected semi-simple Lie group, and more generally semi-simple algebraic groups. Before continuing with the pointwise theory of radial averages, however, let us note the following fundamental problems related to mean ergodic theory and equidistribution, as well as certain pointwise convergence problems, all of which are unresolved.

(1) **The Mean Ergodic Theorem**.
Given any strongly continuous isometric representation $\tau : G \to Iso(V)$ to the isometry group of a Banach space $V$, we can consider the operator $\tau(\mu) = \int_G \tau(g) d\mu(g) \in \text{End}(V)$, which constitutes a convex average of the isometric operators $\tau(g)$, $g \in \text{supp}\,\mu$. If there exists a projection $\mathcal{E} : V \to V^I$ where $V^I$ is the closed subspace of $\tau(G)$-invariant vectors, then a result establishing that $\|\tau(\mu_t)f - \mathcal{E}f\| \to 0$ as $t \to \infty$ is called a mean ergodic theorem for the family $\mu_t$ in the representation $\tau$. The study of such results for averages on $G = \mathbb{R}$ or $G = \mathbb{Z}$ is an extensive field, which forms one of the classical themes of operator ergodic theory and fixed point theory. It was also to some extent pursued for Følner averages on amenable groups, see e.g. [Ol]. We note however that this problem is largely unresolved for many of the natural averages $\mu_t$ on groups $G$ with exponential volume growth, even in the case of unitary representation in Hilbert spaces. For example, even for connected amenable Lie groups with exponential volume growth and the ball averages w.r.t. an invariant Riemmanian metric this problem is completely open.

(2) **The equidistribution problem**.

When $X$ is a compact metric space and $G$ acts continuously, it is of course natural to ask under what conditions can we conclude that equidistribution holds for *all* orbits w.r.t. the invariant measure. Namely, when does

$$\lim_{t \to \infty} \pi(\mu_t)f(x) = \int_X f \, dm$$

for every continuous function $f \in C(X)$ and *every* $x \in X$ hold ? When $\mu_t$ are Følner averages on an amenable group equidistribution holds if and only if $m$ is the unique $G$-invariant probability measure on $X$. This fact, which can be proved in the same way as in the classical case of $\mathbb{Z}$, clearly accounts for at least some of the popularity that Følner averages enjoy as an averaging method along the orbits. Results of this type have not been established, or disproved, for any non-amenable group, for any family of averages $\mu_t$. The main source of the difficulty lies of course in the fact that taking the weak* limit of a subsequence of the measures $\mu_t * \delta_x$, it is usually not possible to show that the limiting measure is invariant under the group. In the special case of homogeneous spaces of connected Lie groups, an extensive theory of equidistribution has been developed. For a recent comprehensive discussion of equidistribution of orbits of non-amenable groups acting on homogeneous algebraic varieties we refer to [GW].

(3) **Amenable groups of exponential volume growth**.

The ball averaging problem is completely open when the group $G$ is a connected amenable Lie group of exponential volume growth. Namely, it is even unknown whether $\beta_t$, defined w.r.t. an invariant Riemannian metric, converges pointwise in $L^2$ (or in any $L^p$). We note that since ball averages do not form a regular family, do not have the Følner property, and in general are not comparable to the product averages of Theorem 7.17, none of the methods discussed thus far applies. In fact as noted in (1) above even the mean ergodic theorem has not been established.

(4) **Pointwise convergence : The $L^1$-problem**.

It is unknown if for a left-invariant Riemannian metric on *any* connected Lie group of exponential volume growth, the normalized ball averages $\beta_t$ satisfy the pointwise ergodic theorem in $L^1$. This has not been established even in one case, as far as we know. Furthermore, in Theorem 7.10 the averages to which the Dunford-Zygmund method applies are shown to con­verge only for $f \in L(\log L)$. Thus in the class of exponential solvable Lie groups for example, the only pointwise ergodic theorem in $L^1$ established so far is for a tempered sequence of Følner averages - see §6.3.

## 8. SPHERICAL AVERAGES

In the following three chapters, we will again concentrate on radial analysis, namely consider averages on balls associated with an invariant metric on $G$. How­ever, we will now concentrate on the case when $G$ has exponential volume growth, and our prime examples will be non-compact semisimple Lie and algebraic groups. These groups being non-amenable, they do not admit an asymptotically invariant sequence and no transfer principle has been established for them. The volume of balls obviously does not satisfy the doubling property, and so the Wiener covering argument, even for convolutions, does not apply. Thus none of the arguments that were useful in the polynomial volume growth case is relevant here.

Nevertheless it is possible to develop a systematic theory of radial averages which elucidates the basic analytic facts about them (such as maximal inequalities for con­volutions) and to establish ergodic theorems satisfied in general measure-preserving actions of the group.

Two key ideas that will be employed in our analysis are as follows. First, for simplicity of exposition, let us assume that the volume of the balls has exact ex­ponential growth in terms of the radius (as is indeed the case for split rank one groups). Clearly the volume of a shell of unit width, namely $B_{t+1} \setminus B_t$, occupies a fixed proportion of the volume of the ball $B_{t+1}$. It then follows that the ball averages and the shell averages have equivalent maximal operators, so that we can restrict the discussion to the shell averages. On the other hand, the maximal func­tion for the shell averages in Euclidean space is in fact equivalent to the maximal function for the *singular* spherical averages (see [N7] for more details). A s we shall see below, the first key idea is that it is the ideas and techniques of classical singular integral theory, and particularly Tauberian theory, that are most suitable for the analysis of shell (and hence ball) averages on semisimple groups with exponential volume growth. The second key idea is to apply spectral methods based on the uni­tary representation theory of the group in question to prove maximal inequalities and pointwise convergence theorems. This is indeed possible in the case of radial averages on semisimple Lie and algebraic groups, but also in many other cases, namely whenever the singular sphere averages all commute under convolution. Our methods will thus apply in principle to all lcsc groups that admit a radial commu­tative convolution structure, even to amenable ones, and in fact give rise to some interesting results regarding singular averages also on groups of polynomial volume growth.

AS we shall see below, the geometric reduction from ball averages to shell aver­ages together with the use of spectral methods from singular integral theory and the unitary representations theory of the groups restricted to the commutative con­volution subalgebra will serve to replace the growth and Følner conditions on the group used in the polynomial volume growth case, and the existence of the algebraic semi-direct product structure and the Dunford-Zygmund method used otherwise. Before turning to a discussion of the analytic tools involved, we present some basic examples which will serve to motivate our analysis and demonstrate its scope, and formulate some of the pertinent results which will be proved later on. We begin with the following fundamental result on singular averages on $\mathbb{R}^n$.

## 8.1. Euclidean Spherical averages.
Our first example of an ergodic measure preserving action of $\mathbb{R}$ was the flow given by an irrational line on the 2-torus (see Example 3.1). Let us now note that there is an equally basic problem related to geometric averaging on the plane, as follows.

### Example 8.1. Circles in $\mathbb{R}^2$ and spheres in $\mathbb{R}^n$.
Let us denote by $S_t$ a circle of radius $t$ and center 0 in $\mathbb{R}^2$. Let $\sigma_t$ denote the normalized rotation-invariant measure on $S_t$, and we consider the radial averages (=spherical means) that $\sigma_t$ define on $\mathbb{R}^2$. If $f$ a function on $\mathbb{R}^2$, define :

$$\pi(\sigma_t)f(v) = \int_{w \in S_t} f(v + w)d\sigma_t(w)$$

In analogy with Example 3.1 in §3.1, we focus on the action of the spherical means as operators on function spaces on $\mathbb{T}^2$, namely the action on $\mathbb{Z}^2$-periodic functions. We can then consider the natural problems of equidistribution, pointwise almost everyehere convergence to the space average, and (singular !) differentiation. These problems have all been resolved, and we recall the following results, which we formulate for the action of the spherical means on the $n$-torus, for use in later comparisons.

(1) **Equidistribution of Spheres**.
   $\forall v \in \mathbb{T}^n$, $\forall f \in C(\mathbb{T}^n)$

$$\lim_{t \to \infty} \pi(\sigma_t)f(v) = \int_{\mathbb{T}^2} fdm \quad \text{(Exercise !)}$$

(2) **Maximal inequality for sphere averages.**
   $\forall f \in L^p(\mathbb{T}^n)$, $p > \frac{n}{n-1}$,

$$\left\| \sup_{t>0} |\pi(\sigma_t)f| \right\|_{L^p(\mathbb{T}^n)} \leq C_p(n) \|f\|_{L^p(\mathbb{T}^n)}$$

(3) **Singular spherical differentiation**.
   $\forall f \in L^p(\mathbb{T}^n)$, $p > \frac{n}{n-1}$,

$$\lim_{t \to 0} \pi(\sigma_t)f(v) = f(v) \,, \text{ for almost all } v \in \mathbb{T}^n \,.$$

(4) **Pointwise Ergodic Theorem for sphere averages.**
   $\forall f \in L^p(\mathbb{T}^n)$, $p > \frac{n}{n-1}$,

$$\lim_{t \to \infty} \pi(\sigma_t)f(v) = \int_{\mathbb{T}^n} fdm \,, \text{ for almost all } v \in \mathbb{T}^n \,.$$

*Remark* 8.2.     (1) The identification of the limit in (1) and (4) as the space average of the function (namely the mean ergodic theorem), is a simple exercise in spectral theory as will also be seen more generally below. The equidistribution theorem can be proved using a variant of the classical argument of Weyl for the action of a translation on $\mathbb{T}^n$ - it requires only the vanishing at infinity of the characters of the convolution algebra of radial averages, which is a consequence of the Riemann-Lebesgue lemma.

(2) The fundamental result underlying the pointwise convergence theorems (3) and (4) is the maximal inequality (2), which is due to E. Stein (see [SW]) for $n \geq 3$ and to J. Bourgain [Bn] for $n = 2$. The range of $p$ stated is best possible.

(3) As usual, the pointwise convergence in the ergodic theorem (and the differentiation theorem) follows if it holds on a dense subspace. The existence of such a subspace was first established for $n \geq 3$ in [Jo], and for $n = 2$ in [La].

(4) The spherical averages $\sigma_t$ being singular measures on $\mathbb{R}^n$ accounts for the non-trivial restriction on the range $p$ for which the pointwise ergodic theorem holds. This phenomenon was not encountered in the discussion of the absolutely continuous averages that appeared thus far in the previous chapters.

(5) The discrete analog of the spherical means are the averages over integer points lying on a discrete sphere $S_k \cap \mathbb{Z}^n$. Maximal inequalities for the convolution operators defined by such averages on $\mathbb{Z}^n$ were recently established for $n \geq 5$ in [MSW].

*Remark* 8.3. **Measurability of singular maximal functions.**

Since the averages $\sigma_t$ are singular, it is not clear why the maximal function $f_\sigma^*$ is well-defined and measurable even for one function class in $L^p(\mathbb{T}^n)$, incuding the function class 0. The maximal inequality and pointwise convergence should be interpreted as asserting, in particular, that for any two representatives $f$ and $f'$ of a given function class, there exists a co-null set such that for $v$ in this set, both $\pi(\sigma_t)f(v)$ and $\pi(\sigma_t)f'(v)$ exist *for all* $t > 0$ simultaneously and are equal. Then the supremum in the maximal inequality and the limit in the ergodic theorem are indeed well-defined.

This material problem is discussed in detail e.g. in [S2, Ch. XI, §3.5] or alternatively in [Co2, II.4], and we will content ourselves here with noting that $\sup_{t>0} |\sigma_t f(v)|$ is clearly defined for all $v$ and constitutes a measurable function if $f$ is continuous on $\mathbb{T}^n$. The strong $L^p$-maximal inequality for continuous functions (which is sometimes refered to as an a-priori inequality) can then be used to define and prove the measurability and the $L^p$ boundedness of the maximal function for all $L^p$-functions. This remark applies to all other maximal inequalities that will appear in the sequel, since the spaces we consider can always be assume to be compact metric with the $G$-action continuous - see [N2] for details.

Another basic set-up in which spherical averages appear naturally is furnished by the Heisenberg groups.

**Example 8.4. $\mathbb{C}^n$-Spheres in the Heisenberg group**.

Let $H_n = \mathbb{C}^n \times \mathbb{R}$ denote the Heisenberg group, and let $\sigma_t$ denote the normalized rotation invariant measure on the sphere $S_t \subset \mathbb{C}^n$ with center 0 (which we call $\mathbb{C}^n$-spheres). Let $H_n(\mathbb{Z})$ be the discrete subgroup of integer points of the Heisenberg group. Consider the homogeneous space given by $U_n = H_n(\mathbb{Z}) \setminus H_n$, which is compact nilmanifold with a transitive (right) $H_n$-action. We then have the following results

(1) **Equidistribution of $\mathbb{C}^n$-spheres in $U_n$.**

$\forall v \in U_n, \forall f \in C(U_n)$

$$\lim_{t \to \infty} \pi(\sigma_t) f(v) = \int_{U_n} f \, dm \, .$$

(2) **Maximal inequality for $\mathbb{C}^n$-sphere averages**. $\forall f \in L^p(U^n)$, $n > 1$, $p > \frac{2n}{2n-1}$,

$$\left\| \sup_{t>1} |\pi(\sigma_t) f| \right\|_{L^p(U^n)} \le C_p(n) \, \|f\|_{L^p(U^n)} \ .$$

(3) **Singular $\mathbb{C}^n$-sphere differentiation on $U^n$**.
$\forall f \in L^p(U^n)$, $p > \frac{2n}{2n-1}$, $n > 1$

$$\lim_{t \to 0} \pi(\sigma_t) f(v) = f(v) \, , \ \text{for almost all } v \in U^n \, .$$

(4) **Pointwise Ergodic Theorem for $\mathbb{C}^n$-sphere averages**. $\forall f \in L^p(U^n)$, $n > 1$, $p > \frac{2n}{2n-1}$,

$$\lim_{t \to \infty} \pi(\sigma_t) f(v) = \int_{U_n} f \, dm \, , \ \text{for almost all } u \in U_n \, .$$

*Remark* 8.5.    (1) The mean ergodic theorem and also the equidistribution theorem here can be proved spectrally, again in a manner analogous to Weyl's classical equidistribution theorem on $\mathbb{T}^n$. A necessary ingredient in this approach is thus the classification of all the characters of the commutative convolution algebra generated by the $\mathbb{C}^n$-sphere averages on the Heisenberg group.

(2) The $p$-range in (2) and (4) are the best possible for the action on $U^n$, and in fact for general probability-preserving actions of the reduced Heisenberg group. This result was established in [NT].

(3) The maximal inequality for convolutions on the Heisenberg group itself, namely for the operator $\sup_{t>0} |f * \sigma_t|$, was recently established in [MS] for $p > \frac{2n}{2n-1}$, excluding the case $n = 1$. This implies the differentiation theorem, and by the transfer principle for singular averages stated in Theorem 6.2 the same maximal inequality holds for any probability-preserving action of the Heisenberg group, and the pointwise theorem holds as well. Another proof of this result was given by [NaT]. The range $p > \frac{2n-1}{2n-2}$ for $n > 1$ was established earlier in [NT].

(4) Note that in the result above $\sigma_t$ are singular averages supported on *subvarieties of codimension* 2, rather than co-dimension one as in the case of ordinary spheres (see [Co2] for the corresponding result for the spheres of codimension one associated with the natural homogeneous norm). Results for even more singular averages on certain nilpotent Lie groups appear in [MS].

Let us now pass from the considerations above regarding spherical averages in the familiar setting of nilpotent groups, and consider the ergodic theory of spherical averages for lcsc groups which are of exponential volume growth, and non-amenable. We will start with the simplest example, namely that of the isometry groups of hyperbolic space, acting on homogeneous spaces with finite volume.

## 8.2. Non-Euclidean spherical averages.

### Example 8.6. Spheres in Hyperbolic space.

Let $M$ be a compact (or finite volume) Riemann surface, and $\Gamma = \pi_1(M)$ its fundamental group. $\Gamma$ is naturally identified with a lattice subgroup of the isometry group $G = PSL_2(\mathbb{R}) = Iso(\mathbb{H}^2)$ of the hyperbolic plane $\mathbb{H}^2$. $G$ acts by translations on the homogeneous space $\Gamma \setminus G$, and if $K \subset G$ is the maximal compact subgroup of rotations fixing a point $x$, then $M$ can be identified with the double coset space $\Gamma \setminus G/K$. Denote the unique probability measure supported on a sphere of radius $t$ and center $p$ in hyperbolic space $\mathbb{H}^2 = G/K$, which is invariant under the rotations fixing $x$, by $\tilde{\sigma}_t(x)$. Given any (continuous, say) function on $M$, we can lift it to a $\Gamma$-periodic function on $\mathbb{H}^2$ and average it w.r.t. $\tilde{\sigma}_t(x)$. We denote the result of this operation by $\pi(\sigma_t)f(x)$.

More generally, if $\mathbb{H}^n$ denotes hyperbolic $n$-space, and $\Gamma$ is a lattice in $G = Iso(\mathbb{H}^n)$, we can consider the homogeneous space $M = \Gamma \setminus G/K$, $K$ a maximal compact subgroup fixing $x \in \mathbb{H}^n$. $M$ is an $n$-dimensional Riemannian manifold of constant negative sectional curvature provided $\Gamma$ is torsion free. We then define the averaging operators $\pi(\sigma_t)f(x)$ corresponding to the spherical means $\tilde{\sigma}_t(x)$, acting on $\Gamma$-periodic functions on $\mathbb{H}^n = G/K$. We denote the unique $G$-invariant probability measure on $\Gamma \setminus G$ by $m$, as well as its projection onto $M = \Gamma \setminus G/K$. We can now state the following results (for general dimension $n$).

(1) **Equidistribution of spheres in compact hyperbolic homogeneous spaces**.

When $M$ is compact, $\forall x \in M$, $\forall f \in C(M)$

$$\lim_{t \to \infty} \pi(\sigma_t)f(x) = \int_M f\,dm \ .$$

(2) **Maximal inequality for sphere averages in hyperbolic finite-volume homogeneous spaces.**

For any $M$ of finite volume, $\forall f \in L^p(M)$, $p > \frac{n}{n-1}$, $n \geq 2$, and for almost every $x \in M$

$$\left\| \sup_{t > 0} |\pi(\sigma_t)f| \right\|_{L^p(M)} \leq C_p(n)\,\|f\|_{L^p(M)} \ .$$

(3) **Singular spherical differentiation in hyperbolic homogeneous spaces.**

$\forall f \in L^p(M)$, $p > \frac{n}{n-1}$, $n \geq 2$ and for almost every $x \in M$

$$\lim_{t \to 0} \pi(\sigma_t)f(x) = f(x)\,.$$

(4) **Pointwise Ergodic Theorem for sphere averages in hyperbolic finite-volume homogeneous spaces.**

For any $M$ of finite volume, $\forall f \in L^p(M)$, $p > \frac{n}{n-1}$, and almost every $x \in X$

$$\lim_{t \to \infty} \pi(\sigma_t)f(x) = \int_M f\,dm\,.$$

*Remark* 8.7.    (1) The equidistribution of spheres was first proved by G. A. Margulis, who used the mixing property of the geodesic flow. Also the same result holds for any hyperbolic homogeneous space of finite volume, not only compact ones, provided we restrict the function $f$ to be continuous of compact support. Further proofs of these results were also given by [K],[S] and [Kn1]. We note that mixing of the geodesic flow in the case of constant negative curvature is fairly straightforward, as it follows from the general fact that the matrix coefficients of irreducible non-trivial unitary representations of the isometry group of hyperbolic space vanish at infinity. This fact was noted already by Fomin and Gelfand [FG], and was later generalized to the Howe-Moore mixing theorem [HM].

(2) The $L^2$-maximal inequality in (2) and the pointwise ergodic theorem in (4) were proved in [N2] for $\dim M > 2$. Interpolation arguments were used in [N3] and [NS2] to establish the range as $p > \frac{n}{n-1}$, which is best possible. The case of $n = 2$ is treated in [N8], and is based on the results of [I].

(3) The differentiation theorem and the maximal inequality for convolutions on $\mathbb{H}^n$, $n > 2$ is due to [EK], and for $n = 2$ to [I]. Using the transfer principle of Theorem 6.2 (but for the local operator $\sup_{0 < t \le 1} |\pi(\sigma_t)f(x)|$ only !) this result holds also on $M$.

### 8.3. Radial averages on free group.

### Example 8.8. Spheres in the free group.

Let $\mathbb{F}_2$ be the free group on two generators $\mathbb{F}_2 = \mathbb{F}_2(a, b)$, where $a$ and $b$ are two free generators. The symmetric generating set $S = \{a, b, a^{-1}, b^{-1}\}$ determines a word length on $\mathbb{F}_2$ given by (see the discussion in §4.1)

$$|w|_S = \min \{k; \ w = s_{i_1} \cdots s_{i_k}, \ s_{i_j} \in S, 1 \le j \le k\} \ .$$

$d(u, v) = |u^{-1} \cdot v|_S$ is a metric on $\mathbb{F}_2$, which is invariant under the action of $\mathbb{F}_2$ on itself by left translations. The word metric determines the corresponding spheres and balls, and thus also the normalized averaging operators $\sigma_n$ and $\beta_n$ on them.

Now consider the sphere $\mathbb{S}^d \subset \mathbb{R}^{d+1}$, where $d \ge 2$, and let $A$ and $B$ be two orthogonal linear transformations on $\mathbb{S}^d$. Clearly, the assignment $a \mapsto A$ and $b \mapsto B$ extends uniquely to a homomorphism $\mathbb{F}_2 \mapsto O_d(\mathbb{R})$, which defines an action of $\mathbb{F}_2$ by orthogonal transformations on $\mathbb{S}^d$. This action preserves, in particular, the rotation-invariant probability measure $m$ on the sphere. Let us assume for simplicity that the image of $\mathbb{F}_2$ is contained in the connected component of the identity in $O_d(\mathbb{R})$.

(1) **Equidistribution of spherical averages on free group orbits in the unit sphere**.
$\forall x \in \mathbb{S}^d$, $\forall f \in C(\mathbb{S}^d)$

$$\lim_{t \to \infty} \pi(\sigma_t)f(x) = \int_{\mathbb{S}^d} f \, dm \ .$$

(2) **Mean ergodic theorem for spherical averages on free group orbits in the unit sphere.**
$\forall f \in L^2(\mathbb{S}^d)$

$$\lim_{t \to \infty} \left\| \pi(\sigma_t)f - \int_{\mathbb{S}^d} f \, dm \right\|_{L^2(\mathbb{S}^d)} = 0$$

(3) **Maximal inequality for spherical averages on free group orbits in the unit sphere.**
$\forall f \in L^p(\mathbb{S}^d)$, $p > 1$, and for almost every $x \in \mathbb{S}^d$

$$\left\| \sup_{t > 0} |\pi(\sigma_t)f| \right\|_{L^p(\mathbb{S}^d)} \le C_p \|f\|_{L^p(\mathbb{S}^d)} \ .$$

(4) **Pointwise Ergodic Theorem for spherical averages on free group orbits in the unit sphere.**

$\forall f \in L^p(\mathbb{S}^d)$, $p > 1$, and almost every $x \in X$

$$\lim_{t \to \infty} \pi(\sigma_t) f(x) = \int_{\mathbb{S}^d} f \, dm \, .$$

*Remark* 8.9.    (1) The equidistribution result was proved by V. Arnold and A. Krylov in [AK], where in fact they consider an arbitrary connected homogeneous space of a compact Lie group, rather than just the sphere. Furthermore, the question of generalizing von-Neumann mean ergodic theorem and Birkhoff's pointwise ergodic theorem for the sphere and ball averages on the free group is raised explicitly in [AK]. The analogous problems for averages on the group of isometries of hyperbolic $n$-space in also raised there.

(2) The mean ergodic theorem for general probability-preserving action of the free group on $r$ generators was proved by Y. Guivarc'h [Gu2]. The formulation is slightly different, asserting the convergence (in the strong operator topology) of the operators $\sigma_n' = \frac{1}{2}(\sigma_n + \sigma_{n+1})$. This modification is necessary because in a general ergodic action a function $f$ might satisfy $\pi(\sigma_n) f = (-1)^n f$, a situation that does not arise on $\mathbb{S}^d$ because of our density assumption. For more on this periodicity phenomenon, see §10.5 below.

(3) The maximal inequality and the pointwise ergodic theorem for $\sigma_n'$ acting on $L^2$ functions were established in [N1], for all measure-preserving ergodic actions of the free groups. The extension to $L^p$, $p > 1$ was established in [NS1].

All the foregoing examples fall under our general theme of study, which is the analysis of averaging operators arising from families of probability measures $\mu_t$ on an lcsc group $G$ in a probability-measure-preserving action of $G$.

(1) In Example 3.1 the group of course is the real line $\ell \cong \mathbb{R}$, acting by translations on $\mathbb{T}^2 = \mathbb{R}^2 / \mathbb{Z}^2$, $\mu_t = \beta_t =$ ball averages on $\mathbb{R}$.

(2) In Example 8.1, the group is $\mathbb{R}^2$, again acting by translation on $\mathbb{T}^2 = \mathbb{R}^2 / \mathbb{Z}^2$, and $\mu_t = \sigma_t =$ the normalized circle averages on $\mathbb{R}^2$.

(3) In Example 8.3 the group is the Heisenberg group $H_n = \mathbb{C}^n \times \mathbb{R}$, acting by tranlations on the homogeneous space $U_n = H_n(\mathbb{Z}) \setminus H_n$ and $\mu_t = \sigma_t$ are the $\mathbb{C}^n$-sphere averages supported in $\mathbb{C}^n \subset H_n$.

(4) In Example 8.5 the group is $Iso(\mathbb{H}^n)$, the measures are the unique bi-$K$-invariant measures $\sigma_t$ on $G$ projecting to the normalized rotation-invariant measure $\tilde{\sigma}_t$ on a sphere of radius $t$ and center $[K]$ in $\mathbb{H}^n = G/K$, and the action is on the homogeneous space $\Gamma \setminus G$, with its unique invariant probability measure.

(5) In example 8.8 the group is $\mathbb{F}_2$ and the space is $\mathbb{S}^d$ with its unique isometry-invariant probability measure $m$.

## 9. The spectral approach to maximal inequalities

We now turn to an exposition of a spectral approach to ergodic theorems for certain lcsc groups, including semisimple Lie and algebraic groups, and some of their lattices. We will start by demonstrating the method for the basic case of ball and sphere averages in the most accessible connected Lie groups, as follows.

9.1. **Isometry groups of hyperbolic spaces.** Our basic set-up and notation will be as follows :

(1) $\mathbb{H}^n$ =Hyperbolic $n$-dimensional space, with connected isometry group $G = Iso^0(\mathbb{H}^n)$.

(2) $S_t(x)$ = Sphere of radius $t$ with center $x \in \mathbb{H}^n$.

(3) $\tilde{\sigma}_t(x)$ =normalized measure on $S_t(x)$, invariant under the group of rotations fixing $x$.

(4) $\sigma_t$ =the spherical averages on $G$. These are given by $\sigma_t = m_K * \delta_{a_t} * m_K$, where $\{a_t, t \in \mathbb{R}\}$ satisfies $d(a_t o, o) = |t|$, namely its orbit through $o = [K]$ in the symmetric space $\mathbb{H}^n = G/K$ is a geodesic. $\sigma_t$ is the unique bi-$K$-invariant probability measure on $G$ projecting onto $\tilde{\sigma}_t$.

(5) $(X, m)$ a compact metric space with a continuous $G$-action, where $m$ is a $G$-invariant probability measure.

(6) Radial averages (= spherical means) are defined by, for $f \in C(X)$,

$$\pi(\sigma_t)f(x) = \int_G f(g^{-1}x)d\sigma_t(g).$$

The basic example discussed in §8.2 for the set-up above arises when choosing $\Gamma$ to be a discrete group of isometries with fundamental domain of finite volume. Then for $f \in C(\Gamma \setminus G)$, namely a $\Gamma$-periodic function on $\mathbb{H}^n$, the averages above give

$$\pi(\sigma_t)f(x) = \text{ average of f on } S_t(x) \text{ w.r.t. } \tilde{\sigma}_t(x).$$

*Remark* 9.1. We note that the assumption above that our basic Borel $G$-space with $G$-invariant probability measure is a compact metric $G$-space is without loss of generality, as noted in [N2]. In that case for each $F \in C_c^\infty(G)$, and every $f \in C(X)$, the function $h = \pi(F)f$ has the property that $g \mapsto h(g^{-1}x)$ is a $C^\infty$-function on $G$, for every $x \in X$. Clearly the space of such functions is dense in $C(X)$ in the uniform norm, and also in every $L^p(X)$, in the $L^p$-norm, $1 \le p < \infty$. We will thus consider below differentiation operators applied to fuctions in $C_c^\infty(G) * C(X)$ in a general action without further comment.

9.2. **Commutativity of spherical averages.** To analyze the spherical means, we use the following basic observation, originating with Gelfand and Selberg in the 1950's.

**Proposition 9.2.** *The spherical averages $\sigma_t = \sigma_t(o)$ on $Iso(\mathbb{H}^n)$ commute with one another under convolution.*

*Proof.* Consider a one-parameter group of isometries $A = \{a_t, t \in \mathbb{R}\}$, whose orbit through a given point $o$ forms a geodesic in $\mathbb{H}^n$, namely $d(a_t o, o) = |t|$. Now the probability measure $m_K * \delta_{a_t} * m_K$ on $G$ projects under $G \to G/K$ to the unique rotation invariant probability measure on $G/K$, supported on a sphere of radius $t$ and center $o = [K]$. Since $d(a_t o, o) = d(o, a_{-t}o)$, and $K$ is transitive on each sphere with center $o$, we have $a_t = ka_{-t}k'$, and so $Ka_tK = Ka_{-t}K$. Hence the inversion map $g \mapsto g^{-1}$ restricts to the identity on bi-$K$-invariant sets, functions and measures, but also reverses the order of convolution. Now $\sigma_t * \sigma_s$ is bi-$K$-invariant since $\sigma_t$ and $\sigma_s$ are, and hence

$$(\sigma_t * \sigma_s)^\vee = \sigma_t * \sigma_s = \sigma_s^\vee * \sigma_t^\vee = \sigma_s * \sigma_t$$

<div style="text-align: right">□</div>

**Notation**. We denote the algebra of bounded complex bi-$K$-invariant Borel measures on $G$ by $M(G, K)$.

Given any strongly continuous unitary representation $\pi$ of $G$ on a Hilbert space $\mathcal{H}$, each element $\mu$ of $M(G, K)$ is mapped to a bounded operator $\pi(\mu)$. The map $\mu \mapsto \pi(\mu)$ is a continuous algebra homomorphism, commuting with the involutions on $M(G, K)$ and on $\operatorname{End} \mathcal{H}$. When $M(G, K)$ is commutative, we denote by $\mathcal{A}$ the closure in the operator norm topology of $\pi(M(G, K))$, which is a commutative algebra closed under the adjoint operation, and so a commutative $C^*$-algebra. Thus we can appeal to the following fundamental result.

**Spectral theorem.** Let $\mathcal{A}$ be a commutative norm-closed algebra of bounded operators on a Hilbert space $\mathcal{H}$, closed under taking adjoints. Consider the $*$-spectrum $\Sigma^*(\mathcal{A})$ of $\mathcal{A}$, consisting of continuous complex $*$-characters of $\mathcal{A}$, with the $w^*$-topology inherited from the dual $\mathcal{A}^*$ of $\mathcal{A}$. Every $f \in \mathcal{H}$ determines a spectral measure $\nu_f$ on $\Sigma^*(\mathcal{A})$, and the action of an operator $\mu \in \mathcal{A}$ is given by the formula :

$$\langle \pi(\mu)f, f \rangle = \int_{\varphi \in \Sigma^*(\mathcal{A})} \varphi(\mu) d\nu_f(\varphi) \,.$$

**Functional calculus.** The spectral theorem for $\mathcal{A}$ implies that for any bounded measurable function $F$ on $\Sigma^*(\mathcal{A})$, $\pi(F)$ can be interpreted as a bounded operator on $\mathcal{H}$, given by the formula

$$\langle \pi(F)f, f \rangle = \int_{\varphi \in \Sigma^*(\mathcal{A})} F(\varphi) d\nu_f(\varphi)$$

Furthermore, the functional calculus can be extended to more general distributions on the spectrum, including measures and also derivative operators. We will make extensive use of these facts below.

### 9.3. Littlewood-Paley square functions. We now turn to a proof of the following

**Theorem 9.3. Pointwise ergodic theorem in $L^2$ for sphere and ball averages on $Iso(\mathbb{H}^n)$, $n > 2$ [N2]**. *The sphere averages $\sigma_t$ and the ball averages $\beta_t$ on $Iso(\mathbb{H}^n)$ satisfy the pointwise ergodic theorem and the strong maximal inequality in $L^2(X)$, if $n > 2$, for any probability measure-preserving action of $G$.*

For simplicity of exposition, we will consider the following

<u>Model case</u>: **Proof of the pointwise ergodic theorems for $\sigma_t$ and $\beta_t$ in $SL_2(\mathbb{C})$-actions on compact metric spaces.**

The proof proceeds along the following steps (we suppress the notation $\pi$ for the representation for ease of notation) :

(1) First, consider the *uniform average* $\mu_t$ of the spherical measures $\sigma_s$, $0 < s \leq t$, and use Proposition 9.2 to write :

$$\mu_t = \frac{1}{t} \int_0^t \sigma_s ds = m_K * \frac{1}{t} \int_0^t \delta_{a_s} ds * m_K$$

Since $\frac{1}{t}\int_0^t \delta_{a_s}\,ds$ are the Birkhoff averages on $\mathbb{R} \cong A$, they satisfy a strong maximal inequality in every $L^p$, $1 < p < \infty$. It follows immediately that also $\left\|f_\mu^*\right\|_p \leq C_p \left\|f\right\|_p$.

(2) Next, compare $\sigma_t$ to their uniform average $\mu_t$. Using Remark 9.1, for every function $f \in C_c^\infty(G) * C(X)$, the function $g \mapsto f(g^{-1}x)$ is a $C^\infty$ function on $G$ (and thus $s \mapsto \sigma_s f(x)$ is $C^\infty$ on $\mathbb{R}_+$), and we can write :

$$\sigma_t f(x) - \mu_t f(x) = \frac{1}{t}\int_0^t s \frac{d}{ds}\sigma_s f(x)\,ds$$

By the Cauchy-Schwarz inequality :

$$|\sigma_t f(x) - \mu_t f(x)| \leq$$

$$\leq \frac{1}{t}\left(\int_0^t s\,ds \int_0^t s \left|\frac{d}{ds}\sigma_s f(x)\right|^2 ds\right)^{1/2}$$

(3) Estimating, we have :

$$\sup_{t \geq 0} |\sigma_t f(x)| \leq \sup_{t \geq 0} |\mu_t f(x)| + R(f,x) \ ,$$

where we define :

$$R(f,x)^2 = \int_0^\infty s \left|\frac{d}{ds}\sigma_s f(x)\right|^2 ds$$

$R(f,x)$ **is called the Littlewood-Paley square function**.

(4) We have by (3) :

$$\|f_\sigma^*\|_2 \leq \left\|f_\mu^*\right\| + \|R(f,\cdot)\|_2$$

We now compute the norm of the square function by the spectral theorem, namely by going over to the Fourier-Gelfand transform side. We obtain, recalling that $\Sigma^*$ denotes the $*$-spectrum of $\mathcal{A} = M(G,K)$, and $\nu_f$ the spectral measure determined by $f$ on $\Sigma^*$ :

$$\|R(f,\cdot)\|_2^2 = \int_X \int_0^\infty s \left|\frac{d}{ds}\sigma_s f(x)\right|^2 ds\,dm(x) =$$

$$\int_0^\infty s \left\|\frac{d}{ds}\sigma_s f\right\|_{L^2(X)}^2 ds = \int_0^\infty s \int_{\Sigma^*} \left|\frac{d}{ds}\varphi_z(\sigma_s)\right|^2 d\nu_f(z)$$

Here we have obtained a spectral expression for the distribution $\frac{d}{ds}\sigma_s$, using the functional calculus in the commutative algebra $\mathcal{A}$ of spherical averaging operators. We refer to [N2, §6, Lemma 4] for more on this argument.

(5) We can conclude that if the expression :

$$\Phi(z) = \int_0^\infty s \left|\frac{d}{ds}\varphi_z(\sigma_s)\right|^2 ds$$

has a *uniform spectral estimate*, namely a bound independent of $z$, as $\varphi_z$ varies over the spectrum $\Sigma^*$, then the strong $L^2$-maximal inequality is proved.

(6) We recall (see [He, Ch. IV, §5] for full details) that for $G = SL_2(\mathbb{C})$ the characters of $M(G, K)$ are given by

$$\varphi_z(\sigma_t) = \frac{\sinh(zt)}{z \sinh(t)} , \quad \varphi_0(\sigma_t) = \frac{t}{\sinh t}$$

The continuous (i.e. bounded) $*$-characters are parametrized by :
   (i) $z = i\lambda$, $\lambda$ real ; *Principal series*, or
   (ii) $z = a$ where $a$ is real and $0 < a \leq 1$ ; *Complementary series.*

(7) For the principal series, part (6) implies immediately that the $*$-characters decay exponentially in the distance, *uniformly in $\lambda$*, with fixed rate. Furthermore for $SL_2(\mathbb{C})$ the explicit expression above yields

$$\left| \frac{d}{dt} \varphi_{i\lambda}(\sigma_t) \right| \leq C(1 + t) \exp(-t)$$

Therefore $\Phi(i\lambda) \leq C\Phi(0) < \infty$, $\forall \lambda \in \mathbb{R}$.

Note however the foregoing estimate fails for the second derivative, namely the second derivative is *not* bounded uniformly in $\lambda$.

(8) For the complementary series, the $*$-characters decay arbitrarily slowly, and in fact, if $a = 1 - \varepsilon$, then $\varphi_a(\sigma_t) \cong c_\varepsilon \exp(-\varepsilon t)$. Furthermore, it can easily be proved directly from the formula in (6) that here

$$\left| \frac{d}{dt} \varphi_a(\sigma_t) \right| \leq \varepsilon \exp(-\varepsilon t)$$

and therefore

$$\int_0^\infty s \left| \frac{d}{ds} \varphi_a(\sigma_t) \right|^2 ds \leq \varepsilon^2 \int_0^\infty s \exp(-2\varepsilon s) ds \leq C < \infty$$

For future reference we also note that in fact for every derivative of the complementary series characters

$$\left| \frac{d^k}{dt^k} \varphi_{1-\varepsilon}(\sigma_t) \right| \leq C_k \varepsilon^k \exp(-\varepsilon t) .$$

(9) Thus we have established that $\Phi(z) \leq C < \infty$ as $\varphi_z$ ranges over the spectrum $\Sigma^*$, and this suffices to prove the (a-priori) strong $L^2$-maximal inequality for the sphere averages $\sigma_t$. Of course, it follows immediately that the ball averages $\beta_t$ satisfy the same maximal inequality, being convex averages of the sphere averages.

(10) The mean ergodic theorem for $Iso(\mathbb{H}^n)$ is of course a consequence of the Fomin-Gelfand result [FG], or more generally the Howe-Moore mixing theorem [HM], which in particular establishes decay of non-trivial continuous characters of the algebra $\mathcal{A}$. In other words, $\lim_{t\to\infty} \varphi_z(\sigma_t) = 0$, and by the spectral theorem it follows immediately that when $G$ is ergodic

$$\langle \pi(\sigma_t)f, f \rangle = \int_{\varphi \in \Sigma^*} \varphi(\sigma_t) d\nu_f(\varphi) \longrightarrow \int_X f dm$$

In the $SL_2(\mathbb{C})$-case, the latter conclusion also follows upon inspection of the explicit form of the characters. The mean ergodic theorem for the balls is an easy consequence.

(11) According to the recipe of §2.3, the last step left is to prove the pointwise convergence of $\pi(\sigma_t)f(x)$ for a dense set of functions $f$. Here, since the balls, and certainly the spheres, are not Følner sets, the variants of Riesz's argument used in the amenable case cannot be applied. A more technical argument has to employed, described briefly as follows. Decompose the spectrum as a union $\Sigma' = \cup_{\varepsilon > 0}\Sigma_\varepsilon$, where $\Sigma_\varepsilon = \{\varphi_z \,;\, |z| \leq 1 - \varepsilon\}$. We can then use the fact that each of the characters appearing in $\Sigma_\varepsilon$ has (exponential) decay in $t$ of a fixed positive rate $\delta(\varepsilon)$. Consider now functions $f$ which are sufficiently $L^2$-smooth, and whose spectral measure $\nu_f$ is supported in $\Sigma_\varepsilon$. Such functions can be shown to satisfy the desired conclusion, namely $\pi(\sigma_t)f(x)$ converges almost everywhere, using Sobolev-space arguments. Finally, the set of all such functions as $\varepsilon \to 0$ is dense in the set of $K$-invariant functions, and this completes the proof in the case of sphere averages, and the case of ball averages easily follows. We refer to [N2] for more details on these arguments.

### 9.4. Exponential volume growth : ball versus shell averages.

Theorem 9.3 established the strong maximal inequality and pointwise ergodic theorem in $L^2$ for sphere and ball averages. Of course, for the absolutely continuous ball averages one would expect a similar result in $L^1$, or at least in $L^p$, $1 < p < \infty$. The $L^1$-problem is still open, and we now describe the proof of the following.

**Theorem 9.4. Pointwise ergodic theorem in $L^p$, $1 < p < \infty$ for ball averages on $Iso(\mathbb{H}^n)$[N3] [NS2].** *The ball averages $\beta_t$ satisfy the pointwise ergodic theorem and and strong maximal inequality in $L^p(X)$, $p > 1$, for all dimensions $n \geq 2$.*

<u>Model case</u>: **Proof of the pointwise ergodic theorems for $\beta_t$ on $SL_2(\mathbb{C})$.**

(1) The ball averages $\beta_t$ satisfy, for $t \geq 1$ :

$$\beta_t = \frac{\int_0^t (\sinh s)^2 \sigma_s ds}{\int_0^t (\sinh s)^2 ds} \leq C_1 e^{-2t} \int_0^t e^{2s} \sigma_s ds \ .$$

(2) Consider the shell average, for $t \geq 1$, and the corresponding maximal function, given by :

$$\gamma_t = \int_0^1 \sigma_{t-s} ds \,, \quad f_\gamma^*(x) = \sup_{t \geq 1} |\pi(\gamma_t)f(x)| \,.$$

(3) We now use the exponential volume growth of balls in $G$, in order to bound $f_\beta^*$ by $f_\gamma^*$. Since $\gamma_t$ is the uniform average of spheres with radius in $[t-1, t]$ (here $t \geq 1$), this amounts to comparing the average on a ball of radius $r$ to the maximum of the averages on annuli of width one and radii bounded by $r$. Thus, using the foregoing estimate of the densities, for $f \geq 0$ we have

$$\beta_t f(x) \leq C_1 e^{-2t} \int_0^t e^{2s} \sigma_s f(x) ds \leq$$

$$\leq C_1 e^{-2t} \Big( \sum_{k=0}^{[t]-1} \int_{t-k-1}^{t-k} e^{2s} \sigma_s f(x) ds + \gamma_1 f(x) \Big) \leq$$

$$\leq C_1 \sum_{k=0}^{[t]-1} e^{-2k} \gamma_{t-k} f(x) + C_1 \gamma_1 f(x) \leq$$

$$\leq C_1 \left( \sum_{k=0}^{\infty} e^{-2k} \right) \sup_{1 \leq s \leq t} \gamma_s f(x) + C_1 \gamma_1 f(x) \leq 4 C_1 f_\gamma^*(x) \,.$$

Hence it suffices to prove the maximal inequality for the shell averages. One advantage that $\gamma_t$ offers is that the exponential density that weighs the sphere averages in $\beta_t$ no longer appears, and we can use classical methods of Fourier analysis to estimate the Gelfand transform without difficulty . In the next section we will estimate the rate of decay of the transform, and then apply some analytic interpolation techniques to estimate certain maximal functions associated with the (regularized) shell averages $\gamma_t$. These techniques will allow us to convert $L^2$-boundedness results for square functions associated with the derivative of (regularized) shell averages, to $L^p$-boundedness results for the maximal function associated with the shells themselves.

This completes the proof of Theorem 9.4 for ball averages, provided we prove the corresponding result for the shell averages. $\qquad\square$

**9.5. square functions and analytic interpolation.** We will presently show how to pass from the norm boundedness of the square functions in $L^2$ to best-possible $L^p$ results, via analytic interpolation. This technique is most easily implemented for the shell averages, a fact that motivates their introduction to our discussion. However, by the discussion of §9.4, to complete the proof of Theorem 9.4, it suffices indeed to consider the shell averages and prove the following.

**Theorem 9.5. Pointwise ergodic theorem for shell averages** $\gamma_t$ **on** $Iso(\mathbb{H}^n)$, $n \geq 2$[N3] [NS2]. *The shell averages* $\gamma_t$ *satisfy the pointwise ergodic theorem and the strong maximal inequality in* $L^p(X)$, $p > 1$, *for all dimensions* $n \geq 2$.

<u>**Model case**</u>: **Proof of the pointwise ergodic theorems for** $\gamma_t$ **on** $SL_2(\mathbb{C})$.

(1) First, let us smooth the shell averages by the usual procedure, and define, for $t \geq 1$ :

$$\tilde{\gamma}_t = \int_{\mathbb{R}} \psi(t-s) \sigma_s ds$$

where $\psi$ is a positive smooth function identically one on $[0,1]$, vanishing outside $[-1,2]$.

(2) We would like now to consider all the derivative operators $\frac{d^k}{ds^k} \tilde{\gamma}_s$, and bound the expressions

$$\Phi_k(z) = \int_1^\infty s^{2k-1} \left| \frac{d^k}{ds^k} \varphi_z(\tilde{\gamma}_s) \right|^2 ds$$

independently of $z$, namely uniformly on the spectrum $\Sigma^*(\mathcal{A})$, as in §9.2.

(3) For shell averages (unlike the case of the singular sphere averages $\sigma_s$) the smoothing allows uniform control (as $\lambda \to \infty$) of any given derivative of *principal series* characters. Indeed, recall that for $SL_2(\mathbb{C})$ these characters are given by $\frac{\sin \lambda t}{\lambda \sinh t}$. Thus the desired boundedness is an easy consequence of classical 1-dimensional Fourier theory, since here the problem reduces to estimating the decay in $\lambda$ of the Fourier transform of a smooth compactly-supported function on the line. As is well-known, the decay in this case is faster than any polynomial in $\lambda$, and the boundedness of the integral above follows.

(4) Arbitrarily high derivatives of the *complementary series* characters can also be controlled, in fact even when evaluated on $\sigma_t$, and therefore also for $\tilde{\gamma}_t$. This was already noted in §9.3, part (8), and can be proved by differentiating the explicit form of the characters.

(5) Therefore, using the functional calculus for the distributions corresponding to higher derivatives, the $k$th order square functions

$$R_k(f,x)^2 = \int_1^\infty s^{2k-1} \left| \frac{d^k}{ds^k} \tilde{\gamma}_s f(x) \right|^2 ds$$

have an $L^2$-norm bound.

(6) Now use the Riemann-Liouville fractional integral family of operators (see e.g. [S1]), and embed $\tilde{\gamma}_t$ in an analytic family of operators $T_z$, $z \in \mathbb{C}$. By the analytic interpolation theorem, an $L^2$-norm bound for the derivative of $\tilde{\gamma}_t$ appearing in a square function can be converted to an $L^p$-norm bound for the operator $\sup_{t \geq 1} \tilde{\gamma}_t$. This is done by interpolating against the maximal inequality that the uniform averages $\frac{1}{t} \int_0^t \tilde{\gamma}_s ds$ satisfy in every $L^p$, $p > 1$. The latter result was noted for the averages $\mu_t = \int_0^t \frac{1}{t} \sigma_s ds$ in §9.3, part (1), and of course it follows in exactly the same way in the present case. We refer to [N3] and [NS2] for the details.

**9.6. The $L^p$-theorem for sphere averages on $Iso(\mathbb{H}^n)$.** It is of course natural to complete also the discussion of the maximal inequality for the sphere averages, whose boundedness was established in $L^2$, and prove the best possible results in $L^p$. Indeed, a variant of the method of analytic interpolation via the Riemann-Liouville fractional integrals can be applied to $\sigma_t$ also, in order to embed it in an analytic family of operators. This method yields the following :

**Theorem 9.6. Pointwise ergodic theorem for the sphere averages $\sigma_t$ on $Iso(\mathbb{H}^n)$, $n > 2$ [N3][NS2].** *The sphere averages $\sigma_t$ on $Iso(\mathbb{H}^n)$, $n > 2$ satisfies the pointwise ergodic theorem and strong maximal inequality in $L^p$, $p > \frac{n}{n-1}$, which is the best possible range.*

**Problem 9.7.** The pointwise ergodic theorem and $L^p$, $p > 2$ maximal inequality for sphere averages in general measure-preserving actions of $Iso(\mathbb{H}^2)$ is an open problem.

*Remark* 9.8.     (1) The constraint on the $L^p$-range of the maximal inequality for spheres in hyperbolic space is the same as the Euclidean constraint. The constraint is determined by the rate of decay (in the spectral variable) of the Fourier-Gelfand transform of the spherical measure on a sphere of radius one. We note that the counter-example in the Euclidean case for boundedness of the maximal operator in $L^{\frac{n}{n-1}}$ can be taken to be a local one [SW]. Namely, it is given by a function of compact support and a singularity at the origin (say), and thus the $L^p$-constraint arises already from the local operators $\sup_{0 < t \leq 1} \sigma_t$. Now since hypebolic spheres (in the ball model) are just off-center Euclidean spheres, the Euclidean local counter-example is also a hyperbolic counter-example [N3, §5.4]. The point is thus to show that there are no further constraints in the hypebolic case.

(2) The complementary series poses a serious challenge in the analysis employed to prove the maximal inequality $\sup_{1 < t < \infty} \sigma_t$ for the sphere averages on the groups $Iso(\mathbb{H}^n)$. This is the result of the arbitrarily slow rate of decay of the complementary series characters. To control the norm of square functions, it is necessary to prove derivative estimates for the spherical functions of the complementary series which are significantly better than the standard Harish Chandra estimates. This problem takes a sizable part of the effort in [N3] and [NS2].

(3) Let us emphasize that the behaviour of the averages $\gamma_t$ is in marked contrast to the Euclidean case. In the hyperbolic set-up, the maximal operators associated with $\gamma_t$ and $\beta_t$ are in fact equivalent. But the Euclidean shell averages satisfy the same maximal inequalities as the sphere averages, and not the same maximal inequalities as the ball averages. This fact is a reflection of the difference between polynomial and exponential volume growth, and for more on this matter we refer to [N6] and [N7].

(4) Recalling our comments in the introduction to chapter 8, we note that the only available proof of the maximal inequality for the ball averages $\beta_t$ in $L^p$, namely the proof described in §9.3 - §9.5 above, makes use of differentiation theory of the singular averages $\sigma_t$. This is also in marked contrast to the Euclidean (or polynomial volume growth) case, where the singular sphere averages did not play any role.

(5) A more geometric approach to obtain the $L^p$, $p > 1$ maximal inequality for balls in every dimension greater than two, is to start with $SL_2(\mathbb{C})$ and use an analog of the "method of rotations". Namely, embed $Iso(\mathbb{H}^3) \subset Iso(\mathbb{H}^n)$ as the stability group of a totally geodesic subspace. Then the Cartan polar coordinates decompositions $G = KAK$ in the two groups can be aligned. Since $A$ is one dimensional, the maximal inequalities for the ball averages in $Iso(\mathbb{H}^n)$ follows from those of $\sigma_t$ in $Iso(\mathbb{H}^3)$. However, this leaves out $Iso(\mathbb{H}^2)$, where the spherical functions estimate are the hardest case. For spheres, the range of $p$ where the strong maximal inequality holds improves with dimension, so this method does not give optimal results.

## 10. Groups with commutative radial convolution structure

The methods outlined in the previous section have a wide scope of applications, and have been developed into a systematic spectral approach to the proof of pointwise ergodic theorems for radial averages on lcsc groups admitting a commutative radial convolution structure. In the present section we will indicate some of the results obtained in this direction, and comment on some of the open problems.

### 10.1. Gelfand pairs. Let us recall the following well known definition :

**Definition 10.1. Gelfand Pairs**. A Gelfand pair $(G, K)$ consists of an lcsc group $G$, and a compact subgroup $K \subset G$, such that the algebra $M(G, K)$ of bounded Borel measures on $G$ which are bi-$K$-invariant is commutative, or equivalently, the convolution algebra $L^1(G, K)$ of bi-$K$-invariant $L^1$-functions is commutative. We remark that $G$ is then necessarily unimodular.

### Example 10.2. Some examples of Gelfand pairs.

(1) $G$ a connected semisimple Lie group, $K$ a maximal compact subgroup. This extensive family includes :

(i) $G = SO(n,1) = Iso(\mathbb{H}^n)$ the isometry group of the simply connected Riemannian $n$-manifold of constant negative curvature, $K = SO(n)$ the group of rotations fixing a point. Thus here $M(G,K)$ is the usual algebra of radial averages on hyperbolic space.

Similarly, $G = SU(n,1)$, the isometry group of complex hyperbolic space, $K$ the maximal compact subgroup fixing a point.

(ii) $G = O(p,q)$, the isometry group of simply-connected pseudo-Riemmanian manifold of signature $(p,q)$ and constant curvature, $K = O_p(\mathbb{R}) \times O_q(\mathbb{R})$.

(iii) $G = SL_n(\mathbb{C})$, the general Linear group, which is the isometry group of the space of positive definite matrices, $K = SU_n(\mathbb{C})$ the unitary group.

(iv) $G = Sp_n(\mathbb{R})$, the symplectic group, $K = Sp(n)$.

(2) $S = K \ltimes \mathfrak{p}$, a Cartan motion group. Here $K$ is a maximal compact subgroup of a semisimple Lie group $G$, and $\mathfrak{p}$ the fixed-point-subspace of a Cartan involution on the semisimple Lie algebra $\mathfrak{g}$. Examples include :

(i) The Euclidean motion group $S = O_n(\mathbb{R}) \ltimes \mathbb{R}^n$, $K = O_n(\mathbb{R})$. Here $M(G,K)$ is the usual algebra of radial measures on $\mathbb{R}^n$.

(ii) The Heisenberg motion group, $S = U_n(\mathbb{C}) \ltimes H_n$, $K = U_n(\mathbb{C})$. Here $M(G,K)$ is the algebra of radial measures on $H_n$ generated by the $\mathbb{C}^n$-spheres.

(iii) $S = O_n(\mathbb{R}) \ltimes Sym_n(\mathbb{R})$, $K = O_n(\mathbb{R})$, where $Sym_n(\mathbb{R})$ is the space of $n \times n$ symmetric matrices, and the action of $O_n(\mathbb{R})$ is by conjugation.

(3) $G$ a connected semisimple algebraic Chevalley group over a locally compact non-discrete field, $K =$ a good compact open subgroup. Examples include :

(i) $G = PGL_2(\mathbb{Q}_p)$, $K = PGL_2(\hat{\mathbb{Z}}_p)$. Here $G/K$ is a $(p+1)$-regular tree, $M(G,K)$ the algebra of radial avegaes on the regular tree.

(ii) $G = SL_n(\mathbb{Q}_p)$, $K = SL_n(\hat{\mathbb{Z}}_p)$. Here $G$ acts by isometries of an affine building, and $M(G,K)$ is a commutative subalgebra of the convolution Hecke algebra of double cosets of an Iwahori subgroup.

(4) $G = Aut(T_{r_1,r_2})$ the automorphism group of the semi-homogeneous bipartite tree of valencies $r_1$ and $r_2$, $K_i$ the maximal compact subgroup fixing a vertex of valency $r_i$. Further examples of Gelfand pairs $(G,K)$ arise from certain closed non-compact boundary-transitive subgroups of the group of automorphisms of a finite product of such trees.

The foregoing list is a very partial one (for a discussion of some more examples of amenable Gelfand pairs see e.g. [BJR1][BJR2]). However already the groups mentioned (together with other groups possesing a radial convolution structure, e.g. the free groups) give rise to a large collection of interesting measure-preserving actions. We mention briefly only the following.

**Example 10.3. Some examples of measure-preserving actions.**

(1) $SL_2(\mathbb{R})$ acts (transitively) on the unit tangent bundle of a compact Riemann surface $M = \pi_1(M) \setminus SL_2(\mathbb{R})$, preserving a volume form of finite total mass.

(2) $SL_n(\mathbb{R})$ acts (transitively) on the space of unimodular lattices $\mathcal{L}_n$ in $\mathbb{R}^n$, namely $\mathcal{L}_n = SL_n(\mathbb{Z}) \setminus SL_n(\mathbb{R})$, preserving a volume form of finite total mass. Of course, ane subgroup of $SL_n(\mathbb{R})$ also acts on $\mathcal{L}_n$.

(3) More generally, any subgroup $H$ of a semisimple algebraic group $G$ acts on the probability space $\Gamma \setminus G$, where $\Gamma$ is a lattice subgroup, e.g.

(i) $G = SL_2(\mathbb{R}) \times SL_2(\mathbb{R})$ and $\Gamma = SL_2(\mathbb{Z}[\sqrt{d}])$ under the skew diagonal embedding $\gamma \mapsto (\gamma, \tau(\gamma))$, $\tau$ the Galois automorphism of $\mathbb{Z}[\sqrt{d}]$, $d$ a square free positive integer.

(ii) $G = SL_2(\mathbb{R}) \times SL_2(\mathbb{Q}_p)$, $\Gamma = SL_2(Z[\frac{1}{p}])$ an irreducible lattice in $G$.

(4) $\Gamma = \mathbb{F}_t$ the free group, for example embedded as a lattice in $G = PGL_2(\mathbb{Q}_p)$, as $\mathbb{Z}[\frac{1}{p}]$-points in an appropriate quaternion algebra, and $X$ any action of $G$, e.g. on a compact locally symmetric space.

(5) $\Gamma = SL_2(\mathbb{Z})$, and the action by group automorphisms of $\mathbb{T}^2$, or on $SL_2(\mathbb{R})/\Gamma$, $\Gamma$ a lattice.

Again, this list is very partial, but it already clearly demostrates that establishing a pointwise (or mean) ergodic theorem (preferably with error term), a strong maximal inequality or a differentiation theorem gives rise to diverse applications, depending on the family $\mu_t$ to which it applies, and the group and action involved. We will mention here some applications only very briefly, just to motivate our discussion, and without attempting to explain them further. Rather, we will concentrate below on the proof of the ergodic theorems themselves.

**Example 10.4. Some applications of ergodic theorems and maximal inequalities.**

(1) Boundedness properties of natural singular integrals on homogeneous spaces.

(2) Integral geometry on locally symmetric spaces, e.g. singular spherical differentiation, and pointwise equidistribution of spheres and other singular subvarieties in homogeneous spaces.

(3) Evaluation of the main term in counting lattice points on homogeneous algebraic varieties.

(4) Estimating error terms in problems of diophantine approximation on homogeneous algebraic varieties and homogeneous spaces.

10.2. **Pointwise theorems for commuting averages : General method.** We now assume $(G, K)$ is a Gelfand pair, and as usual $(X, m)$ denotes an ergodic probability measure preserving action of $G$. The spectral approach outlined in §9 to pointwise convergence is based on applying certain geometric properties of the convolution structure of $M(G, K)$, together with the tools of harmonic analysis on Abelian Banach algebras, in order to prove pointwise convergence when the elements of the algebra are represented as averaging operators on $L^p(X)$. Since we consider representations arising from a measure-preserving and thus unitary action of $G$, the homomorphism $\pi : M(G, K) \to \mathrm{End}\, L^2(X)$ is a $*$-homomorphism. Furthermore, an algebra representation arising from a unitary representation of the group, gives rise to $*$-characters (also called $(G, K)$-spherical functions), which can all be identified with *positive-definite* continuous bounded functions on the group - see [GV, Ch I] and [He, Ch IV]. The basic measures in $M(G, K)$ are $\sigma_g = m_K * \delta_g * m_K$, $g \in G$. Every other bi-$K$-invariant probability measure on $G$ is a convex combination of these. Note that in the discussion of §9, the one-parameter group $A = \{a_t, \; t \in \mathbb{R}\}$ was a fixed group of hyperbolic isometries, and there it was enough to consider $\sigma_t = m_K * \delta_{a_t} * m_K$. Many other possibilities for the choice of the averages arise in practice, and it is not neccesary to restrict to averages associated with one-parameter groups, or even to one parameter family of averages - we refer to [N6] for some examples and more details. Nevertheless, for simplicity of exposition, let us choose a family of bi-$K$-invariant probability measures $\nu_t$, $t \in \mathbb{R}$ (not necessarily of the form $m_K * \delta_{a_t} * m_K$), and explain briefly the ingredients sufficient for a proof of the ergodic theorems and maximal inequalities in this case. We will also comment briefly in the next section on the various problem that arise along the way. The recipe for the proof proceeds along the following steps.

(1) Identify the positive-definite *-spectrum $\Sigma^*(\mathcal{A})$. Prove that for every non-trivial positive-definite *-character $\varphi_z \in \Sigma^*(\mathcal{A})$ (i.e. every non-constant positive-definite spherical function)

$$\lim_{t \to \infty} \varphi_z(\nu_t) = 0$$

Conclude that $\nu_t$ satisfy the **mean ergodic theorem** :

$$\lim_{t \to \infty} \left\| \nu_t f - \int_X f \, dm \right\|_{L^2(X)} = 0$$

(2) Analyze the behaviour in $z$ of $\varphi_z(\nu_t)$ and $\frac{d^k}{dt^k} \varphi_z(\nu_t)$ as $t \to \infty$ (for a pointwise ergodic theorem) and as $t \to 0$ (for a differentiation theorem).

(3) Establish the existence of a dense set of functions in $L^2$ where $\nu_t f(x)$ converges pointwise almost everywhere, using the spectral decomposition of $\mathcal{A}$ in $L^2(X)$, estimates of the positive-definite *-characters, and Sobolev space arguments (see [N2]).

(4) Establish a maximal inequality in $L^p$, $p > 1$ for an averaged version $\mu_t$ of $\nu_t$, for example the uniform averages $\mu_t = \frac{1}{t} \int_0^t \nu_s \, ds$. This may be achieved using a number of methods, depending on the case at hand, as follows.

I) First possibility [N2]: Use a maximal inequality for the radial components of $\nu_t$, namely the averages $\nu_t$ determine on the Abelian group $A$, when a Cartan polar decomposition $G = KAK$ is available. Namely use the representation of $\nu_t$ as a convex combination of the basic measures $\sigma_g = m_K * \delta_a * m_K$, $a$ the Cartan component of $g$.

II) Second possibility [NS1] : Use a central limit theorem for the transient random walk associated with $\nu_1$ and conclude that (for fixed positive constants $c$ and $C$)

$$\mu_n \le \frac{C}{cn+1} \sum_{k=0}^{cn} \nu_1^{*k}$$

Then use the Hopf-Dunford-Schwarz maximal inequality for uniform averages of $\nu_1^{*k}$. Finally, argue that for some constant $B$, $\mu_t \le B \mu_{[t]+1}$, $[t]$ the integer part of $t$ [MNS].

III) Third possibility [N1]: Establish the following subadditive convolution inequality given by (for fixed positive constants $c$ and $C$)

$$\mu_t * \mu_s \le C \left( \mu_{ct} + \mu_{cs} \right) ,$$

using convolution estimates on the group. The subadditive convolution inequality in the group algebra id sufficient to deduce strong maximal inequality in $L^2$, in every action of the group (see also §10.5 below). A strong maximal inequality in $L^p$, $p > 1$ can be deduced if an iterated form of the subadditive inequality is established, for the products $\mu_{t_1} * \mu_{t_2} * \cdots * \mu_{t_n}$.

(5) Estimate the difference $|\sigma_t - \mu_t|$ using an appropriate Littlewood-Paley square function. For the particular case of uniform averages $\mu_t$, as noted

in §9.3 :

$$\nu_t - \mu_t = \frac{1}{t} \int_0^t s \frac{d}{ds} \nu_s \, ds$$

and

$$|\nu_t f(x)| \leq |\mu_t f(x)| + \left( \int_0^\infty s \left| \frac{d}{ds} \nu_s f(x) \right|^2 ds \right)^{1/2}$$

So $f_\nu^*(x) \leq f_\mu^*(x) + R(f, x)$.

(6) Transfer the estimate of the $L^2$-norm of the square function $R(f, x)$ to the Fourier-Gelfand transform side, using the functional calculus in the commutative algebra $\mathcal{A}$. Namely, show that

$$\|R(f, \cdot)\|_2^2 \leq \|f\|_2^2 \sup_{z \in \Sigma^*(\mathcal{A})} \Phi(z)^2 \ ,$$

where :

$$\Phi(z)^2 = \int_0^\infty s \left| \frac{d}{ds} \varphi_z(\nu_s) \right|^2 ds \ .$$

(7) Use the estimates of the characters and their derivative in (2) to show that $\Phi(z)$ has a bound *independent of* $z$. Do the same for the square functions $R_k(f, x)$ corresponding to higher derivatives.

(8) To convert $L^2$-norm bounds for the square functions associated with the derivative operators $\frac{d^k}{ds^k} \nu_s$, to an $L^p$-norm bound for the maximal function $f_\nu^*(x) = \sup_{t>0} |\nu_t f(x)|$, embed $\nu_t$ and their derivatives (and integrals) in an analytic family of operators. For example, we have utilized in §9.4(6) the Riemann-Liouville fractional integral operators in the case of $\mu_t$ (for the use of other families see [NS2]). This allows the use of the analytic interpolation theorem to interpolate between the maximal inequality for the derivative operators in $L^2$, and the $L^p$-maximal inequality for the uniform averages $\mu_t$, for $p > 1$. The latter maximal inequalities are a consequence of the methods indicated in part (4) of the present recipe.

## 10.3. Pointwise theorems and the spectral method : Some open problems.
Each ingredient in the recipe outlined in §10.2 above poses certain difficulties, depending on the Gelfand pair and the family of averages under consideration.

We indicate briefly some of the open problems that arise, taking the example of semisimple Lie and algebraic groups.

(1) The classification of positive-definite spherical functions (namely the $*$-characters of $L^1(G, K)$ that arise in the representations under consideration) on semisimple Lie (and algebraic) groups is far from complete, for many infinite families of groups.

(2) Even when the classification is complete, the best estimates for the derivatives of the spherical functions are far from sufficient to prove the ergodic theorems for (say) the sphere averages. It is thus necessary to establish decay estimates uniformly over the positive-definite $*$-spectrum, which considerably improve Harish Chandra's classical estimates, for example.

(3) In semisimple Lie (or algebraic) groups of real (or split) rank greater than two radial averages more singular that Riemannian spheres occur very naturally. For very singular averages of high codimension, optimal results require very precise estimates of the spherical functions, including along the different singular directions in the radial variable $H \in \mathfrak{a} = Lie(A)$, where $G = KAK$ is the Cartan decomposition. Such estimates are usually not available for higher rank groups.

(4) The latter phenomenon manifests itself even in the case of product groups, where spheres and balls are a convex average of spheres and balls on the component groups, taken with exponential weights. The complexity of the convolution structure makes it difficult to estimate the spherical functions and their derivatives, when evaluated on these averages. In particular, the best possible range of $L^p$ maximal inequalities for sphere averages is not known.

(5) For sphere averages on $Iso(\mathbb{H}^2)$, the maximal inequality in not valid in $L^2$ (as in the case of the Euclidean plane). A similar problems arises for Cartan motion groups, for example the Heisenberg group $H_1$. This makes the spectral methods much less effective, and indeed in these cases, the pointwise ergodic theorem for sphere averages in general actions of the group has not been established. Note that the convolution case for $Iso(\mathbb{H}^2)$ has been settled in [I], but in the absence of a transfer principle for non-amenable groups, this result has no bearing on the case of general ergodic actions.

(6) Spectral methods and analytic interpolation theory do not give maximal inequalities and ergodic theorems in $L^1$. This is an open problem even for the ball averages, on all semisimple groups, and in particular, for $G = SL_2(\mathbb{C})$. Note that a general weak-type $(1,1)$-maximal inequality for ball averages on semisimple groups, acting by convolutions on the symmetric space, has been established in [Str]. Again, however, the absence of a transfer principle renders this result irrelevant for the case of general ergodic actions.

While the list of problems above certainly poses some formidable challenges, in many interesting situations these challenges can be surmounted, and the general spectral method outlined in §10.2 can be implemented (for sphere averages, say). This is true particularly in the case of the action by convolution, which is much more explicit than a general action. This fact is demonstrated by the examples discussed in §8.1, namely the maximal inequalities for convolution with spheres in Euclidean spaces [SW], Heisenberg groups [Co2] [NT][NaT], and hyperbolic spaces of dimension $n \geq 3$ [EK], or dimension $n = 2$ [I]. By the discussion in Chapter 9, maximal inequalities for sphere averages can also be established for general ergodic actions of $G = Iso(\mathbb{H}^n)$, $n > 2$. More generally, the spectral method gives the best possible range of $p$ where the maximal inequality and pointwise ergodic theorem for sphere averages hold, for any simple Lie group of real rank one [NS2] (except $SL_2(\mathbb{R})$).

We will present in the following two sections further examples where the spectral method can be fully or partially implemented for sphere averages, and in the following chapter, further examples where it can be implemented for ball averages. We first turn to the case of singular averages on higher real-rank semisimple groups, which is far from completely solved, and exhibits many of the problems refered to in §10.3. The only available results on singular averages on higher-real rank groups were obtained in the case of complex group, which we discuss in section 10.4. In section 10.5 we will consider some totally disconnected Gelfand pairs, as well as some of their discrete lattice subgroups, to which the spectral method applies. In particular we will prove the ergodic theorems of spheres in the free groups $\mathbb{F}_k$, noted in §8.3.

10.4. **Sphere averages on complex groups.** We begin by a brief reminder of the basic relevant set-up and notation. Let $G$ denote a connected semi-simple Lie group with finite center and without non-trivial compact factors, $\mathfrak{g}$ its Lie algebra. Let $\theta$ denote a Cartan involution on $G$ and $\mathfrak{g}$, and let $\mathfrak{g} = \mathfrak{k} \oplus \mathfrak{p}$ be the corresponding Cartan decomposition, so that $\mathfrak{k}$ is the Lie subalgebra corresponding to a connected maximal compact subgroup $K$. Let $\mathfrak{a} \subset \mathfrak{p}$ denote a maximal Abelian subalgebra, and let $\Phi(\mathfrak{a}, \mathfrak{g}) = \Phi \subset \mathfrak{a}^*$ denote the (real) root system of $\mathfrak{a}$ in $\mathfrak{g}$. Let $\mathfrak{g}_\alpha$ denote the root space corresponding to $\alpha \in \Phi$, and $\mathfrak{g} = \mathfrak{m} \oplus \mathfrak{a} \oplus \sum_{\alpha \in \Phi} \mathfrak{g}_\alpha$ the root space decomposition. Fix a system of simple roots $\Delta \subset \Phi$, the corresponding ordering of $\mathfrak{a}^* = \hom(\mathfrak{a}, \mathbb{R})$ and the system of positive roots $\Phi_+$, and let $\rho$ denote half the sum of the positive roots. Let $W = W(\mathfrak{a}, \mathfrak{g})$ denote the Weyl group of the root system, $\mathfrak{a}_+$ the positive Weyl chamber, and $\overline{\mathfrak{a}_+}$ its closure. Let $A = \exp \mathfrak{a}$ denote the Lie subgroup corresponding to $\mathfrak{a}$, $A_+ = \exp \mathfrak{a}_+$ and $\overline{A_+}$ its closure. The Cartan (or polar coordinates) decom position in $G$ is given by $G = K \overline{A_+} K$ and $g = k_1 e^{H(g)} k_2$, where $H(g)$ is the $\overline{\mathfrak{a}_+}$ component of $g$. Let $\langle , \rangle$ denote the Killing form on $\mathfrak{g}$, and let $d$ denote the induced Riemannian metric on the symmetric space $G/K$. Then the restriction of $\langle , \rangle$ to $\mathfrak{a}$ is an inner product, and we have $d(\exp(H)o, o) = \sqrt{\langle H, H \rangle}$ for all $H \in \mathfrak{a}$, where $o = [K]$ denotes our choice of origin in $G/K$. We recall that the Cartan polar coordinates decomposition yield the following integration formula for Haar measure on $G$ (See [He, p. 186] or [GV]) :

$$\int_G f(g) dm_G(g) = \int_K \int_{\overline{\mathfrak{a}_+}} \int_K f(ke^H k') \, \xi(H) \, dm_K(k) dH \, dm_K(k')$$

Here $\xi(H) = \prod_{\alpha \in \Phi_+} (\sinh \alpha(H))^{m_\alpha}$, $H \in \overline{\mathfrak{a}_+}$, $m_\alpha = \dim_{\mathbb{R}} \mathfrak{g}_\alpha$, and $m_G$, $m_K$ denotes Haar measures on $G$ and $K$, $dH$ denotes Lebesgue measure on $\mathfrak{a}$.

We can now formulate the following result for complex smisimple Lie groups.

**Theorem 10.5. Pointwise ergodic theorem for sphere averages on complex groups** [CN]. *Let $G$ be a connected complex semisimple Lie group with finite center. Fix a regular direction $H \in \mathfrak{a}$, and let $\sigma_t^H = m_K * \delta_{\exp tH} * m_k$. Let $\sigma_t$ denote the Riemannian sphere averages. In every measure-preserving action of $G$, we have*

(1) *The averages $\sigma_t^H$ satisfy the pointwise ergodic theorem, strong maximal inequality and the singular differentiation theorem in $L^p$, $p > p_G$. Here $p_G$ is an explicit computable constant, and e.g. for $G = SL_n(\mathbb{C})$, $p_G = \frac{2n-1}{2n-2}$.*

(2) *The Riemannian sphere averages $\sigma_t$ satisfy the pointwise ergodic theorem, maximal inequality and singular differentiation theorem in $L^2$.*

We remark that it is not known whether the range of $p$ stated in Theorem 10.5 (1) is optimal.

Theorem 10.5 is proved using the general method described in §10.2. A key ingredient is thus to etablish the right spectral estimates that allow control of the Littlewood-Paley square function. The relevant spectral estimate are given as follows.

**Theorem 10.6. Uniform spectral estimates for positive–definite spherical functions on complex groups** [CN]. *For $G$ as in Theorem 10.5, the following holds for non-trivial positive-definite spherical functions $\varphi_\lambda$, and for $H \in \overline{\mathfrak{a}_+}$*

$$\left| \frac{d^k}{dt^k} \varphi_\lambda(\exp(tH)) \right| \leq$$

$$\leq C_k(G, H)(1 + \|\lambda\|)^{k-\gamma} \exp(-\kappa t \rho(H))$$

*Here $\gamma = \gamma_G$ depend only on $G$. Furthermore $\kappa$ depends only on $G$ and is strictly positive, provided $G$ has no $SL_2(\mathbb{C})$-factors (equivalently, $G$ is complex and satisfies property $T$). Otherwise $\kappa$ is still strictly positive but depends also on $\lambda$.*

Theorem 10.6 implies that the $k$-th Littlewood-Paley square function associated with the operators $\sigma_t^H$ satisfies $\|R_k(f, \cdot)\|_2 \leq C_k \|f\|_2$, as long as $k \leq \gamma$. Thus the spectral method of §10.2 applies, and this proves Theorem 10.5.

The proof of Theorem 10.6 uses the explicit form of the spherical functions given by Harish Chandra's formula in the complex case. It is based on a method of descent, which allows writing a spherical function on $G$ as a sum of multiples of spherical functions on lower-dimensional complex subgroups. Such a decomposition is defined for every $\lambda \in \mathfrak{a}^*$, and the multiples that occur depend on $\lambda$. Choosing the optimal decomposition for each $\lambda$, a rate of decay in $\lambda$ is obtained, which depends only on the root system. We remark that the classical Harish Chandra estimate only provides an estimate which amounts to *polynomial growth* in the spectral parameter, rather than *polynomial decay* as in Theorem 10.6. Thus the Harish Chandra estimate cannot be used to establish the uniform spectral estimates necessary in order to bound the Littlewood-Paley square functions, that the spectral method presented above calls for.

The exponential decay in $t$ established in Theorem 10.6 holds uniformly for all positive-definite spherical functions, and is a consequence of the results of M. Cowling and R. Howe on matrix coefficients of unitary representations with a spectral gap of semisimple groups, which we will describe further in the next chapter.

Let us now turn to discuss maximal inequalities and ergodic theorems in the completely different set-up of totally disconnected lcsc groups, and some of their lattice subgroups. In the next section we consider groups of tree automorphisms, and in 11.3 we will discuss higher rank groups and lattices.

## 10.5. Radial structure on lattice subgroups : A generalization of Birkhoff's theorem.
The general spectral method presented in section 10.2 has some remarkable applications to certain countable groups which are not in themselves Gelfand pairs. In particular these application leads to a very natural generalization of Birkhoff's pointwise ergodic theorem as well of Hopf's maximal inequality. To explain it, consider the regular tree $T_k$ of constant valency $k \geq 3$, and its group of graph automorphisms $G = Aut(T_k)$, which acts transitively on the tree. The group $K$ of automorphisms stabilizing a given vertex $o$ is compact, and $T_k$ can be identified with $G/K$. The orbits of $K$ in $T_k = G/K$ are precisely the spheres $S_n(o)$ with $o$ as a center. A simple direct computation shows that if $\mathcal{S}_n$ denotes the operator of averaging a function on the tree $T_k$ on a sphere of radius $n$, then this sequence of operators satisfy the recurrence relation $\mathcal{S}_1 \mathcal{S}_n = \frac{1}{k} \mathcal{S}_{n-1} + \frac{k-1}{k} \mathcal{S}_{n+1}$. It follows immediately that each $\mathcal{S}_n$ is a polynomial in $\mathcal{S}_1$, and thus the algebra generated by these operators is cyclic and commutative. It is easy to see that each operator $\mathcal{S}_n$ can be identified with a (right) convolution operator on $G$, namely with convolution by the double coset $KgK$ corresponding the sphere $S_n(o)$ (a $K$-orbit in $G/K$). Thus it follows that $(G, K)$ is a Gelfand pair, and since $L^1(G, K) = M(G, K)$ is cyclic, and satisfies the second order recurrence relation with constant coefficients given by the identity above, it is not difficult to give its characters and $*$-characters in explicit form, as we shall see below.

Now note that the when $k = 2r$ is even, the Cayley graph $X(\mathbb{F}_r, S)$ of the free group on $r$ generators determined by a set $S$ of $r$ free generators and their inverses, has the structure of a $2r$-regular tree. The free group in question acts as a group $\Gamma$ of automorphisms of the Cayley graph, via its action by left translations and thus embeds in $G$. Furthermore the $\Gamma$-action is of course simply transitive, and $\Gamma$ intersects trivially with the stability group of every vertex. It is easily seen that the operator $\mathcal{S}_n$ of averaging on a sphere with radius $n$ in the Cayley graph can be given as an operator of (right) convolution on the group $\Gamma$, namely convolution by the measure $\sigma_n$, the uniform probability measure on a sphere of radius $n$ in the free group, w.r.t. the word metric defined by $S$. Indeed :

$$\rho_\Gamma(\sigma_n)f(w) = f * \sigma_n(w) = \frac{1}{|S_n|} \sum_{|x|_S = n} f(wx^{-1}) = \frac{1}{|S_n|} \sum_{d(w,u)=n} f(u) = \mathcal{S}_n f(w) \ .$$

Thus we conclude that the algebra generated by the averaging operators on spheres in the tree $T_{2r}$ embeds as a *commutative subalgebra* $\mathcal{A}$ of the convolution algebra $\ell^1(\mathbb{F}_r)$ (which in itself is of course, as non-commutative as a group algebra can be). It follows that every unitary representation of $\mathbb{F}_r$ gives rise to a $*$-representation of the commutative algebra $\mathcal{A}$, and thus the general spectral method explained in §10.2 applies.

Note that in fact the identification explained above of the algebra generated by the operators of averaging on spheres in the Cayley graph, with the algebra generated in $\ell^1(\Gamma)$ by the operators of right convolution by $\sigma_n$ holds more generally for every discrete group $\Gamma$ with generating set $S$. Thus the method of §10.2 applies whenever the automorphism group $G$ of the Cayley graph, together with the stability group $K$ of a vertex form a Gelfand pair.

Continuing with the case of the free groups, we recall that the spectrum of $\mathcal{A}(\mathbb{F}_r)$ is given as follows (where we define $q = 2r - 1$). The solutions to the second order recurrence relation satisfied in the algebra are

$$\varphi_z(\sigma_n) = \mathbf{c}(z)q^{-nz} + \mathbf{c}(1-z)q^{-n(1-z)}, \quad z \neq \frac{1}{2} + \frac{ij\pi}{\log q},$$

$$\mathbf{c}(z) = \frac{q^{1-z} - q^{z-1}}{(q+1)(q^{-z} - q^{z-1})}$$

and:

$$\varphi_z(\sigma_n) = \left(1 + n\frac{q-1}{q+1}\right)(-1)^{jn}q^{-\frac{n}{2}}, \quad z = \frac{1}{2} + \frac{ij\pi}{\log q}$$

A necessary and sufficient condition for $\varphi_z$ to be a continuous character (on the closed subalgebra of $\ell^1(\mathbb{F}_r)$ generated by the spheres) is that it be bounded, and this condition is equivalent to $0 \leq \operatorname{Re} z \leq 1$. The unitary representation of $\mathbb{F}_r$ in $L^2(X)$, extended to $\ell^1(\mathbb{F}_r)$, assigns to $\sigma_1$ a self-adjoint operator. Consequently, the values $\varphi_z(\sigma_1) = \gamma(z)$ are real, for those $\varphi_z$ that occur in the spectrum of $\sigma_1$ in $L^2(X)$. It is easily verified that $\gamma(z)$ is real iff $\operatorname{Re} z = \frac{1}{2}$, or $\operatorname{Im} z = \frac{ij\pi}{\log q}$.

The image of this set under $\gamma$ is the $*$-spectrum of $A(\mathbb{F}_r)$. Note that for $z$ and $1 - z$ the same character obtains, so we can assume that $0 \leq \operatorname{Re} z \leq \frac{1}{2}$. Note also that the characters corresponding to $z = s$ and to $z = s + \frac{ij\pi}{\log q}$ differ by sign only: $\varphi_{s + \frac{ij\pi}{\log q}}(\sigma_n) = (-1)^{jn} \varphi_s(\sigma_n)$. In particular, the sign character $\varepsilon$, given by $\varepsilon(\sigma_n) = (-1)^n$, is obtained at the points $z = \frac{i(2j+1)\pi}{\log q}$.

The $*$-spectrum is naturally divided to the principal series characters $\varphi_z$ where $\operatorname{Re} z = \frac{1}{2}$, and the complementary, where $z = s + \frac{ij\pi}{\log q}$, where $s$ is real and $0 \leq s \leq \frac{1}{2}$.

The comparison with the case of the spherical functions on the group $SL_2(\mathbb{C})$ is evident (see §9.3(6)), save of course for the fact that $\mathcal{A}(\mathbb{F}_r)$ has a unit, and thus its spectrum is compact.

In order to demonstrate some of the phenomena that arise here, let us consider in addition to the free group also the groups $\Gamma(r, h) = G_1 * G_2 \cdots * G_r$, the free product of $r$ finite groups each of order $h$, where $r \geq 2$, $h \geq 2$, $r + h > 4$, with generating set $S = \bigcup_{i=1}^{r} G_i \setminus \{e\}$. Here we define $q(\Gamma(r, h) = (r - 1)(h - 1)$.

The sphere averages on $\Gamma(r, h)$ also commute, and also satisfy a second order recurrence relation with constant coefficients. Indeed, again this algebra is closely related to a Gelfand pair, this time associated with the group of automorphisms of a semi-homogeneous tree (see [N0] and [N1] for more details).

The spectrum of $\Gamma(r, h)$ can thus be analyzed similarly, with one significant difference compared to that of $\mathbb{F}_r$. At the point $i\pi / \log q = i\zeta$ the special character that obtains is of the form (see [N1, §2.3])

$$\varphi_{i\zeta}(\sigma_n) = c_{i\zeta}(-1)^n + c_{1-i\zeta}(-1)^n q^{-n}.$$

Here $c_{i\zeta} = 1$ and $c_{1-i\zeta} = 0$ if and only if $h = 2$. Indeed $c(i\zeta) = \frac{q - (r-1)^{-1} - h + 2}{r(h-1)(1 - q^{-1})}$, and so $c(i\zeta) = 1$ iff $\Gamma = \Gamma(r, 2)$ or $\Gamma = \mathbb{F}_k$. We define $c_{i\zeta} = c(\Gamma)$.

Given a probability-preserving action of one of the groups $\Gamma(r, h)$ or $\mathbb{F}_r$, Let us denote by $\mathcal{E}'$ the orthogonal projection on the subspace of $L^2(X)$ given by $\ker(\pi(\sigma_1) - \varphi_{i\zeta}(\sigma_1)I)$. Note that the latter subspace consists of functions satisfying: $\pi(\sigma_n)f = (-1)^n f$ in the case of the free groups and $\Gamma(r, 2)$, but not otherwise.

As usual let $\mathcal{E}$ be the projection on the space of $\Gamma$-invariant functions ($\mathcal{E}f = \int_X f \, dm$ in the ergodic case). Recall that we defined $\sigma_n' = \frac{1}{2}(\sigma_n + \sigma_{n+1})$, and let $(X, m)$ be a $\Gamma$-space where $\mathcal{E}' \neq 0$.

We can nor formulate the following result.

**Theorem 10.7. Ergodic theorems for radial averages on free products**[N1].
*For $\Gamma = \Gamma(r, h)$ or $\Gamma = \mathbb{F}_k$ and $(X, m)$ (as above), the following holds :*

(1) $\sigma_n$ and $\beta_n$ *satisfy the maximal inequality in $L^2(X)$, but are not mean (and hence not pointwise) ergodic sequences in $L^2(X)$.*

(2) *The sequences $\sigma_n'$ is a pointwise ergodic sequence in $L^2(X)$.*

(3) $\sigma_{2n}$ *converges to $\mathcal{E} + c(\Gamma)\mathcal{E}'$, which is a conditional expectation operator w.r.t. a $\Gamma$-invariant sub-$\sigma$-algebra iff $\Gamma = \Gamma(r, 2)$ or $\Gamma = \mathbb{F}_k$.*

(4) $\beta_{2n}$ *converges to $\mathcal{E} + c(\Gamma)\frac{q(\Gamma) - 1}{q(\Gamma) + 1}\mathcal{E}'$, which is not a conditional expectation operator on a $\Gamma$-invariant sub-$\sigma$-algebra.*

*The convergence is for each function $f \in L^2(X)$, pointwise almost everywhere and in the norm of $L^2(X)$.*

*Remark* 10.8. We have chosen to focus on $L^2$ for simplicity of exposition. The strong maximal inequality for $\sigma_n$ (and thus $\sigma'_n$ and $\beta_n$) holds in fact in every $L^p$, $1 < p < \infty$, as can be verified using the argument in [NS1].

*Remark* 10.9. **The ball averaging problem : some counterexamples.**

Recalling the ball averaging problem in ergodic theory stated in §4.3, we see that Theorem 10.7 paints a rather complicated picture of the possibilities. First, note that the ball averages $\beta_n$ do *not* form a mean (and of course, pointwise) ergodic sequences, and neither does the sequence of spheres $\sigma_n$. This fact was originally observed in [Be1], by constructing an ergodic action of $\mathbb{F}_2$ where the special character $\sigma_n \mapsto (-1)^n$ of $\mathcal{A}$ is realized by an joint eigenfunction of the algebra $\mathcal{A}$ acting in $L^2$. Thus the ball average problem does not have a positive solution for general groups. However, Theorem 10.7 shows that this periodicity phenomenon is the only obstruction, and in fact $\sigma_{2n}$ and $\beta_{2n}$ do converge pointwise. Furthermore, note that the limit of $\sigma_{2n}$ is the conditional expectation w.r.t. the $\Gamma$-invariant $\sigma$-algebra of sets invariant under a subgroup of index at most two in the case where $\Gamma$ is a free group, but this is not the case for the groups $\Gamma(r, h)$, $h > 2$. For the ball averages $\beta_{2n}$ the limit is not a projection operator (and thus not a conditional expectation) at all. Thus it appears that the identification of the exact possible limits of subsequences of $\sigma_n$ and $\beta_n$ is a rather delicate problem, which seems inaccessible at this time, for general word-hyperbolic groups, say. The only exception thus far are those groups for which the spectral considerations above (or some variants, see [N1]) apply.

The situation just described is of course a reflection of non-amenability, since the balls do not have the Følner property of asymptotic invariance, and so a limit of a subsequence of $\beta_n$ need not posses any invariance properties w.r.t. the $\Gamma$-action.

We remark that it was conjectured in [Bu3, §9] that for every word-hyperbolic group $\Gamma$, and every symmetric set of generators the averages $\sigma_{2n}f$ converge pointwise to a function $f'$ invariant under the group $\Gamma_2$ generated by all words of even length. However consider the groups $\Gamma = \Gamma(r, h)$, $h > 2$, and a function $f \in L^2(X)$ which realizes the character $\varphi_{i\zeta}$ of the algebra $\mathcal{A}$. Then according to the formula above for $\varphi_{i\zeta}$

$$\lim_{n \to \infty} \pi(\sigma_{2n})f(x) = \lim_{n \to \infty} \varphi_{i\zeta}(\sigma_{2n})f(x) = c_{i\zeta}f(x) \,.$$

Thus the limit does exist, but $f$ is not a function invariant under $\Gamma_2$, being an eigenfunction of each $\pi(\sigma_{2n})$ with eigenvalue different than 1.

*Remark* 10.10. **Generalization of Birkhoff's theorem.**

Given an arbitrary invertible measure preserving transformation $T$ on a probability space $X$, Birkhoff's pointwise ergodic theorem asserts that for any $f \in L^1(X)$, the averages of $f$ along an orbit of $T$, namely the expressions $\frac{f(T^{-n}x) + \cdots + f(T^n x)}{2n+1}$ converge, for almost all $x \in X$, to the limit $\tilde{f}(x)$, where $\tilde{f}$ is the conditional expectation of $f$ w.r.t. the $\sigma$-algebra of $T$-invariant sets. Part of our quest to establish ergodic theorems for group actions can thus be motivated by the following obvious and natural question which presents itself. Given *two* arbitrary invertible measure preserving transformations $T$ and $S$, find a geometrically natural way to average a function $f$ along the orbits of the group generated by $T$ and $S$, so as to obtain the same conclusion.

Of course if $T$ and $S$ happen to commute, then, according to the discussion in Chapter 5, the expressions $\frac{1}{(2n+1)^2} \sum_{-n \le n_1, n_2 \le n} f(T^{n_1} S^{n_2} x)$ converge for almost all $x \in X$, for any $f \in L^1(X)$, and again the limit is the conditional expectation of $f$ w.r.t. the $\sigma$-algebra of sets invariant under $T$ and $S$. In other words, the pointwise ergodic theorem holds for finite-measure-preserving actions of the free Abelian group on two generators, namely $\mathbb{Z}^2$. However, it is clear that when choosing generically two volume-preserving diffeomorphisms of a compact manifold, or two orthogonal transformations of the Euclidean unit sphere, or in general two measure preserving maps of a given measure space, the group generated by them is not Abelian, and in fact, it is generically free.

The answer to the problem above is then to find an averaging sequence satisfying a pointwise ergodic theorem for finite-measure-preserving actions of the free non-Abelian group on two generators. The first choice that one would consider by direct analogy with Birkhoff's and Wiener's theorems (for $\mathbb{Z}$ and $\mathbb{Z}^d$), would be the normalized ball averages w.r.t. a set of free generators. For the free group this problem has been settled, using the spectral methods described above, by the following result.

**Theorem 10.11. Generalization of Birkhoff's theorem**[NS1]. *Consider the free group $\mathbb{F}_r$, $r \ge 2$, with symmetric free generating set $S$. Let $(X, m)$ be a probability-preserving ergodic action. Then*

(1) *The sequence $\sigma'_n = \frac{1}{2}(\sigma_n + \sigma_{n+1})$ satisfies the strong maximal inequality and is a pointwise ergodic sequence in $L^p$, for $1 < p < \infty$.*

(2) *The sequence $\beta_n$ satisfies the pointwise ergodic in $L^p$, $1 < p < \infty$ if and only if $L^2(X)$ does not contains a non-zero function $f_0$ satisfying $\pi(w)f_0 = (-1)^{|w|} f_0$ for every $w \in \mathbb{F}_r$.*

(3) *If such an eigenfunction $f_0$ is present then it is unique, has constant absolute value, and $\beta_n f_0$ does not converge. For any $f \in L^p(X)$, $1 < p < \infty$*

$$\lim_{n \to \infty} \beta_{2n} f(x) = \int_X f \, dm + \frac{r-1}{r} \int_X f \overline{f_0} \, dm \cdot f_0(x)$$

*pointwise and in the $L^p$-norm.*

**Problem 10.12.** We note that the weak-type maximal inequality in $L^1$ for the sphere averages (or equivalently, the ball averages) on the free group is an open problem.

*Remark* 10.13. Note that the pointwise ergodic theorem for the free finitely generated group $\mathbb{F}_k$ implies a corresponding one for any factor group of $\mathbb{F}_k$, namely for any finitely generated group, just as Wiener's theorem for $\mathbb{Z}^d$ implies the corresponding result for any finitely generated Abelian group. However, the weights that must be taken on the factor group are those induced by the canonical factor map, and these usually bear little resemblance to the intrinsic ball and sphere averages on the factor group.

*Remark* 10.14. **Sphere averages on free algebras.** We note that the same spectral methods that were employed for sphere averages in the group algebra of the free group can be used more generally for other free algebras in various varieties. For example, consider the free associative algebra on $r$ non-comutative elements.

This algebra has of course a natural length function, and clearly the algebra of radial elements is commutative (under convolution), and satisfy a first-order recuurence relation. It is thus a simple exercise to develop the spectral theory of the ∗-representations of the subalgebra of radial elements. These include representation where each generator is mapped to bounded self-adjoint contraction on a Hilbert space, and we thus obtain a pointwise ergodic theorem for the powers of the self-adjoint operator which is the uniform average of the $r$ contractions. . Similarly, we can consider the free algebra on $r$ non-commuting idempotents, where again we have a commutative subalgebra of radial elements. The ∗-spectrum can be determined here too, again by solving a second-order recurrence relation with contant coefficients. The ∗-representations include those where each generator is mapped to a self-adjoint projection. We thus obtain in particular a pointwise ergodic theorem for the radial averages of (non-commuting) conditional expectations on a probability space.

In Theorem 10.7 we already considered the convolution algebras of the free products $\mathbb{Z}_p * \cdots * \mathbb{Z}_p$ which are the free groups generated by $r$ elements of order $p$. The ∗-representations here are given by the unitary representations of the groups, and the spherical functions can again be explicitly determined from a second order recurrence relation (see e.g. [N1] for their description).

There are further examples of free algebras in other varieties, where the radial elements form a commutative subalgebra whose ∗-spectrum can be determined using a recurrence relation. In all of these cases a mean and pointwise ergodic theorem for sphere averages in ∗-representations is obtained, with the sole obstruction given by periodicity phenomena, when they occur.

*Remark* 10.15. **Boundary transitive subgroups of tree automorphisms.** Let us note that the group algebra of the algebraic group $PGL_2(\mathbb{Q}_p)$ also contains an isomorphic copy of the algebra $\mathcal{A}$ of even radial averages on the tree $T_{p+1}$. This follows from the fact that the group has a faithful representation as a group of automorphisms of the regular tree, which is transitive on the boundary. Again the averaging operators on the tree can be represented as convolution operators on $PGL_2(\mathbb{Q}_p)$. However, here the Howe-Moore mixing theorem [HM] applies, namely the matrix coefficients of unitary representations without invariant unit vectors vanish at infinity for an algebraically connected simple group. This implies that the sign character of the algebra $\mathcal{A}$ cannot appear in ergodic actions of $PSL_2(\mathbb{Q}_p)$, and both $\beta_{2n}$ and $\sigma_{2n}$ converge pointwise and in norm to the ergodic mean, in every $L^p$, $1 < p < \infty$, as follows from [N1] and [NS1]. A similar analysis holds for many closed non-compact boundary-transitive subgroups of the group of automorphisms of a semi-homogeneous tree [N1]. The analog of the Howe-Moore theorem here was proved in [LM].

*Remark* 10.16. **Non-commutative Hecke algebras** We note that further natural algebra structures appear in the setting of groups of automorphisms of semi-homogeneous tree, and more generally groups of automorphisms of the Bruhat-Tits buildings of semisimple algebraic groups over locally compact totally disconnected non-discrete fields. These are the Hecke algebras $L^1(Q \backslash G/Q)$ of double cosets of a compact open subgroup $Q$ under convolution, which are non-commutative in general, but have the property that all of their irreducible ∗-representations have a uniformly bounded degree [Ber]. Particularly significant among these algebras is the Iwahori algebra, consisting of double cosets of an Iwahori subgroup. In the case of the group $Aut(T_{r_1,r_2})$, $(r_1 \neq r_2)$ for example, the Iwahory subgroup can be identified with the stability group of an edge, and thus it is contained in the two non-conjugate maximal compact open subgroup stabilizing one of the vertices of the edge. Thus the commutative algebra associated with a Gelfand pair structure on the group is contained as a subalgebra of the Iwahori algebra in this case (and others). The spectrum of the Iwahori algebra on a semi-homogeneous tree can be easily determined (see e.g. [N0]), and it appears naturally when analyzing the spectrum of some natural convolution algebras in certain lattices of $Aut(T_{r_1,r_2})$ $(r_1 \neq r_2)$. These include for example the groups $\mathbb{Z}_p * \mathbb{Z}_q$, and so also the group $PSL_2(\mathbb{Z}) = \mathbb{Z}_2 * \mathbb{Z}_3$. Thus it is possible also to use non-commutative harmonic analysis on Hecke algebras to derive ergodic theorems for discrete groups. This possibility was explored in the case of semi-homogeneous trees in [N1], but it is natural to expect that it can be developed much further.

*Remark* 10.17. We note that an alternative proof of the ergodic theorems for spheres on the free group was developed by A. Bufetov [Bu3]. The method is based on the theory of Markov processes rather than on spectral theory, and will be discussed further in §12.4, together with some other ergodic theorems on free groups and other Markov groups.

## 11. Actions with a spectral gap

We now turn to a discussion of a fundamental phenomenon that appears in the study of non-amenable algebraic groups, and which has no analog in the theory of amenable groups. Utilizing it, we will be able to greatly expand the scope of the radial ergodic theorems on semisimple algebraic groups, obtain quantitative estimates in the pointwise ergodic theorems, and also obtain a host of results on a diverse array of non-radial averages.

The phenomenon in question is the existence of properly ergodic (i.e., non-transitive) actions with a spectral gap. We define the latter property in the form most convenient for our purposes here, as follows.

**Definition 11.1. Spectral gaps**.

(1) A strongly continuous unitary representation $\pi$ of an lcsc group $G$ is said to have a spectral gap if $\|\pi(\mu)\| < 1$, for some (or equivalently, all) absolutely continuous symmetric probability measure $\mu$ whose support generate $G$ as a group.

(2) Equivalently, $\pi$ has a spectral gap if the Hilbert space does not admit an asymptotically-$G$-invariant sequence of unit vectors, namely a sequence satisfying $\lim_{n \to \infty} \|\pi(g)v_n - v_n\| = 0$ uniformly on compact sets in $G$.

(3) A measure preserving action of $G$ on a $\sigma$-finite measure space $(X, m)$ is said to have a spectral gap if the unitary representation of $G$ in the space orthogonal to the space of $G$-invariant functions has a spectral gap. Thus in the case of an ergodic probability-preserving action, the representation in question is on the space $L_0^2(X)$ of function of zero integral.

(4) An lcsc group $G$ is said to have *Kazhdan's property T* [Kaz] provided every strongly continuous unitary representation which does not have $G$-invariant unit vectors has a spectral gap.

*Remark* 11.2.       (1) The equivalence between (1) and (2) is a standard argument in spectral theory and can be found e.g. in [MNS].

(2) The phenomenon of spectral gaps does not occur for properly ergodic probability-preserving actions of amenable groups. Indeed, in any such action $X$, there exists a non-trivial sequence of sets $A_n \subset X$ whose measures satisfy $0 < c < m(A_n) < C < 1$ for all $n$, which is asymptotically invariant, namely $\lim_{n \to \infty} m(gA_n \Delta A_n) = 0$ uniformly on compact sets in $G$ [Ro]. It follows immediately that $\|\pi(\mu)\|_{L^2_0(X)} = 1$ for *every* probability measure $\mu$ on $G$.

(3) We recall the well-known fact that amenability of an lcsc group (namely the existence of a Følner sequence) can be characterized by the condition that the left regular representation $\lambda_G$ satisfies $\lambda_G(\mu) = 1$, for at least one (or equivalently, all) absolutely continuous symmetric probability measure $\mu$ whose support generate $G$ as a group.

In the realm of non-amenable algebraic groups, actions preserving a $\sigma$-finite measure which have a spectral gap are quite ubiquitous, and we briefly indicate some examples.

**Example 11.3. Some examples of actions with a spectral gap**.

(1) $G$ a non-amenable group, $X = G$, and the action is by left translations, w.r.t. Haar measure.

(2) $G$ a connected semisimple Lie group, and the action is by isometries of the Riemannian symmetric space $G/K$. More generally, the action on a reductive symmetric space of the form $G/H$, where $H$ is the fixed-point-group of an involutive automorphism of $G$.

(3) $G$ a simple algebraic group, and the action by translation on the homogeneous space $G/L$, where $L$ is a proper algebraic unimodular subgroup [Gu1], e.g. the action of $SL_n(\mathbb{R})$ on the space of symmetric matrices, or the action of $Sp(n, \mathbb{R})$ on $\mathbb{R}^{2n}$.

(4) Similarly, for $G$ a simple algebraic group the action by translations on the homogeneous space $G/\Lambda$ where $\Lambda \subset G$ is a discrete subgroup which is not Zariski dense (see e.g. the discussion in [IN]).

(5) $\tau : G \to H$ a representation of a semisimple algebraic group $G$ in a simple algebraic group $H$, $\Delta \subset H$ a lattice subgroup, and $G$ acts on $X = \Delta \setminus H$, a locally symmetric space of finite volume via $\tau$.

(6) If $G$ is a simple algebraic group of split-rank at least two then $G$ has property $T$ (see e.g. [M]), and then of course *any* unitary representation without invariant unit vectors has a spectral gap.

**11.1. Pointwise theorems with exponentially fast rate of convergence.** The utility of a spectral gap in a given representation $\pi$ is in the fact that typically, given a natural family $\mu_t$ of probability measures on $G$, not only do we have $\|\pi(\mu_t)\| < 1$ for each $t > 0$, but in fact (as we shall see presently) the far stronger conclusion that the norms decay exponentially in $t$ holds, namely :

$$\|\pi(\mu_t)\| \leq C_\mu \exp(-\delta_\mu t),$$

where $\delta_\mu > 0$ depends on the family $\mu_t$.

When $G = SL_2(\mathbb{C})$ for example, the spectral gap condition is equivalent with the the condition that the $*$-spectrum that arise in the representation $\pi$ of $M(G, K)$ contains only complementary series characters $\varphi_a$ with parameter satisfying $a \leq 1 - \theta$, $\theta = \theta(X) > 0$ (except for the trivial character). The exponential decay of the operator norms of $\sigma_t$ for example follows immediately upon evaluating the characters $\varphi_z$ on the sphere $S_t$. It also follows easily for $\gamma_t$ and $\beta_t$ upon integrating $\varphi_z$ against these measures, as we will see below.

The exponentially decaying norm estimate above is a most useful fact, which as we shall see gives rise to an interesting new phenomenon in ergodic theory, namely pointwise ergodic theorems with an explicit exponentially fast rate of convergence to the ergodic mean, for properly ergodic actions. The validity of the norm estimate follows from spectral estimates that we will consider in more detail in §11.2. But to illustrate the pooint, we now turn to our first use of such estimates, namely the following pointwise ergodic theorem with an error term for the bi-$K$-invariant averages $\beta_t$ on a simple Lie group, and to its proof. We recall that the averages $\beta_t$ we shall consider are defined as the $K$-invariant lifts to $G$ of the $K$-invariant probability measures on balls $B_t$ w.r.t. the Killing form on the symmetric space $G/K$, with center $[K]$ and radius $t$. This family generalizes the case of hyperbolic space considered in Chapter 9.

**Theorem 11.4. Pointwise ergodic theorem with exponentially fast rate of convergence for ball averages on semisimple Lie groups in actions with a spectral gap** [MNS]. *Let $G$ be connected non-compact semisimple Lie group with finite center, and let the $G$-action on $X$ have a spectral gap. Then the ball averages $\beta_t$ converge pointwise exponentially fast to the ergodic mean. More precisely $\forall f \in L^p(X)$, $p > 1$ and for almost every $x \in X$ :*

$$\left| \beta_t f(x) - \int_X f \, dm \right| \leq B_p(f, x) \exp(-\theta_p t) \ , \ \theta_p > 0 \ .$$

*where $\theta_p$ depends on $p$, $G$ and $X$.*

*Furthermore, the integer-radius sphere averages $\sigma_n$ also converge pointwise exponentially fast to the ergodic mean.*

<u>**Model case**</u>: **Proof for ball averages $\beta_t$ in ergodic actions of $SL_2(\mathbb{C})$ with a spectral gap.**

Recall that we are assuming here that the $*$-spectrum determined by the representation $\pi_0$ of $G = SL_2(\mathbb{C})$ and $M(G, K)$ on $L^2_0(X)$ (the space of functions with zero integral), contains only complementary series characters $\varphi_a$ with parameter satisfying $a \leq 1 - \theta$, $\theta = \theta(X) > 0$. We denote the spectrum by $\Sigma^*(\pi_0)$. The proof proceeds along the following steps.

(1) Since $a \leq 1 - \theta < 1$, we have, using the explicit form of the characters (see §9.3) :

$$|\varphi_a(\sigma_t)| \leq B \exp(-\theta t)$$

and hence we have the following exponential decay estimate on the norm of the sphere averages

$$\|\sigma_t\|_{L^2_0 \to L^2_0} = \sup_{z \in \Sigma^*(\pi_0)} |\varphi_z(\sigma_t)| \leq B \exp(-\theta t) \ .$$

(2) Similarly, for the ball averages, using the estimate of their density :

$$\|\beta_t\|_{L^2_0 \to L^2_0} \leq C \exp(-2t) \int_0^t \exp(2s) \|\sigma_s\|_{L^2_0 \to L^2_0} \, ds \leq$$

$$\leq C B \exp(-2t) \int_0^t \exp((2 - \theta)s) ds \leq C' \exp(-\theta t) \ .$$

(3) For integer radius balls (and similarly, spheres) acting in $L_0^2$, we have :

$$\sum_{n=0}^{\infty} \left\| \exp(\frac{\theta}{2}n)\beta_n f \right\|_2^2 \le \sum_{n=0}^{\infty} \exp(-\theta n) \|f\|_2^2 < \infty$$

equivalently :

$$\sum_{n=0}^{\infty} \int_X \left| \exp(\frac{\theta}{2}n)\beta_n f(x) \right|^2 dm < \infty$$

and so

$$\sup_{n \ge 0} \left| \exp(\frac{\theta}{2}n)\beta_n f(x) \right|^2 \le$$

$$\le \sum_{n=0}^{\infty} \left| \exp(\frac{\theta}{2}n)\beta_n f(x) \right|^2 = C(f,x)^2 < \infty$$

For almost all $x \in X$.

(4) We can therefore conclude :

A. $\beta_n$ (and $\sigma_n$) satisfies the *exponential-maximal inequality* in $L_0^2(X)$ stated in (3), and

B. $\beta_n$ (and $\sigma_n$) *converges pointwise exponentially fast* to the ergodic mean for $f \in L^2$, i.e. for almost every $x \in X$ :

$$\left| \beta_n f(x) - \int_X f dm \right| \le C(f - \int_X f dm, x) \exp(-\frac{\theta n}{2}) .$$

(5) The natural generalization of statements A and B in (4) above are of course true in every $L_0^p$, $p > 1$, by the Riesz-Thorin interpolation theorem.

(6) Now note that the foregoing arguments show in fact a more general result, namely that the same conclusion holds for every sequence $\beta_{t_k}$ with $\sum_{k \in \mathbb{N}} \exp(-\frac{1}{2}\theta t_k) < \infty$. Fix such a sequence $t_k$, and given a point $t$, choose the closest point to it (which we assume is at a distance at most one), and denote it by $t_n$. Now write *for $f$ bounded* :

$$\left| \beta_t f(x) - \int_X f dm \right| \le$$

$$\le |\beta_t f(x) - \beta_{t_n} f(x)| + \left| \beta_{t_n} f(x) - \int_X f dm \right|$$

The second term is bounded (using part (3) above, and replacing the integers by the sequence $t_k \in \mathbb{R}$) by

$$C'(f - \int_X f dm, x) \exp(-\frac{1}{2}\theta t_n)$$

where

$$\left\| C'(f - \int_X f dm, \cdot) \right\|_2 \le B \|f\|_2 \le B \|f\|_\infty .$$

As to the first term, note that the family $\beta_t$ is uniformly (locally) Lipschitz continuous (w.r.t. the $L^1(G)$-norm), and the function $f$ is bounded. It follows that the first term is bounded by $|t - t_n| \|f\|_\infty$.

(7) Let us choose the sequence $t_k$ fine enough, so that it satisfies $|t - t_n| \leq \exp(-\frac{1}{4}\theta n)$. This can be achieved by dividing the interval $[n, n+1]$ to $2 + [\exp(\theta n/4)]$ equally spaced points, and the resulting sequence still satisfies the condition $\sum_{k \in \mathbb{N}} \exp(-\frac{1}{2}\theta t_k) < \infty$ stated in (6). We can then use the argument of (6) for the sequence $t_k$ to estimate both the first and the second term.

The estimate in (6) immediately implies the following, which we call an $(L^\infty, L^2)$-exponential maximal inequality :

$$\left\| \sup_{t>0} \exp(\frac{1}{4}\theta t) \left| \beta_t f(x) - \int_X f \, dm \right| \right\|_{L^2} \leq C_2 \|f\|_\infty$$

(8) We can now conclude that for every bounded function $f$ the expression

$$\sup_{t>0} \exp(\frac{1}{4}\theta t) \left| \beta_t f(x) - \int_X f \, dm \right| = B(f, x)$$

is finite almost everywhere, and the conclusion of Theorem 11.4 is thus established for $f$ in $L^\infty$.

(9) We now note that the strong maximal inequality for ball averages of Theorem 9.4 namely

$$\left\| \sup_{t>0} \beta_t f \right\|_p \leq C_p \|f\|_p \quad , \ p > 1$$

can be established very easily for balls with exact exponential growth in actions with a spectral gap. Indeed, clearly exact exponential volume growth implies that for a fixed constant $B$ we have $\beta_t \leq B\beta_{[t]+1}$, $t \geq 1$ as measures on $G$. Thus the maximal inequality for the family of all balls follows from its validity for the sequence of balls with integer radii. The boundedness of the latter maximal function is an elementary conclusion of exponential decay of the operator norm, as follows from an obvious variation on the arguments presented in part (3) and (4).

(10) Finally, we can use the analytic interpolation theorem again. This time we interpolate between the exponential maximal inequality stated in (7) (from $L^\infty$ to $L^2$) and the strong maximal inequality (from $L^p$ to $L^p$, $p > 1$) for the ball averages, stated in (9). We then obtain an $(L^p, L^r)$-exponential maximal inequality (in the obvious notation) and hence exponential pointwise convergence to the ergodic mean, for every $f \in L^p$, $1 < p < \infty$. We refer to [MNS] for the details. This concludes the outline of the proof of Theorem 11.4.

The proof above is complete for actions of $SL_2(\mathbb{C})$ with a spectral gap, and similar arguments also yields the general case of ball averages on semisimple groups in actions with a spectral gap. One uses spectral estimates of spherical functions, resulting in the exponential decay of the operator norm (see the discussion in the following section). Further, the monotonicity property $\beta_t \leq C\beta_{[t]+1}$ of the ball averages is valid here and follows from strict $t^q \exp ct$-volume growth of the balls, which is a relatively starightforward consequence of the structure theory of semisimple Lie groups. The monotonicity is utilized to deduce the strong maximal inequality in $L^p$, $1 < p < \infty$ of the operator $\sup_{t>0} \beta_t$ from the discrete version $\sup_{n \in \mathbb{N}} \beta_n$ (see more on this argument in §12). It is also necessary to establish that the ball averages are uniformly locally Lipchitz continuous in the $L^1(G)$-norm. We refer to [MNS][N4] for the details.

## 11.2. The spectral transfer principle.

Spectral estimates for spherical functions are crucial to the successful implementation of the spectral approach to maximal inequalities and pointwise ergodic theorems that was outlined in the preceding section. The uniform derivative estimates necessary for the proof of the maximal inequalities for spheres (and other singular averages) are usually very difficult and have not been established in general. However the basic exponential decay estimate for the operator norm of radial averages such as spheres and balls (acting in an irreducible unitary representation, say) hold in great generality and depend only on basic structural features that hold for all simple algebraic groups over locally compact non-discrete fields.

These estimates were developed in various forms by M. Cowling [Co1] and R. Howe [H], as well as C. C. Moore [HM], U. Haagerup [CHH] and Borel-Wallach [BW]. An elegant exposition to the case of $SL_n(\mathbb{R})$ appears in [HT], and the strategy outlined there was used by H. Oh [Oh][Oh1] to obtain definitive quantitative results for semisimple algebraic groups. We summarize some of these results as follows.

## Theorem 11.5. Decay estimate and $L^p$-integrability of matrix coefficients

(1) [Co1] *Let $G$ be a simple non-compact connected Lie group with finite center. For every irreducible non-trivial unitary representation $(\pi, \mathcal{H})$ of $G$, and every two $K$-finite vectors $u, v \in \mathcal{H}$, the associated matrix coefficient $\psi_{u,v}(g) = \langle \pi(g)u, v \rangle$ has the following two properties.*

  (a) *The matrix coefficient satisfies an exponential decay estimate along the group $G$ :*

$$|\psi_{u,v}(g)| = |\langle \pi(g)u, v \rangle| \leq C_{u,v} \exp(-\delta_\pi d(K, gK))$$

  *where $d$ denotes the distance function associated with the Riemmanian metric given by the Killing form on the symmetric space $G/K$.*

  (b) *The matrix coefficient $\psi_{u,v}(g) = \langle \pi(g)u, v \rangle$ belongs to $L^p(G)$, for some $p = p(\pi) < \infty$.*

(2) [Oh, Thm 1.1] *The same estimate holds for infinite-dimensional irreducible unitary representations of any simple algebraic group over a locally compact non-discrete field $F$ with Char $F \neq 2$, where $K$ is a good maximal compact subgroup of $G$ (see [Car, §3.5]) and $d(gK, K)$ the metric induced on $G/K$ by its inclusion in the Bruhat-Tits building of $G$.*

(3) [Co1][Oh1, §5.7] *When $G$ has property $T$, the same estimates hold in both cases for $K$-finite vectors in any unitary representation of $G$, which does not have $G^+$-invariant unit vectors.*

*Remark* 11.6.       (1) We note that in this set-up, $G$ can be defined to have property $T$ if and only if $p(\pi) \leq p(G) < \infty$ for all irreducible non-trivial unitary representations, i.e. the spectral estimate depends only on $G$ and holds uniformly for all the representations. Equivalently $\inf_\pi \delta_\pi > 0$ for all representations $\pi$ with a spectral gap.

(2) Property $T$ holds for all simple algebraic groups of split-rank at least two, but also for some simple real Lie groups of real-rank one - see [M] for a discussion.

(3) We remark that even more precise quantitative estimate for the decay of positive definite spherical functions have been developed and refer the reader to the reference cited above.

(4) Finally, $G^+$ is the group generated by the split unipotent subgroups of $G$ (see [M, Ch. I, §§1.5, 2.3] for a discussion). It is normal and co-compact in $G$, and is of finite index in $G$ whenever the characteristic of the field is zero. In particular, it coincides with $G$ when $G$ is a connected semisimple Lie group without compact factors.

Note that for any $n \geq p/2$, Theorem 11.5 implies that $\psi_{u,v}(g)^n \in L^2(G)$. This fact implies, much as in the Peter-Weyl theorem for compact groups, the following result due to M. Cowling in the real semisimple case, and R. Howe and C. C. Moore in general.

**Theorem 11.7. Spectral Transfer Principle**[Co1][HM]. *Let $G$ be as in Theorem 11.5. If the representation $\pi$ of $G$ has a spectral gap, there exists $n = n(\pi)$ such that*

$$\pi^{\otimes n} \subset \infty \cdot \lambda_G$$

*where $\infty \cdot \lambda_G$ denotes the direct sum of countably many copies of the regular representation of $G$, $\pi^{\otimes n}$ the $n$-fold tensor power of the representation $\pi$, and $\subset$ denote a unitary isomorphism onto a subrepresentation. If $G$ has property $T$ then there exists a uniform bound $n(\pi) \leq n(G) < \infty$ for all representations $\pi$ with spectral gap (and conversely).*

This result has the following explicit spectral estimate as a corollary :

**Theorem 11.8. Uniform norm estimate of arbitrary measures on $G$ in arbitrary representations with a spectral gap** [N4, Thm 1.1]. *Let $G$ be a group as in Theorem 11.5. Let $\pi$ be any unitary representation of $G$ with a spectral gap. Let $\mu$ be any probability measure on $G$. Then*

$$\|\pi(\mu)\| \leq \|\lambda_G(\mu)\|^{1/n(\pi)}$$

*In particular, if $(X, m)$ is a probability-preserving ergodic action of $G$ with a spectral gap, then*

$$\|\pi(\mu)\|_{L^2_0(X)} \leq \|\lambda_G(\mu)\|^{1/n(\pi_0)} \ .$$

Note that it follows easily from the spectral estimate of Theorem 11.5(1) that typically, given a family $\mu_t$ of radial probability measures on $G$ satisfying mild natural growth conditions, we have an *exponential decay estimate on the convolution norm* :

$$\|\lambda_G(\mu_t)\| \leq \exp(-\delta t) \ , \ \delta = \delta(\mu) > 0 \ .$$

In particular this visibly holds for sphere averages, balls averages, a host of variations of shell averages [N6], and many other radial averages $\mu_t$ on $G$.

As we shall see in §12, it is possible to use Theorem 11.8, in conjunction with the arguments of §11.1 to establish exponential-maximal inequalities for a wide class of non-radial families $\mu_t$, as well as exponential pointwise convergence to the ergodic mean. But before turning to non-radial averages let us mention an application of the spectral transfer principle in the radial case, namely establish the pointwise ergodic theorems in the case of totally disconnected simple algebraic groups.

**11.3. Higher rank groups and lattices.** We have encountered in §10.5 totally disconnected Gelfand pairs which appeared as groups acting on semi-homogeneous trees. A more general natural set of geometries to consider is that of affine Bruhat-Tits buildings, in particular those associated with semisimple algebraic groups over locally compact totally disconnected non-discrete fields. The case of totally disconnected simple algebraic groups of (split) rank one reduces to that of closed boundary-transitive group of automorphisms of semi-homogeneous trees. Thus the pointwise ergodic theorems for sphere and ball averages on them are completely resolved by Theorem 10.7 together with Remark 10.8 and Remark 10.15.

Consider now the case of a simple algebraic group over locally compact totally disconnected non-discrete field $F$ of split rank at least two. Such a group $G$ has a good compact open subgroup $K$ giving rise to a Gelfand pair structure (see e.g. [Car, §3.5 and Thm 4.1] and the references there). Furthermore, as follows immediately from Theorem 11.5, the spherical functions associated with $(G, K)$ decay exponentially with a uniform bound, when the group is simple of split rank greater than one. More precisely, positive-definite spherical functions associated with irreducible non-trivial representations are bounded by a function which decays exponentially fast to zero as a function of $d(gK, K)$, where $d$ is the $G$-invariant distance on the building, as is the case when $G/K$ is a Riemannian symmetric space. The same holds true for matrix coefficients of unitary representations without invariant unit vectors. It follows that the spectral norm of the convolution operators associated with the bi-$K$-invariant ball averages $\beta_n$ on $G$ decays exponentially in $n$. Together with the spectral transfer principle for these groups, which follows from the estimate just described, Theorem 11.8 implies that the norm of $\pi_0(\beta_n)$ in $L_0^2(X)$ for any ergodic action $X$ decays exponentially also, with a fixed bound, independent of $X$. Using the arguments brought in the proof of Theorem 11.4 (for integer radius only, this time !) the exponential norm decay estimate gives a proof of the following result.

**Theorem 11.9. Pointwise ergodic theorem with exponentially fast rate of convergence for simple algebraic groups**. *Let $G$ and $K$ be as in Theorem 11.5. Let $\beta_n$ be the bi-$K$-invariant ball averages and $\sigma_n$ the sphere averages. Assume $G$ has split rank at least two, or more generally, property $T$. Then for any action of $G$ on a probability space which is ergodic under $G^+$, for all $f \in L^p(X)$, $1 < p < \infty$ and almost every $x \in X$*

$$\left| \pi(\beta_n)f(x) - \int_X f \, dm \right| \leq C_p(f, x) \exp(-n\delta_p)$$

*where $\delta_p > 0$ depends only on $G$ and $p$. The same result holds also for the sphere averages $\sigma_n$. Furthermore $\|C_p(f, \cdot)\|_{L^p(X)} \leq B_p \|f\|_{L^p(X)}$, and so the maximal functions associated with the spheres (and balls) satisfy an exponential-maximal inequality in $L_0^p(X)$.*

*Finally, the same holds true for any action of a group as in Theorem 11.5, provided that it has a spectral gap, but with $\delta_p$ here depending also on the action.*

Thus we see that the ball and sphere averages converge exponentially fast to the ergodic mean from almost every starting point, with a fixed rate independent of the starting point as well as the action.

Recall that our discussion in §10.5 of totally disconnected Gelfand pairs extended also to some of their lattices. Note that the free group appeared in our discussion there as a group of automorphisms of the tree which acts simply transitively on the vertices. The fact that the group algebra of the free group contained an isomorphic copy of the algebra of radial averaging operators on the tree followed without difficulty.

The discussion is not limited of course to groups acting on semi-homogeneous trees, and this set-up can be expanded considerably, as follows. In [CMSZ] the authors have constructed groups acting simply transitively on the vertices of certain affine buildings of rank 2, namely the $\tilde{A}_2$-buildings. Furthermore, such groups have been constructed for $\tilde{A}_n$-buildings (which are of rank $n$) for any $n$ in [CaS]. As noted already, a simple algebraic groups $G$ always contains a good compact open subgroup such that $(G, K)$ is a Gelfand pair (see e.g. [Car]), and thus for appropriately chosen lattices one can expect to find an isomorphic copy of the Gelfand pair algebra $L^1(G, K)$ in the group algebra $\ell^1(\Gamma)$ of the lattice. In [CMS] the authors have succeeded in showing that such a commutative algebra does occur in the group algebras of some of the discrete groups constructed in [CMSZ]. Furthermore they analyzed its structure and $*$-representations, even in the case where the discrete groups do not arise as lattices in simple algebraic groups. An interesting new feature that arise here is the fact that the discrete groups in question satisfy property $T$, a fact that is proved in [CMS] directly from the representation theory of the commutative algebra in question. These results are very interesting in the context of ergodic theory, as they demonstarte the following rather remarkable phenomenon (based on the results of [CMS] and [RRS]).

**Theorem 11.10. Uniform pointwise ergodic theorem for some groups with property** $T$. *There exists a discrete group $\Gamma$ with a finite generating set $S$, with the following property. For some fixed positive $\delta > 0$ depending only on $(\Gamma, S)$, and every ergodic probability preserving action of $\Gamma$ on $(X, m)$, for all $f \in L^2$, and almost etevy $x \in X$*

$$\left| \pi(\sigma_n)f(x) - \int_X f dm \right| \leq C_2(f, x) \exp(-n\delta)$$

*In fact, there are infinitely many $\tilde{A}_2$ groups satisfying the property above. The same result holds of course in $L^p$, $1 < p < \infty$, with $\delta_p > 0$ (but this is an open problem in $L^1$ !).*

This phenomenon is of course, in striking contrast to the behaviour of averages on discrete amenable groups in classical ergodic theory. It raises the following intriguing problem, which however seems quite inaccesible at this time.

**Problem 11.11.** Does the phenomenon described in Theorem 11.10 occur for *every* discrete group with property $T$, and for every set of generators on it ?

## 12. Beyond radial averages

The spectral methods introduced in Chapter 11, together with some further arguments, can be used to prove a diverse variety of (not necessarily radial) pointwise ergodic theorems on semisimple Lie and algebraic groups. When the action has a spectral gap, it is possible to establish exponentially fast rate of convergence for general families of probability measures $\mu_t$, which are decidedly non-radial. Such families often have a great deal of intrinsic geometric interest and occur naturally in applications. We now turn to an exposition of some of these results, refering for more details to [N4] and [N6].

## 12.1. Recipe for pointwise theorems with rate of convergence.

12.1.1. *Estimating convolution norms.* When discussing the ergodic theory of a general family $\mu_t$ of probability measures on a semisimple Lie group $G$, the first step is to establish exponential decay estimates on the norms of the convolution operators $\lambda_G(\mu_t)$. This will be then converted by the spectral transfer principle (Theorem 11.7 and Theorem 11.8) to norm decay estimates in an arbitray measure-preserving action with a spectral gap.

There exists a fundamental estimate of the convolution norms when the averages $\mu_t$ are absolutely continuous, and it is given in terms of the $L^p(G)$-norms of the densities of $\mu_t$. The validity of such an estimate is called the Kunze-Stein phenomenon, established in [KSt] for $G = SL_2(\mathbb{R})$ and in [Co] in general. The precise formulation is as follows.

**Theorem 12.1. Kunze-Stein phenomenon [KSt][Co].**
   *Given a connected semisimple Lie group $G$ with finite center, for every $1 \le p < 2$ there exists a constant $K_p$ satisfying : $\|F * f\|_2 \le K_p \|F\|_p \|f\|_2$, for every $F \in L^p(G)$ and $f \in L^2(G)$.*

We will refer to any lcsc group satisfying the estimate of Theorem 12.1 a Kunze-Stein group.

12.1.2. *Monotonicity, Hölder continuity, and norm decay.* Suppose then that indeed for some $1 < r < 2$ we have $\|\mu_t\|_{L^r(G)} \le C \exp(-\theta t)$, $\theta > 0$ for a family $\mu_t$ of absolutely continuous measures on a semisimple Lie group $G$. It follows from Theorem 12.1 and Theorem 11.8 that in any probability-preserving action of $G$ with a spectral gap, the sum $\sum_{n=0}^{\infty} \|\mu_n\|_{L_0^2}^2$ is finite. It follows easily that the operator $\sup_{n \in \mathbb{N}} |\mu_n f(x)|$ satisfies the strong maximal inequality in $L^p$, $1 < p < \infty$. (In fact, it satisfies an exponential maximal inequality in $L_0^2(X)$, see §11.1). We would like to use the boundedness of the maximal function for the sequence $\mu_n$ in order to prove a strong $L^p$-maximal inequality for $\mu_t$, $t \in \mathbb{R}$. To that end, recall the use we made in §11.1 of the estimates for the volume growth of the balls, namely the fact that we could dominate the measure $\beta_t$ by $C\beta_{[t]+1}$, $C$ fixed. Thus it is natural to introduce the following conditions on general probability measures $\mu_t$ [MNS][N4], generalizing the conditions used already for the ball averages in §11.

**Definition 12.2.** A family $\mu_t$ of probability measures on an lcsc group $G$ is

   (1) Monotone, if $\mu_t \le C\mu_{[t]+1}$, as measures on $G$ (where $C$ is fixed, independent of $t > 0$).
   (2) Uniformly locally Hölder continuous, if for $F \in L^\infty(G)$ and $t > 0$
   $$|\mu_{t+s}(F) - \mu_t(F)| \le Cs^a \|F\|_\infty \ , \ 0 < s \le 1 \ .$$

We can generalize the arguments used in §11. and formulate a recipe for proving pointwise convergence with exponentially fast rate, using the monotonicity, local Hölder continuity and exponential decay of norms, as follows. (We refer to §11.1 for an example of the method, and [MNS][N4] and [N6] for more details on its use and applications).

**Recipe for pointwise ergodic theorems with exponentially fast rate of convergence.** To prove the strong maximal inequality for $\mu_t$, we follow the following steps.

(1) Let $f \in L^2(X)$ be a non-negative function on $X$. Then, since $\pi(\mu_t)f(x) \leq C\pi(\mu_{[t]+1})f(x)$, we have

$$\left\| \sup_{t>0} \pi(\mu_t)f(x) \right\|_{L^2(X)}^2 \leq C^2 \left\| \sup_{n\in\mathbb{N}} \pi(\mu_n)f(x) \right\|_{L^2(X)}^2 \leq C^2 \|f\|_{L^2(X)}^2$$

(2) The previous argument clearly extends to every $L^p$, $1 < p < \infty$, using the Riesz-Thorin interpolation theorem.

(3) Using the estimate $\|\pi(\mu_t)\|_{L_0^2(X)} \leq C \exp(-\theta t)$, the argument in (1) can be used to prove an exponential maximal inequality for the operators $\exp(\frac{1}{2}\theta t_k)\mu_{t_k}$ in $L_0^2$, where $t_k$ is a sequence such that the sum of the norms converges. Repeat the argument in $L_0^p(X)$.

(4) Now distribute $\exp(\frac{1}{4}\theta n)$ equally spaced points in the interval $[n, n+1]$. Then approximate $\pi(\mu_t)f$ by $\pi(\mu_{t_n})f$ using the closest point $t_n$ to $t$ in the sequence $t_k$. Estimate the difference using the exponential maximal inequality for the entire sequence $\mu_{t_k}$, and the local Hölder regularity of the family $\mu_t$, applied when $f$ is a bounded function.

(5) The previous argument gives an $(L^\infty, L^2)$-exponential maximal inequality, which says that the exponential-maximal function for $f$ bounded has an $L^2$ norm bound in term of the $L^\infty$-norm of $f$. Now interpolate against the usual strong maximal inequality in $L^p$ proved in step (1), using the analytic interpolation theorem.

Thus the recipe above establishes the following result.

**Theorem 12.3. Pointwise ergodic theorem with exponentially fast rate of convergence for general averages on semisimple Lie groups** [N4]. *Let $G$ be a connected semisimple Lie group with finite center (or any Kunze-Stein group). Let $f_t \in L^1(G)$ satisfy $f_t \geq 0$, and $\int_G f_t(g)dg = 1$. Assume that the family of probability measures $\mu_t$ with density $f_t$ form a monotone and uniformly locally Hölder continuous family. Assume that $\|f_t\|_{L^r(G)} \leq Ce^{-\theta t}$ for some $1 < r < 2$ and some $\theta > 0$. Let $(X, m)$ be a probability-preserving action of $G$, and assume that the unitary representation $\pi_0$ of $G$ on $L_0^2(X)$ satisfies $\pi_0^{\otimes n} \subset \infty \cdot \lambda_G$. Then $\pi(\mu_t)$ satisfies a pointwise ergodic theorem with exponentially fast rate of convergence to the ergodic mean in $L^p$, $1 < p < \infty$. In particular, for every $f \in L^p(X)$, and for almost every $x \in X$,*

$$\left| \pi(\mu_t)f(x) - \int_X f dm \right| \leq B_p(x, f) \exp\left(-\frac{\theta_p}{2n}t\right) \ , \ \theta_p > 0 \, .$$

*Furthermore an exponential $(L^p, L^r)$-maximal inequality holds in every $L^p$, $1 < p < \infty$.*

As we shall now demonstrate, Theorem 12.3 applies to some interesting geometric averages, as follows.

**12.2. horospherical averages.** Theorem 12.3 requires in its assumptions an estimate on the convolution norm $\|\lambda_G(\mu_t)\|$. Let us therefore note that the majorization principle due to C. Herz ([He2], see [Co3] for a discussion) has a corollary which is very useful in this regard, due to M. Cowling, U. Haagerup and R. Howe [CHH].

The corollary in question allows us to estimate the norm of the convolution operator $\lambda_G(f)$ on $L^2(G)$, by radialization, as follows.

**Theorem 12.4. Estimating convolution norms by radialization**[CHH]. *Let $G$ be a (non-compact) semisimple algebraic group over a locally compact non-discrete field. Then the following holds for every measurable function $F$*

$$\|\lambda_G(F)\| \leq \int_G \left( \int_K \int_K |F(kgk')|^2 \, dk \, dk' \right)^{1/2} \Xi(g) \, dg \ ,$$

*where $\Xi(g)$ is the Harish Chandra $\Xi$ function, namely the fundamental positive-definite positive spherical function on $G$.*

*Remark* 12.5.    (1) For $SL_2(\mathbb{C})$, $\Xi(a_t) = \varphi_0(a_t) = \frac{t}{\sinh t}$, and in particular, it decays exponentially in the distance $t$ in hyperbolic space.
(2) In general, the Harish Chandra $\Xi$-function on connected semisimple Lie groups has the same behaviour : it decays exponentially in the distance on the symmetric space $G/K$. More precisely $\Xi(g) \leq C \exp(-c\,|g|)$ $(c > 0)$, where $|g| = d(gK, K)$, $d$ the invariant distance on $G/K$ derived from the Riemannian structure given by the Killing form. The same holds true for general semisimple algebraic groups, where instead of the Riemannian distance on the symmetric space we consider the natural distance on the Bruhat-Tits building.

The estimate of Theorem 12.4 yields an estimate of the convolution norm $\lambda_G(f)$ of a non-radial function in terms of an associated radial function, which is much easier to control. In particular, together with Theorem 12.3 it gives a simple and explicit integral criterion for a family $f_t \in L^1(G)$, $t \in \mathbb{R}_+$ to satisfy a strong exponential-maximal inequality in every $L^p(X)$, $1 < p \leq \infty$.

A particularly interesting example to consider is that of horospherical averages, defined as follows.

Let $G = KAN$ be an Iwasawa decomposition of a connected semisimple Lie group with finite center. Recall from §10.4 that Haar measure $m_G$ can be normalized so that in horospherical coordinates it is given by (see [GV] or [He, Ch. I, Prop. 5.1])

$$\int_G f(g) dg = \int_{K \times A \times N} f(ke^H n) e^{2\rho(H)} dk \, dH \, dn$$

where $\rho$ is half the sum of the positive roots. Now let $h_t$ denote the absolutely continuous probability measure on $G$ whose density is given by $\chi_{U_t}/m_G(U_t)$, where $U_t = \{ke^H n \mid k \in K, \|H\| \leq t, n \in N_0\}$. Here $N_0$ is a fixed compact neighborhood of the identity in $N$. Appying Theorem 12.4 and Theorem 12.3, we obtain the following sample result.

**Theorem 12.6. Pointwise ergodic theorem with exponentially fast rate of convergence for horospherical averages in action with a spectral gap** [N4]. *Let $G$ be a connected semi-simple Lie group with finite center, $h_t$ the horospherical averages. Let $(X, m)$ be a probability-preserving action whose unitary representation $\pi_0$ in $L_0^2(X)$ has a spectral gap. Then*

(1)  $\|\lambda_G(h_t)\| \leq C \exp(-\theta t)$, $\theta > 0$.

(2) *If $\pi_0^{\otimes n} \subset \infty \cdot \lambda_G$, then for every $f \in L^p(X)$, $1 < p < \infty$ and almost every $x \in X$*

$$\left| \pi(h_t)f(x) - \int_X f\,dm \right| \leq B_p(x,f)\exp(-\frac{\theta_p}{2n}t)$$

*where $\theta_p > 0$.*

*Remark* 12.7.      (1) Again $h_t$ actually satisfies a more precise estimate, namely an exponential $(L^p, L^r)$-maximal inequality (see [N4] for details).

(2) Theorem 12.6 can be established also for other averages defined by horospherical coordinates. As an example, consider the case of a real-rank-one group. Then the group $A = \{a_t \,;\, t \in \mathbb{R}\}$ figuring in the Iwasawa decomposition is one dimensional, and let $I_1$ denote the unit interval in $\mathbb{R}$. Then $t + I_1$ is a unit interval with center $t$, and let $J_t = \{ka_sn \mid k \in K, s \in t + I_1, n \in N_0\}$. Then the Haar-uniform averages $j_t$ on $J_t$ satisfy the conclusion Theorem 12.6.

(3) We note that it is interesting to compare the horospherical averages $h_t$ and $j_t$ to the averages constructed in the proof of Theorem 7.13. Note that here the pointwise ergodic theorem is strengthened to yield an explicit exponentially fast rate of convergence to the ergodic mean.

(4) Another interesting comparison is between $j_t$ and the averages $m_K * \delta_{a_t}$ appearing in the mean ergodic theorem of W. Veech [V].

(5) Similar results can of course be established for semisimple algebraic groups.

## 12.3. Averages on discrete subgroups.

We have commented throughout our discussion on the problem of comparing between the ball averages on an lcsc group, and the discrete ball averages on its lattice subgroups. In the case homogeneous of nilpotent groups such as of $\mathbb{Z}^d \subset \mathbb{R}^d$ and $H_n(\mathbb{Z}) \subset H_n$, the comparison was completely straightforward : proofs of maximal inequalities on the Lie group can be easily adapted to prove maximal inequalities on the discrete lattice, and then the transfer principle implies that they hold in *every* measure-preserving action (see chapter 5 for the details).

In §10.5 and §11.3 we have encountered a select group of extra symmetric lattices in some simple algebraic groups over local fields, and other totally disconnected Gelfand pairs $(G, K)$. In such a lattice $\Gamma$, the group algebra $\ell^1(\Gamma)$ contains an isomorphic copy of the commutative convolution algebra $M(G, K)$, and most of the basic problems in spectral theory and ergodic theory pertaining to these averages on $\Gamma$ are reducible to the corresponding problems for $M(G, K)$. This gives a satisfactory solution to the basic questions in spectral and ergodic theory for the corresponding averages on $\Gamma$, provided that the representation theory of the commutative algebra $M(G, K)$ is sufficiently well understood.

However, in the set-up of general lattices $\Gamma$ in non-amenable Gelfand pairs, and even semisimple algebraic groups, there is usually no direct connection between $M(G, K)$ and $\ell^1(\Gamma)$, and the basic problems in ergodic theory cannot be resolved using this method. Note that the absence of a transfer principle implies that a result on a maximal inequality for convolutions on the discrete lattice has no bearing on the case of a general action. Furthermore, the natural discrete ball averages (w.r.t. a word metric) on the lattice do not commute in general, and so there is no natural commutative algebra whose spectral theory can be used to establish even a mean ergodic theorem for the averages. Thus establishing a pointwise, or even mean, ergodic theorem for the discrete uniform ball averages in *arbitrary* measure preserving action of $\Gamma$ is a challenging goal, outside the short (but interesting) list of examples in §10.5 and 11.3. Nevertheless, it is possible to make considerable progress on this problem, at least for certain natural averages on $\Gamma$, although these are usually not comparable to balls w.r.t. a word metric. These very recent results are based on three principle, namely induction of actions from the lattice $\Gamma$ to the group $G$, a duality principle which controls discrete averages on the lattice $\Gamma$ by certain (non-radial) absolutely continuous averages on the group $G$, and ergodic theorems for sufficiently general non-radial averages on on the group $G$. These results will be reported in [GN].

For a general lattice of a semisimple Lie group, we will content ourselves here with the following partial result, which applies the spectral transfer principle, and establishes a pointwise ergdic theorem for the lattice in actions of the group. It exemplifies the natural procedure of estimating discrete convolution operators in terms of absolutely continuous ones, by reducing the problem to geometric comparison for the translation action on $G$, and thus obtaining an estimate of convolution norm.

More precisely, let $\beta_k = \frac{1}{vol B_k} \chi_{B_k}(g)$ (where $B_k$ is the bi-$K$-invariant lift of a ball of radius $k$ in $G/K$). Let $\Gamma$ be a lattice subgroup of $G$, let $B_k(\Gamma) = B_k \cap \Gamma$, and let $b_k = \frac{1}{|B_k(\Gamma)|} \sum_{\gamma \in B_k(\Gamma)} \gamma$. We have :

**Theorem 12.8. Pointwise ergodic theorems for lattice averages in $G$-actions with a spectral gap** [N4]. *Let $\Gamma \subset G$ be a lattice in a connected semisimple Lie group $G$ with finite center. Then the sequence $b_k \in \ell^1(\Gamma)$ satisfies :*

(1) $\|\lambda_\Gamma(b_k)\| \leq C \exp(-\theta k)$. *Here*

$$0 < \theta < \theta_\beta(G) = \lim_{t \to \infty} -\frac{1}{t} \log \|\lambda_G(\beta_t)\| \, .$$

(2) *In any $\Gamma$-action satisfying $\pi_0^{\otimes n} \subset \infty \cdot \lambda_\Gamma$, and in particular in any action of $G$ with a spectral gap, for any $f \in L^p(X)$, $1 < p < \infty$ and for almost every $x \in X$*

$$\left| b_k f(x) - \int_X f dm \right| \leq C(x, f) \exp(-\frac{\theta}{2n} k) \, .$$

Clearly, the spectral estimate for convolutions stated in (1), together with norm estimate provided by the spectral transfer principle (see Theorem 11.8) implies the exponentially fast pointwise convergence in (2). Such spectral estimates can be established for many other discrete groups, and not only for lattices. In fact, the unitary representation of $\Gamma$ in $L^2(G)$ is equivalent with countably many copies of the unitary representation of $\Gamma$ in $\ell^2(\Gamma)$. Thus any sequence of measures $\nu_n \in \ell^1(\Gamma)$ with $\|\lambda_\Gamma(\nu_n)\|_{\ell^2(\Gamma)} \leq \exp(-\theta n)$ will satisfy the conclusion of Theorem 12.8.

The challenge of establishing spectral norm estimates for convolution operators on discrete groups has attracted considerable attention, and one particularly interesting approach was to establish the even stronger property of rapid decay [Ha][J], one of whose formulations is as follows. Assume $d$ is a word metric, and let $L(\gamma) = d(e, \gamma)$. Then for some $s = s(\Gamma, d) \geq 0$, some $C = C(\Gamma, d) > 0$, and every finitely supported function $f \in \ell^1(\Gamma)$,

$$\|\lambda_\Gamma(f)\| \leq C \|f \cdot (1 + L)^s\|_{\ell^2(\Gamma)} \, .$$

From this estimate it follows that if in addition the spheres (or balls) have strict $t^q \exp ct$-exponential growth (with $c > 0$) then the spectral norm decays exponentially fast. We note that strict $t^q \exp ct$-volume growth has been established in the case of word metrics on word-hyperbolic groups in [Coo], so that every discrete hyperbolic subgroup of a connected semisimple Lie group satisfies the conclusions on Theorem 12.8, w.r.t. balls defined by a word metric on it (see [N4, Cor. 3.3]).

The rapid decay property has been established in a number of interesting cases, including all word-hyperbolic groups [J], uniform lattices in $SL_3(\mathbb{R})$ and $SL_3(\mathbb{C})$ [Laf], certain lattices acting on rank-two Bruhat-Tits building [RRS], and groups acting properly on cube complexes [CR]. However strict $t^q \exp ct$-volume growth has not been established in general for uniform lattices or cube complex groups, and constitutes a completely open problem for general discrete (sub)groups.

## 13. Weighted averages on discrete groups and Markov operators

In this chapter we will consider briefly some applications of the general theory of Markov operators to ergodic theorems for group actions. We divide the chapter into three parts, dealing with applications of ergodic and maximal theorems for the following sequences :

(1) uniform averages of powers of a single Markov operator,
(2) subadditive sequences of self-adjoint Markov operators,
(3) powers of a single self-adjoint Markov operator.

13.1. **Uniform averages of powers of a Markov operator.** In our discussion so far we have put a great deal of emphasis on proving ergodic theorems for Haar-uniform averages on geometrically significant sets. A different approach, which has long roots in ergodic theory and the theory of Markov processes is to relax the requirement of unifom averages and allow weighted ergodic theorems. Thus we can consider for a discrete group a sequences of averages of the form $\sum_{\gamma \in \Gamma} \nu(\gamma) \delta_\gamma$, where $\nu$ is a general probability measure on $\Gamma$. One ergodic theorem that can be obtained here is for the sequence of averages $\nu_n = \sum_{k=0}^n \nu^{*n}$, namely the uniform average of convolution powers of a single measure. It follows immediately from the general theory of non-negative contractions developed by Hopf and Dunford-Schwartz (see e.g. [DSI] for an extensive exposition) that the sequence $\pi(\nu_n)$ satisfies the weak-type $(1,1)$ maximal inequality in $L^1$, and converges pointwise to a $\pi(\nu)$-invariant function, in every probability-preserving action of $\Gamma$. In particular, if the support of $\nu$ generates $\gamma$ as a group, then the limit is a $\Gamma$-invariant function, namley $\int_X f \, dm$ when $\Gamma$ acts ergodically.

The ergodic theorem for the uniform averages of powers of a Markov operator can be used to obtain weighted ergodic theorems for actions of several transformations also in another manner, namely using the construction of skew product actions. This useful idea has been introduced already by S. Kakutani in his proof of the random ergodic theorem [Ka] for almost all sequences of transformations chosen independently (or according to a Markov measure) from a finite set (say). Thus for a probability-preserving action of $FS_r$, the free semigroup on $r$ non-commuting elements, one can form the skew product $\Omega_r \times X$, where $\Omega_r$ is the topological Markov chain cosisting of infinite words in the free generators of the free semigroup. Consider the transformation $T(\omega, x) = (S\omega, \omega_1 x)$ where $S : \Omega_r \to \Omega_r$ is the forward shift, and $\omega_1$ the first letter of $\omega$. One takes the natural probability measure $p$ on $\Omega_r$ associated with uniform weights on the generators $S$ (or any other Markov measure on $\Omega_r$), and the measure $p \times m$ on $\Omega_r \times X$. Ergodicity of $T$ and the $FS_r$-action on $X$ are equivalent [Ka], and thus the uniform averages of powers of the operator $T$ satisfy the pointwise ergodic theorem, by the Hopf-Dunford-Schwartz theorem, so that

$$\lim_{n \to \infty} \frac{1}{n+1} \sum_{k=0}^{n} F(S^k \omega, \omega_k \omega_{k-1} \cdots \omega_1 x) = \int_{\Omega_r \times X} F(\omega, x) dp \, dm$$

Now apply the foregoing result to the functions $F(\omega, x) = f(x)$ on the skew product $\Omega_r \times X$ which are lifted from $X$, and take expectations w.r.t. the probability $p$ on $\Omega_r$. Then one obtains a pointwise ergodic theorem in $L^1(X)$ for the uniform averages $\mu_n = \frac{1}{n} \sum_{k=0}^{n} \sigma_k$ of the sphere averages on the free semigroup (or more general weights determined by an arbitrary Markov measure on $\Omega_r$).

Note however that the passage to expected values implies that we can not conclude that the weak-type $(1, 1)$-maximal inequality in valid for $f_\mu^*$.

This application of Kakutani's random ergodic theorem in the context of ergodic theorems for free groups is due to R. Grigorchuk [Gri1] [Gri2]. The idea of using skew products and the theory of Markov operators is developed further in [Gri2] and [Bu1][Bu2]. Indeed, the method is suitable for operator averages for Banach space representations which do not arise from measure-preserving action [Gri1][Bu1], requires the measure $m$ on $X$ to be merely stationary, and not necessarily invariant [Gri2], and applies also when the weights taken along the group orbit depend on the starting point [Gri2]. However, in general these results establish weighted ergodic theorems, which bear no discernible relation to the uniform averages of spheres w.r.t. a word metric. The only exception is in the case where the Markov measure is associated with all the generators in $S$ having equal probability. We remark that the same analysis applies also to the free group, and not only the free semi-group.

We note that the maximal and pointwise ergodic theorem in $L^2$ for the uniform averages of the spheres averages was proved in [N1], and for $f \in L^1$ in [NS1]. The proof in [NS1] is also probabilistic and quite elementary, and does not require spectral theory. It uses the standard estimates of the central limit theorem for convolution powers of a binomial distribution on $\mathbb{N}$, and Hopf's maximal inequality (see a particularly simple proof by Garcia of the latter in [Ga]). Indeed, it is shown that the sum of convolution powers of $\sigma_1$ dominates the uniform average of spheres, or more precisely :

$$\mu_n = \frac{1}{n} \sum_{k=0}^{n} \sigma_k \leq \frac{C}{3n+1} \sum_{k=0}^{3n} \sigma_1^{*k} .$$

Thus it follows immediately that the maximal function $f_\mu^*$ satisfies the weak-type $(1, 1)$-maximal inequality, by the maximal inequality for the average of powers. Furthermore, it is clear that on functions of the form $h = f - \pi(\sigma_1)f$ the sequence $\pi(\mu_n)h$ converges pointwise, and this implies the pointwise ergodic theorem as usual by the recipe of §2.4.

We record the facts described above as follows :

**Theorem 13.1.** *The sequence $\mu_n = \frac{1}{n+1} \sum_{k=0}^{n} \sigma_k$ of uniform averages of spheres on the free group satisfied the weak-type $(1, 1)$-maximal inequality in $L^1$ and is a pointwise ergodic sequence in $L^p$, for all $1 \leq p < \infty$.*

Thus it appears that the sequence of uniform averages of the sphere averages $\mu_n = \frac{1}{n} \sum_{k=0}^{n} \sigma_k$ is a natural sequence of weighted averages to consider, and we may inquire for which discrete groups is Theorem 13.1 satisfied for a word metric. We next turn to consider a method which establishes at least the maximal inequality in $L^2$ for $\mu_n$, for a certain class of groups. This method was first employed in [N1] to prove the maximal and pointwise ergodic theorem in $L^2$ for the averages $\mu_n$ on $\mathbb{F}_k$. It also does not use spectral considerations, and is based on a general subadditive maximal inequality, which we consider below.

## 13.2. Subadditive sequences of Markov operators, and maximal inequalities on hyperbolic groups. Let us introduce the following definition.

**Definition 13.2. Subadditive sequences** (see [N1]). A sequence $T_n$ of operators on $L^2(X)$ will be called a *subadditive sequence of self adjoint Markov operators* if it satisfies the following:

(1) $T_n = T_n^*$, $\|T_n\| \leq 1$ .
(2) $T_n f \geq 0$ if $f \geq 0$, $T_n 1 = 1$.
(3) There exist a constant $C_0 > 0$, a positive integer $k$, and a fixed non-negative bounded operator $B$ on $L^2(X)$ such that:

$$T_n T_m f(x) \leq C_0(T_{kn} f(x) + T_{km} f(x)) + B f(x)$$

for all bounded and nonnegative $f \in L^2$.

We can now state the following subadditive maximal inequality, proved independently in [Bar] and [N1].

**Theorem 13.3. Subadditive maximal inequality** [Bar][N1]. *Let $T_n$ be a subadditive sequence of self adjoint Markov operators. Define $f^*(x) = \sup_{n \geq 0} |T_n f(x)|$. Then $\|f^*\|_2 \leq C \|f\|_2$ for all $f \in L^2$. We can take $C = 2C_0 + \|B\|$.*

We note that the subadditive maximal inequality of Theorem 1 generalizes similar results due to E. M. Stein [S0] [S1] and B. Weiss [We]. In particular, it was applied in [S0][S1] to prove a pointwise ergodic theorem for the even powers of positivity-preserving self-adjoint contractions on $L^2$. Also, it is noted in these references that it implies the pointwise convergence of martingales in $L^2$, as well as Birkhoff's pointwise ergodic theorem in $L^2$. The origin of this maximal inequality is attributed in [S1] [We] to A. Kolmogoroff and G. Seliverstoff [KS], and to R. E. A. C. Paley [Pal].

It is reasonable to expect that for a large class of discrete groups $\Gamma$, the sequence of averages $\mu_n \in \ell^1(\Gamma)$ will satisfy the subadditive inequality given by $\mu_n * \mu_m \leq C(\mu_{kn} + \mu_{km}) + b$. When this inequality holds in $\ell^1(\Gamma)$, then $T_n = \pi(\mu_n)$ is a subadditive sequence of self adjoint Markov operators on $L^2$ in any finite-measure-preserving action, and therefore will satisfy the strong maximal inequality in $L^2$.

This is indeed the case at least in the following context.

**Theorem 13.4. Subadditive inequality for word-hyperbolic groups**[FN]. *Let $\Gamma$ be a non-elementary word-hyperbolic group, $S$ a finite symmetric generating set. Then there exist constants $1 < q < \infty$ and $0 < C < \infty$, depending only on $(\Gamma, S)$, such that the following inequalities hold:*

(1) $\sigma_t * \sigma_s \leq C \sum_{j=0}^{2s} q^{-(s - \frac{1}{2}j)} \sigma_{t-s+j}$ *if $t \geq s$.*
(2) $\mu_n * \mu_m \leq C(\mu_{2n} + \mu_{2m})$.

Combining the foregoing results, we have

**Theorem 13.5. Maximal inequality for word-hyperbolic groups**[FN]. *Let* $\Gamma$ *be a word-hyperbolic group,* $S$ *a finite symmetric generating set. Then the sequence* $\mu_n$ *satisfies the strong maximal inequality in* $L^2(X)$, *i.e.* $\|f_\mu^*\|_2 \leq C(\Gamma, S) \|f\|_2$ *for every* $f \in L^2(X)$.

We remark that an elementary word-hyperbolic group is a finite extension of $\mathbb{Z}$, and the subadditivity of $\mu_n$ in this case (for free generators) is of course easily verified (see [S1] for the case of $\mathbb{N}$).

As usual, given the maximal inequality for $\mu_n$ in $L^2$, to complete the proof of the pointwise ergodic theorem it suffices (see §2.4) to find a dense set of functions $f \in L^2$ where $\pi(\mu_n)f(x) \to \int_X f\, dm$, almost everywhere and in the $L^2$-norm. For hyperbolic groups, a sufficient condition is the existence of a dense set of functions $f \in L^2$ satisfying the following exponential mixing condition: $|\langle \pi(\gamma)f, f\rangle| \leq C_f \exp(-c_f |\gamma|)$, for some $c_f > 0$ (see [FN]). However, in general these issues are far from being resolved, and we thus formulate the following..

**Problem 13.6. Analogs of von-Neumann and Birkhoff theorems for word-hyperbolic groups.**

(1) Is the mean ergodic theorem valid for the averages $\mu_n$ on every word-hyperbolic group ? Is the pointwise ergodic theorem valid for $\mu_n$ in $L^2$, or even $L^1$ ?

(2) Is there any finitely generated group for which the subadditive convolution inequality fails for the averages $\mu_n$ ?

**13.3. The powers of a self-adjoint Markov operators.** For a *self adjoint* Markov operator a result much sharper than the Hopf-Dunford-Schwartz ergodic theorem was subsequently proved independently by E. Stein [N0] and J. C. Rota [Rot], using two entirely different methods. In both cases, the results proved imply as a special case that when the measure $\nu$ is symmetric, then not only does the sequence of uniform averages $\pi(\mu_n) = 1/(n+1)\sum_{k=0}^{n} \pi(\nu_k)$ converge, but already the sequence of powers $\pi(\nu^{*2n})$ converge, pointwise almost everywhere, and in the $L^p$-norm, at least for $1 < p < \infty$. The limit is a $\pi(\nu^{*2})$-invariant function, and thus invariant under the group generated by the support of $\nu^{*2}$. We note that the passage to even powers reflect the same periodicity phenomenon that was encountered in §10.5 with regard to the free groups. Namely it reflects the fact that $\pi(\nu)$ may have an eigenfunction $f_0$ with eigenvalue $(-1)$, or equivalently that the commutative convolution algebra generated by the powers of $\nu$ may have the an eigenfunction $f_0 \in L^2(X)$ realizing the sign character of the algebra $\nu^{*n} \mapsto (-1)^n$.

Both results mentioned above apply to general self-adjoint Markov operators $P$ and not only to convolution powers, but in this generality, it was established by D. Ornstein [Or] that pointwise convergence usually fails to hold in $L^1$ for the powers of $P^2$.

Stein's method is spectral, and can be viewed as a special case of the general method described in §10.2, where the commutative algebra is taken to be simply the convolution algebra $\ell^1(\mathbb{N})$, corresponding to the algebra genrated by powers of a single operator. The continuous *-representation are given here simply by self-adjoint contraction operators, and the *-characters are given by $\varphi_\lambda(k) = \lambda^k$, where $\lambda \in [-1, 1]$. This point of view was developed in [N1], where Stein's method was generalized to apply e.g. to the radial algebra on free groups.

Rota's method depends on probabilistic considerations, and reduces the point-wise convergence theorem for the even powers of a self-adjoint Markov operator to the pointwise convergence theorem for martingales, on an auxiliary probability space constructed from the Markov operator in question. In was shown in [Rot] that in fact pointwise convergence holds for $f \in L(\log L)(X)$.

In [Bu3] A. Bufetov has constructed, for a given action of the free group $\mathbb{F}_r$ on a space $(X, m)$, and a given free generating set $S$, a Markov operator $P$ on the set $X \times S$, which is not self-adjoint, but such that the operators $(P^*)^n P^n$ are comparable to the action of the sphere averages on $\mathbb{F}_r$ on the space $X$. Rota's theorem applies to such sequence and the pointwise convergence of $\sigma_{2n} f$ to a function invariant under words of even length follows, in $L(\log L)(X)$.

## 14. FURTHER DEVELOPMENTS

14.1. **Some non-Euclidean phenomena in higher-rank groups.** In the present section we return to the set-up of connected semisimple Lie groups, and we would like to demonstrate the fact that exponential volume growth on semisimple Lie groups (and more general lcsc groups) can be used to derive a number of non-standard maximal inequalities and pointwise convergence results which have no Euclidean analog. This phenomenon occurs for actions of semisimple algebraic groups of split rank at least two, or for direct products of lcsc groups. To illustrate some of these results, let us consider the simplest case where $G = L_1 \times L_2$ is a product of two real-rank one groups, say $L_1 = L_2 \cong PSL_2(\mathbb{R})$ for definiteness. We attempt to prove maximal inequalities for the product group $G$ based on maximal inequalities of its factors $L_i$. The ball averages on the factor groups $L_i$ is given in this case by (choosing the curvature on $\mathbb{H}^2$ appropriately)

$$\beta_t = \frac{\int_0^t \sinh s \ \sigma_s ds}{\int_0^t \sinh s \ ds}$$

Let us define $D_t = \{(u, v) \, ; \, u^2 + v^2 \leq t^2\}$, and let $\gamma_t^{(i)}$ (resp. $\sigma_t^{(i)}$) be the unit shell (resp. sphere) averages on the real-rank one group $L_i$. Then the ball averages $\beta_t^G$ on the product group $G = L_1 \cdot L_2$ satisfy the following estimate :

$$\beta_t^G = \frac{\int_{D_t} |\sinh u \sinh v| \, \sigma_{|u|}^{(1)} \sigma_{|v|}^{(2)} \, du dv}{\int_{D_t} |\sinh u \sinh v| \, du dv} \leq B \frac{\sum_{D_t \cap \mathbb{N}^2} \exp n \exp m \ \gamma_n^{(1)} \gamma_m^{(2)}}{\sum_{D_t \cap \mathbb{N}^2} \exp n \exp m}$$

On the l.h.s. we have the integral w.r.t. the $G$-invariant Riemannian volume, which in geodesic polar coordinates is given by the integral w.r.t. the density $|\sinh u \sinh v|$ on a disc of radius $t$ in $\mathbb{R}^2$. Using the Weyl group symmetry, this integral is equal to four times its value on $D_t \cap \mathbb{R}_+^2$. On the r.h.s. the integral was further replaced by the weighted sum of integrals (w.r.t. the uniform Euclidean measure) on unit squares whose lower left hand corner is a lattice point of norm at most $t$. The weight attached to each unit square is the obvious one : we estimate the function $\sinh u \sinh v$ on the square whose lower left hand corner is $(n, m)$ by $C \exp n \exp m$. Since the shell averages on $L_i$ satisfy the strong maximal inequality in $L^p$, $1 < p \leq \infty$ by Theorem 9.6, it follows that the same holds for the operators $\beta_t^G = \beta^{L_1 \times L_2}$.

The support of $\gamma_t^{(1)} \gamma_s^{(2)}$, $(t, s) \in \mathbb{R}_+^2$ has radial coordinates in $\mathfrak{a}_+ \cong \mathbb{R}_+^2$ which constitute (the $W$-orbit of) a square of unit side length, whose lower left hand corner is the element $(tH_1, sH_2)$. But now note that already the family of unit square averages $\gamma_t^{(1)} \gamma_s^{(2)}$ (for arbitrary non-negative $s$ and $t$) satisfy a strong maximal inequality in $L^p$, $p > 1$. Indeed, each square average as above is the product of two one-dimensional shell (or "interval", in this context) averages on real rank one groups, and Theorem 9.5 applies to the shell averages on $L_1$ and $L_2$. As a result, the strong maximal inequalities which hold for the unit square averages in an arbitrary measure-preserving action of $G$, also hold for an *arbitrary family of sets which admit a reasonable covering by unit squares*, not only for the discs described above (which correspond to the radial coordinates of ball averages w.r.t. the Riemmanian Killing metric).

It is possible to greatly expand the scope and generality of results of this kind, and apply them to all higher-rank semisimple groups. In particular they yield a pointwise ergodic theorem for uniform averages supported on an arbitrary sequence of bi-$K$-invariant sets which are reasonably covered by "cube averages" of a fixed size, provided they leave eventually any given compact set. In particular these arguments apply to a wide array of of bi-$K$-invariant averages which occupy an exponentially decaying fraction of the volume of a ball. This improves the results we previously described for the shells, which occupy a *fixed* proportion of the volume of the ball. To illustrate thes points, let us define the following simple geometric property of sets in the vector space $\mathfrak{a}$.

**Definition 14.1. Condition** $A_{(c,C)}$. A measurable set $E \subset \mathfrak{a}$ (of finite measure) satisfies condition $A_{(c,C)}$ if for every $H \in E$, there exists $H' \in E$, satisfying the conditions $\|H - H'\| \le C$, and $b_c(H') \subset E$. Here $b_c(H')$ is a ball of radius $c$ and center $H'$ in $\mathfrak{a}$.

Clearly, any union of sets satisfying condition $A_{(c,C)}$ also satisfies it. Equally clearly, all balls of radius at least a fixed constant satisfy condition $A_{(c,C)}$ for some $(c, C)$. The same holds for Euclidean cells (=cubes) of fixed side length. Thus any union of such sets (with fixed $(c, C)$) also satisfies the condition.

Now view $\mathfrak{a}$ as the Lie algebra of a split Cartan subgroup of a connected semisimple Lie group with finite center. For a Weyl group invariant set $E \subset \mathfrak{a}$ let $K \exp(E)K = R(E)$ be the radialization of $E$, and let $\nu_E$ be the normalized average on $R(E)$, namely the normalized restriction of Haar measure to this set. We can now formulate the following

**Theorem 14.2. Pointwise ergodic theorem for general sequences of bi-$K$-invariant averages**[N6]. *Let $G$ be a connected semisimple Lie group with finite center. Then*

(1) *The maximal function*

$$\mathcal{A}^* f(x) = \sup \left\{ |\pi(\nu_E) f(x)| \; ; \; E \subset \mathfrak{a} \text{ , and } E \text{ satisfies condition } A_{(c,C)} \right\}$$

*satisfies the strong maximal inequality in every $L^p$, $1 < p < \infty$.*

(2) *Let $E_n$ be any sequence of measurable sets satisfying condition $A_{(c,C)}$. Assume that $\lim_{n\to\infty} vol_G (E_n \cap b_t) / vol_G (E_n) = 0$, for every $t > 0$, where $b_t$ is a ball of radius $t$ in $\mathfrak{a}$ with center 0. Assuming $G$ simple, the averages $\nu_{E_n}$ satisfy the pointwise ergodic theorem in $L^p$, $1 < p < \infty$. For $G$ semisimple, the same conclusion holds, provided we assume the previous condition also for the projection of $E_n$ to any factor group.*

The phenomenon described in Theorem 14.2 allows for a choice of very general sequences, and does not seem to have any Euclidean analog We refer to [N6] for further details.

## 14.2. Best possible rate of convergence in the pointwise theorem. The pointwise ergodic theorems with exponentially fast speed of convergence for ball averages can be considerably sharpened, and the best possible rate of convergence can be determined. Furthermore, exponentially fast pointwise convergence can be proved even for the sphere averages on real-rank one groups (in the range of $p$ where pointwise convergence actually holds for $L^p$ functions !). This requires spectral arguments which are considerably more elaborate than those we have presented here, and the details can be found in [N8].

## 14.3. Added in proof. We note the following very recent developments that have taken place since the final version of the present survey was submitted.

### 14.3.1. *Exact polynomial volume growth.* E. Breuillard has completed the proof of the following remarkable result.

**Theorem 14.3.** *On every lcsc group of polynomial volume growth, the balls w.r.t. any word metric have exact polynomial volume growth. In fact, the same holds true for every invariant asymptotically geodesic pseudo-metric on the group.*

This result provides an alternative proof of the general case of the localization conjecture, namely the fact that the balls are asymptotically invariant under translations, from which the pointwise ergodic theorem in $L^1$ follows - see the discussion in Sections 4 and 5. We note that in the case of connected Lie groups, precise results are obtained by Breuillard regarding the actual geometric shape of the balls, and not just their volume asymptotics. The passage from the Lie case to the general case uses the results of Guivarc'h [Gu] and Jenkins [J] on growth, and of Gromov [G] and Losert [Lo] on the structure of groups of polynomial volume growth. Thus the theorem considerably sharpens Proposition 5.11 of Section 5.5, where we have used these results to deduce strict polynomial growth. Details can be found in the preprint "The asymptotic shape of metric balls in groups of polynomial growth" by E. Breuillard.

### 14.3.2. *Exact exponential volume growth.* Given a linear representation $\tau : G \to GL_n(\mathbb{R})$ of a semisimple group $G$, and a vector space norm on $M_n(\mathbb{R})$, one can consider (say) the distance function $\log \|g\|$. It was noted in Theorem 4.7 of §4.4 that in [GW] it was shown that the balls of radius $t$ w.r.t. such a distance function have exact $t^q \exp(ct)$ growth. F. Maucourant has developed an alternative approach to this result which yields that in fact after scaling by the volume growth one obtains $w^*$ convergence to a limiting Radon measure on $M_n(\mathbb{R})$. Details can be found in the preprint "Homogeneous asymptotic limits of Haar measure of semisimple linear groups and their lattices", by F. Maucourant.

### 14.3.3. *The ergodic theory of lattice subgroups.* It is possible to develop general pointwise ergodic theorems for *arbitrary actions* of a lattice $\Gamma$ in a semisimple group $G$. Such theorems apply to the discrete averages on lattice points in balls w.r.t. a distance functions on the group $G$. We note here a sample result formulated in the simplest case, and refer to [GN] for a full discussion.

Let $\beta_k = \frac{1}{vol B_k} \chi_{B_k}(g)$, where $B_k$ is the bi-$K$-invariant lift of a ball of radius $k$ in $G/K$. Let $\Gamma$ be a lattice subgroup of $G$, let $B_k(\Gamma) = B_k \cap \Gamma$, and let $b_k = \frac{1}{|B_k(\Gamma)|} \sum_{\gamma \in B_k(\Gamma)} \gamma$. We have :

**Theorem 14.4. Pointwise ergodic theorem for general lattice actions**[GN].
*Let $\Gamma \subset G$ be a lattice in a connected semisimple Lie group $G$ with finite center. Then the sequence $b_k \in \ell^1(\Gamma)$ satisfies, in every ergodic probability preserving action of $\Gamma$, for any $f \in L^p(X)$, $1 < p < \infty$ and for almost every $x \in X$ :*

(1) *$\lim_{k \to \infty} b_k f(x) = \int_X f \, dm$ .*

(2) *In any ergodic action of $\Gamma$ with a spectral gap (and thus in every ergodic action if $\Gamma$ has property $T$), there exists $\delta_p = \delta_p(X) > 0$, such that for any $f \in L^p(X)$, $1 < p < \infty$ and for almost every $x \in X$*

$$\left| b_k f(x) - \int_X f \, dm \right| \le C_p(x,f) \exp(-\delta_p k) \, .$$

*If $\Gamma$ has property $T$, then $\delta_p$ depends only on $\Gamma$, and not on $X$.*

In reference to Theorem 14.4, note that in Theorem 12.8 a similar conclusion is asserted, but only for those actions of $\Gamma$ arising from actions of $G$. However the rate of exponential convergence obtained in Theorem 12.8 is in general faster. Also note that in Theorem 11.10 the best possible rate of exponential convergence is obtained for balls w.r.t. a *word metric* on the lattice, whereas $b_k$ are not associated with a word metric. However in Theorem 11.10 the lattices (and the metrics) are severely restricted.

We note that the same result holds for balls (and many of their subsets) w.r.t. a large class of distance functions on the group, and allows the calculation of the main term as well as an estimate of an error term in many lattice point counting problems. We refer to [GN] for details.

14.3.4. *The ball averaging problem for word metrics on semisimple groups.* As noted in Theorem 4.11, word metrics on semisimple groups are coarsely isometric to norm-like metrics. It is possible to utilize this fact and obtain a solution of the ball averaging problem in this context. We formulate for simplicity the following basic special case.

**Theorem 14.5. Pointwise ergodic theorem for word metric balls on simple Lie groups**. *The balls defined by any word metric on a connected simple Lie group with finite center satisfy the pointwise ergodic theorem in $L^p$, $1 < p < \infty$.*

In fact, a similar result holds for all algebraically connected semisimple algebraic groups, at least in actions where $G^+$ acts ergodically. The convergence is exponentially fast almost surely, if the action has a spectral gap. Details will appear in the paper "On the ball averaging problem in ergodic theory", currently under preparation.

14.3.5. *Further reading.* The present survey aimed for the most part to indicate just the bare outlines of the relevant arguments appearing in the proofs of the ergodic theorems cited, and most details are of course left out. The reader wishing to learn more about some of these arguments is refered to the detailed comprehensive exposition in the forthcoming book "Théorèmes Ergodiques pour les Actions de Groupes". This book is the result of a collaborative effort initiated by the late Martine Babillot, with the participation of C. Anantharaman, J-Ph. Anker, A.

Batakis, A. Bonami, B. Demange, F. Havard, S. Grellier, Ph. Jaming, E. Lesigne. P. Maheux, J-P. Otal, B. Schapira and J-P. Schreiber. The book contains a wealth of information on ergodic theorems for group actions, including an elaboration of a number of the topics mentioned in the survey.

DEPARTMENT OF MATHEMATICS, TECHNION IIT

*Current address*: Department of Mathematics, Technion, Haifa, Israel

*E-mail address*: `anevo@tx.technion.ac.il`